\documentclass[pagebackref,coverfontpercent=100,british,biblatex,biblatexstyle=alphabetic,final,online]{adsphd}

\usepackage[dutch,english]{babel}

\usepackage{nomencl}   
\usepackage[style=long,number=none]{glossary} 

\title{Motivic Integration and Logarithmic Geometry}

\author{Emmanuel}{Bultot}

\supervisor{Prof.~dr.~J.~Nicaise}{}
\president{Prof.~dr.~W.~Veys}
\jury{Prof.~dr.~N.~Budur\\
      Prof.~dr.~A.~Kuijlaars
}
\externaljurymember{Dr.~M.~Raibaut}{Université de Savoie}
\externaljurymember{Prof.~dr.~J.~Sebag}{Université Rennes 1}

\phddegree{Science} 
\faculty{Faculty of Science}
\department{Department of Mathematics}
\researchgroup{Algebra section}
\address{Celestijnenlaan 200B}
\email{emmanuel.bultot@wis.kuleuven.be} 
\website{} 

\date{April 2015}
\copyyear{2015}

\renewcommand{\nomname}{List of Symbols}
\makeatletter
\let\@printnomenclatureorig\@printnomenclature
\def\@printnomenclature[#1]{%
  \cleardoublepage%
  \chaptermark{\nomname}
  \@printnomenclatureorig[#1]
}
\makeatother
\makenomenclature

\newcommand{\glossname}{Abbreviations}
\makeglossary

\let\printglossaryorig\printglossary
\renewcommand{\printglossary}{%
  \renewcommand{\glossaryname}{\glossname}
  \cleardoublepage%
  \printglossaryorig\chaptermark{\glossname}}

\AtEveryBibitem{
	\clearfield{doi}
	\clearfield{issn}
	\clearfield{isbn}
	
	\ifentrytype{misc}{}{
		\clearfield{url}
	}
}

\DeclareFieldFormat{postnote}{#1} 

\usepackage[utf8]{inputenc}
\usepackage[T1]{fontenc}
\usepackage{lmodern}

\usepackage{microtype}

\usepackage{xspace}

\usepackage{eurosym}
\usepackage{etoolbox} 
\newtoggle{draft}

\usepackage{calc}
\usepackage{wrapfig}

\usepackage{csquotes} 
\usepackage{chngcntr} 

\usepackage{amsmath}
\usepackage{amsfonts}
\usepackage{amsthm}
\usepackage{amssymb}
\usepackage{mathrsfs}

\usepackage{tensor} 
\usepackage{stmaryrd} 
\usepackage{commath} 

\usepackage{tikz}

\usetikzlibrary{calc,through,backgrounds,matrix,arrows}

\pgfarrowsdeclarecombine*[-.25ex]
{double angle 90}{double angle 90}{angle 90}{angle 90}{angle 90}{angle 90}

\tikzset{arrow/.style={>=to,->}}
\tikzset{mono/.style={arrow,>->}}
\tikzset{epi/.style={arrow,-double angle 90}}
\tikzset{identity/.style={double distance=1.5pt}}
\tikzset{desc/.style={fill=white,inner sep=2pt}}
\tikzset{cross line/.style={preaction={draw=white, -,
line width=6pt}}}

\swapnumbers 
\theoremstyle{definition}
\newtheorem{num}{\unskip}[section] 

\newtheoremstyle{plainEmpowered}
     {\topsep}
     {\topsep}
     {\slshape}
     {}
     {\bfseries}
     {.}
     { }
     {}

\swapnumbers
\theoremstyle{plainEmpowered}
\newtheorem{prop}[num]{Proposition}
\newtheorem{lem}[num]{Lemma}
\newtheorem{cor}[num]{Corollary}
\newtheorem{thm}[num]{Theorem}
\newtheorem{conj}[num]{Conjecture}

\theoremstyle{definition}
\newtheorem{déf}[num]{Definition}
 
\newtheorem{rmq}[num]{Remark}
\newtheorem{eg}[num]{Example}

\numberwithin{equation}{num}
\numberwithin{figure}{num}

\colorlet{step1}{green!50!black}
\colorlet{step2}{red!50!black}

\colorlet{voir}{green!50!black}

\DeclareMathOperator{\id}{id}
\DeclareMathOperator{\isDef}{\mathrel{\mathop:}=}
\DeclareMathOperator{\defIs}{=\mathrel{\mathop:}}
\DeclareMathOperator{\car}{char}
\DeclareMathOperator{\Frac}{Frac}
\DeclareMathOperator{\rank}{rank}
\DeclareMathOperator{\rk}{rank}
\DeclareMathOperator{\Cok}{Coker}

\newcommand{\rest}[2]{\left.#1\right|_{#2}}
\newcommand{\isoto}{\ensuremath{\cong}}

\newcommand{\norme}[1]{\ensuremath{{\lvert {#1} \rvert}}}

\newcommand{\Hom}{\ensuremath{\mathrm{Hom}}}

\newcommand{\mfk}[1]{\ensuremath{\mathfrak{#1}}}
\newcommand{\mcl}[1]{\ensuremath{\mathcal{#1}}}
\newcommand{\mbb}[1]{\ensuremath{\mathbb{#1}}}
\newcommand{\mrm}[1]{\ensuremath{\mathrm{#1}}}

\newcommand{\ovl}[1]{\ensuremath{\overline{#1}}}
\newcommand{\wtl}[1]{\ensuremath{\widetilde{#1}}}
\newcommand{\tld}[1]{\ensuremath{\tilde{#1}}}
\newcommand{\what}[1]{\ensuremath{\widehat{#1}}}
\newcommand{\wht}[1]{\ensuremath{\widehat{#1}}}

\newcommand{\sset}{\ensuremath{\subseteq}}
\newcommand{\sups}{\ensuremath{\supseteq}}
\newcommand{\supsn}{\ensuremath{\supsetneq}}

\newcommand{\prive}{\ensuremath{\backslash}}

\newcommand{\llb}{\ensuremath{\llbracket}}
\newcommand{\rrb}{\ensuremath{\rrbracket}}

\newcommand{\llpar}{(\negthinspace(}
\newcommand{\rrpar}{)\negthinspace)}

\newcommand{\vphi}{\ensuremath{\varphi}}

\renewcommand{\iff}{\Leftrightarrow}

\mathchardef\mhyp="2D 

\DeclareMathOperator{\Spec}{Spec}
\DeclareMathOperator{\Spf}{Spf}
\DeclareMathOperator{\Sp}{Sp}
\DeclareMathOperator{\spFibre}{sp}

\DeclareMathOperator{\Trace}{Trace}

\DeclareMathOperator{\funcGr}{Gr}
\DeclareMathOperator{\ord}{ord}
\DeclareMathOperator{\ac}{ac}

\DeclareMathOperator{\divr}{div}
\DeclareMathOperator{\dlog}{dlog}

\newcommand{\Gal}{\ensuremath{\mathrm{Gal}}}

\newcommand{\Divr}{\ensuremath{\mathrm{Div}}}

\newcommand{\catAlg}{\ensuremath{\mathrm{Alg}}}

\newcommand{\catMnd}{\ensuremath{\mathrm{Mnd}}}
\newcommand{\catFan}{\ensuremath{\mathrm{Fan}}}

\newcommand{\catAb}{\ensuremath{\mathrm{Ab}}}
\newcommand{\catMndSp}{\ensuremath{\mathrm{MndSp}}}
\newcommand{\catSet}{\ensuremath{\mathrm{Set}}}

\newcommand{\rig}{\ensuremath{\mathrm{rig}}}
\newcommand{\sm}{\ensuremath{\mathrm{sm}}}
\newcommand{\sh}{\ensuremath{\mathrm{sh}}}
\newcommand{\an}{\ensuremath{\mathrm{an}}}

\newcommand{\gp}{\ensuremath{\mathrm{gp}}}

\newcommand{\sat}{\ensuremath{\mathrm{sat}}}

\newcommand{\fs}{\ensuremath{\mathrm{fs}}}
\newcommand{\loc}{\ensuremath{\mathrm{loc}}}
\newcommand{\torsion}{\ensuremath{(\mathrm{torsion})}}
\newcommand{\logr}{\ensuremath{\mathrm{log}}}
\newcommand{\intg}{\ensuremath{\mathrm{int}}}

\newcommand{\tpl}{\ensuremath{\mathrm{top}}}

\newcommand{\equiGroth}[1]{\ensuremath{\mcl{M}_{{#1}_k\times\mbb{G}_m}^{\mbb{G}_m}}}
\newcommand{\groth}[1]{\ensuremath{K_0^R(\mathrm{Var}_{#1})}}
\newcommand{\grothModL}{\ensuremath{\groth{k}/(\mbb{L}-1)}}

\newcommand{\approche}{\ensuremath{\bigl(X,(W_\alpha,f_\alpha)_\alpha\bigr)}}

\newcommand{\logs}[1]{\ensuremath{{{#1}^\dagger}}}

\DeclareMathOperator{\hgt}{ht}

\newcommand{\reg}{\ensuremath{\mathrm{reg}}}

\newcommand{\toisoto}{\ensuremath{\mathbin{{\tikz[baseline] \path[arrow] (0pt,.4ex) edge node[above,inner ysep=.1ex,outer ysep=.15ex,font=\small] {$\isoto$} (4ex,.4ex);}}}}

\DeclareMathOperator{\Vol}{Vol}

\usepackage{todonotes}
\begin{document}


\makefrontcoverXII

\maketitle

\frontmatter 

\includepreface{preface}
\includeabstract{abstract}
\includeabstractnl{abstractnl}



\tableofcontents



\mainmatter 


%
  \graphicspath{{\chaptersdir/introduction/image/}}%
  \chapter*{Introduction}\label{ch:introduction}

\counterwithout{subsection}{section}
\counterwithin{num}{chapter}
\addcontentsline{toc}{chapter}{Introduction}
Motivic integration and logarithmic geometry: the former provides the objects that we study, whereas the latter provides the methods that we use.
The logarithmic viewpoint suggests to work on compactifications of the spaces of interest in order to derive properties about them.
Such an approach is taken in chapter~\ref{ch:serre}, whereas chapter~\ref{ch:log} uses more specifically the theory of logarithmic structures of Fontaine and Illusie.
They both focus on the study of particular motivic objects: chapter~\ref{ch:serre} is centred around the motivic Serre invariant, and chapter~\ref{ch:log} deals with motivic zeta functions.
In this general introduction, we guide the reader through the motivic universe that gave birth to these objects.
This will also serve as motivation for our work.


\subsection{Motivic zeta functions}

Let $f(x)\in \mbb{Z}[x_1,\dots, x_m]$ be a polynomial over~$\mbb{Z}$ and fix a prime number~$p$.
Denote by $N_l$ the number of solutions of the congruence $f(x) \equiv 0 \mod p^l$.
A common procedure to understand the asymptotic behaviour of a sequence of numbers $(a_l)_l$ is to consider the associated generating series $\sum_{l=0}^\infty a_l\,T^l$.
Rationality of the series is equivalent to the sequence $(a_l)_l$ being a linear recurrence sequence.

In order to prove the rationality of the series $P_f(T) = \sum_{l=0}^\infty \frac{N_l}{p^{lm}}\,T^l$, Igusa studied the \emph{Igusa local zeta function} $Z_f^p(s)$ of $f$ introduced by Weil in \cite{172} in 1965.
The function $Z_f^p(s)$ satisfies the relation
\[Z_f^p(s) = P_f(p^{-s}) - \frac{P_f(p^{-s}) - 1}{p^{-s}},\]
hence its rationality in $p^{-s}$ is equivalent to rationality of $P_f(T)$.
The big advantage of $Z_f^p(s)$ over $P_f(T)$ is that it can be expressed as a $p$-adic integral and is thus at the mercy of a well-established theory.
This approach proved fruitful in 1974, when Igusa was able to prove the rationality of $Z_f^p(s)$ using Hironaka's resolution of singularities.
He was followed ten years later by Denef, through completely different methods.

The field attracted renewed interest when Kontsevich, inspired by $p$-adic integration, initiated the theory of motivic integration during a lecture at Orsay in 1995. 
In this new theory, the integration space shifts from a power of the $p$-adic integers $\mbb{Z}_p^m$ to $(\mbb{C}\llb t \rrb)^m$, or, more generally, the arc space $\mcl{L}(X)$ of a $\mbb{C}$-variety~$X$.
The other key difference with $p$-adic integration is that integrals take values in the localized \emph{Grothendieck ring of $\mbb{C}$-varieties} $\mcl{M}_\mbb{C}$ instead of~$\mbb{R}$, whence the name \emph{motivic}.
This new framework allowed naturally Denef and Loeser to define a motivic analogue of $Z_f(s)$ in this new context; the motivic zeta function was born.

We give more details on how the motivic zeta function is defined.
Let $X$ be a smooth irreducible variety of dimension $m$ over a field $k$ of characteristic zero and consider a dominant morphism
\begin{equation}\label{XtoA1}
f\colon X\to \mbb{A}^1_k.
\end{equation}
If $X = \mbb{A}^m_k$, then $f$ corresponds to a nonconstant polynomial over~$k$. 
We set $X_s \isDef f^{-1}(0)$, the fibre over the point $x=0$.
For $d\geq 0$ we denote by $\mcl{L}_d(X)$ the \emph{$k$-scheme of $d$-jets} representing the functor
\[(k\mhyp\text{algebras})\to (\text{Sets}),\quad A\mapsto \Hom_k\bigl(\Spec A[t]/(t^{d+1}), X\bigr).\]
We set
\[\mcl{X}_{d,1} \isDef \Bigl\{\vphi\in \mcl{L}_d(X)\Bigm| f\bigl(\vphi(t)\bigr) = t^d\mod t^{d+1}\Bigr\}\]
and define the \emph{motivic zeta function} of $f$ as
\[Z_f(T) \isDef \sum_{d=1}^\infty\, [\mcl{X}_{d,1}]\, \mbb{L}^{-md}\, T^d \in \mcl{M}_{X_s}\llb T\rrb,\]
where brackets denote classes in the localization $\mcl{M}_{X_s} = K_0(\mathrm{Var}_{X_s})[\mbb{L}^{-1}]$ of the Grothendieck ring of $X_s$-varieties and $\mbb{L}\isDef [\mbb{A}^1_{X_s}]$ (see~\ref{serre:déf:Groth}).
The function $Z_f(T)$ can be seen as a motivic avatar of $Z_f^p(s)$ since it specializes to it for almost all primes $p$, as explained in \cite[6.3]{144}. 
Similarly to the $p$-adic case, Denef and Loeser showed that $Z_f(T)$ is rational.
Their proof yields an explicit formula for $Z_f(T)$ computed on a resolution of singularities for~$f$.

\begin{thm}\label{intro:motzetaDL}
Let $h\colon Y\to X$ be an embedded resolution of singularities for $(X,X_s)$, where $X_s$ is the zero locus of~$f$.
Let $(E_i)_I$ be the irreducible components of the snc divisor $h^{-1}(X_s)$ and $N_i$ their multiplicities.
Then
\[Z_f(T) = \sum_{\emptyset\neq J\sset I} (\mbb{L}-1)^{\norme{J}-1} [\wtl{E_J^\circ}] \prod_J \frac{\mbb{L}^{-\nu_i} T^{N_i}}{1- \mbb{L}^{-\nu_i} T^{N_i}} \in \mcl{M}_{X_s}\llb T\rrb.\]
\end{thm}
We briefly explain the undefined elements in the formula.
For $J\sset I$ we set $E_J = \bigcap_{i\in J} E_i$ and $E_J^\circ = E_J\prive \bigcup_{i\notin J} E_i$.
The variety $\wtl{E_J^\circ}$ is then a certain Galois cover of $E_J^\circ$ of degree $\gcd(N_j)_J$ (see construction in~\ref{serre:EJtilde} and~\ref{log:setUpSnc}).
Finally the integers $\nu_i$ are defined by the equality of divisors $K_{Y/X} = \sum_I (\nu_i-1)\, E_i$, where $K_{Y/X}$ denotes the relative canonical divisor of~$h$.

We also have a local version of $Z_f(T)$.
For every $x\in X_s$, we set
\[\mcl{X}_{d,1,x} \isDef \bigl\{\vphi\in\mcl{X}_{d,1}\bigm| \pi_0(\vphi) = x\bigr\},\] where $\pi_0(x)$ denotes the image under $\vphi\colon \Spec K[t]/(t^{d+1})\to X$ of the maximal ideal $(t)$, where $K$ is some field extension of~$k$.
The \emph{local zeta function of $f$ at~$x$} is
\[Z_{f,x}(T) \isDef \sum_{d=1}^\infty\,[\mcl{X}_{d,1,x}]\, \mbb{L}^{-md}\, T^d \in \mcl{M}_{x}\llb T\rrb.\]
It can be obtained from $Z_f(T)$ by means of the base change morphism $\mcl{M}_{X_s}\llb T\rrb \to \mcl{M}_x\llb T\rrb$.

\subsection{Milnor fibre and the monodromy conjecture}

Motivic zeta functions can be used to produce birational invariants of hypersurface singularities, but what makes them even more fascinating is how their poles are expected to carry topological information on~$f$ through the monodromy conjecture.
Before stating the conjecture, we need to recall a topological construction of Milnor.

Let $f\colon \mbb{C}^m\to\mbb{C}$ be an analytic function and fix a point $x$ of $X_s = f^{-1}(0)$.
Consider an open disc $D = B(0,\eta)\sset \mbb{C}$ around the origin and an open ball $B = B(x,\epsilon)$ around the point $x$.
We set $X' = B\cap f^{-1}(D^\times)$, where $D^\times = D\prive \{0\}$.
An important result of Lê (see \cite[1.1]{180}) states that the restriction $f_x\colon B\cap f^{-1}(D^\times)\to D^\times$ of $f$ is a locally trivial fibration for sufficiently small $\eta$ and $\epsilon$, with $0<\eta\ll\epsilon$.
We call it the \emph{Milnor fibration} at the point~$x$ (see \cite[\S 3.1]{174} for more details on Milnor fibrations).
In particular, all fibres of $f_x$ are homeomorphic to each other.
In order to get a canonical object representing the homotopy type of an arbitrary fibre of~$f_x$, we consider the fibred product
\[F_x \isDef \bigl(B\cap f^{-1}(D^\times)\bigr)\times_{D^\times}\wtl{D^\times},\]
where $\wtl{D^\times}$ denotes the universal cover of the pointed disc $D^\times$. We call $F_x$ the \emph{canonical Milnor fibre} at~$x$.
Since $\wtl{D^\times}$ is contractible, $F_x$ has the same homotopy type as each fibre of $f_x$.
In particular, they have the same singular cohomology.

The fundamental group $\pi_1(D^\times)$, which is isomorphic to the deck transformation group of~$\wtl{D^\times}$, acts on $F_x$.
In particular, the canonical generator of $\pi_1(D^\times)\isoto \mbb{Z}$, given by a positively oriented loop turning once around the origin, induces by functoriality a linear transformation for every $i\geq 0$
\[M_x^i\colon H^i(F_x,\mbb{Z})\to H^i(F_x,\mbb{Z}),\]
where $H^i(F_x,\mbb{Z})$ is the $i$th singular cohomology space of $F_x$.
These operators~$M_x^i$ are called the \emph{monodromy transformations}.
A \emph{local monodromy eigenvalue} of~$f$ is then an eigenvalue of~$M_x^i$ for some $i\geq 0$ and $x\in X_s$.

When $f$ is defined over a number field, Igusa conjectured in 1988 that for almost all primes $p$, $\exp(2\pi i r_0)$ is a local monodromy eigenvalue if $r_0$ is the real part of a pole of Igusa's zeta function~$Z_f^p(s)$.
This conjecture, known as the \emph{monodromy conjecture}, raised much interest because it relates the topology of the fibre $f^{-1}(0)$ with the arithmetic of~$f$ (recall that $Z_f^p(s)$ is tightly linked with the number of solutions of $f(x) \equiv 0\mod p^l$).

In analogy to Igusa's monodromy conjecture, Denef and Loeser stated in \cite[2.4]{74} a similar conjecture for the motivc zeta function.
\begin{conj}
Let $k$ be a subfield of~$\mbb{C}$ and $f\colon X\to \mbb{A}^1_k$ a dominant morphism.
There is a finite subset $\mcl{S}$ of $\mbb{N}_{\geq 1}\oplus \mbb{N}_{\geq 1}$ such that $Z_f(T)$ belongs to the subring
\[\mcl{M}_{X_s}\bigl[T, \tfrac{1}{1-\mbb{L}^{-a}T^b}\bigr]_{(a,b)\in \mcl{S}}\]
of $\mcl{M}_{X_s}\llb T\rrb$ and such that $\exp(2\pi i\, a/b)$ is a local monodromy eigenvalue of $f^\an\colon (X\times_k\mbb{C})^\an \to \mbb{C}$ for every couple $(a,b)\in\mcl{S}$, where $(\cdot)^\an$ denotes the GAGA analytification functor.
\end{conj}

Since the motivic zeta function specializes to Igusa's $p$-adic functions for almost all $p$, interest has been shifted to this more general conjecture.
Although some partial results were obtained, mainly in low dimension (see e.g.\ \cite{175} and \cite{176}), this conjecture still seems out of reach in its full generality, where methods used in solved cases fall short.

\subsection{Analytic Milnor fibre and volume Poincaré series}
In order to get a better feeling for the monodromy conjecture, Nicaise and Sebag introduced in \cite[9.1]{106} the analytic Milnor fibre $\mcl{F}_x$.

Let $k$ be an algebraically closed field of characteristic zero and consider the complete discrete valuation ring $R = k\llb t\rrb$, with fraction field $K = k\llpar t\rrpar$.
For $d\geq 1$, we also consider the totally ramified extension $R(d)\isDef R[T]/(T^d-t)$ and set $K(d) \isDef \Frac R(d)$.
Let $f\colon X\to \mbb{A}^1_k = \Spec k[t]$ be as in~\eqref{XtoA1} and $x\in X_s$ be a point in the zero locus of~$f$.
We denote by $\what{X}$ the $t$-adic completion of $X$.
The \emph{analytic Milnor fibre} of~$f$ at~$x$ is the generic fibre $\mcl{F}_x$ of the formal affine scheme $\Spf \mcl{O}_{\what{X},x}$, in the sense of Berthelot \cite[0.2.6]{BERT} (see \cite[6.1.4]{31} for topological intuition).
It is a smooth separated rigid $K$-variety, whose $\ell$-adic cohomology coincides, when $k = \mbb{C}$, with the singular cohomology of the Milnor fibre $F_x$, and identifies the Galois action on the former to the monodromy action on the latter.
Furthermore, $K(d)$-points of $\mcl{F}_x$ correspond precisely to arcs $\vphi\colon \Spec k\llb t\rrb \to X$ with $f(\vphi) = t^d$ and $\pi_0(\vphi) = x$, so that the analytic Milnor fibre $\mcl{F}_x$ appears to be a key object connecting the two different worlds entering the monodromy conjecture: arithmetic and topology.

This brings us to the question whether $Z_{f,x}(T)$ could be expressed as a Weil-type zeta function, i.e.\ a generating series whose coefficients somehow count $K(d)$-points on~$\mcl{F}_x$.
This has led Nicaise and Sebag to introduce in \cite{106} the \emph{volume Poincaré series} of an stft formal $R$-scheme (where stft stands for separated and topologically of finite type).
Let $\mcl{X}$ be a generically smooth stft formal $R$-scheme of pure relative dimension $m$ and let $\omega$ be a volume form on the generic fibre $\mcl{X}_K$, i.e., a nowhere vanishing differential form of degree $m$.
We write $\omega(d)$ for the inverse image of $\omega$ on the generic fibre of $\mcl{X}(d) \isDef \mcl{X}\times_R R(d)$.
The \emph{volume Poincaré series} of the pair $(\mcl{X},\omega)$ is
\[S(\mcl{X},\omega;T)\isDef \sum_{d\geq 1} \Bigl(\int_{\mcl{X}(d)} \norme{\omega(d)}\Bigr)\, T^d \in \mcl{M}_{\mcl{X}_k}\llb T\rrb\]
(see \ref{serre:déf:motInt} for a definition of the motivic integral).
Its image in $\mcl{M}_k\llb T\rrb$ only depends on $\mcl{X}_K$ and not on~$\mcl{X}$.
It is denoted by $S(\mcl{X}_K,\omega;T)$.

This series can be used in the context of~\eqref{XtoA1} because every volume form~$\phi$ on the smooth variety~$X$ induces canonically a volume form $\frac{\phi}{\dif{f}}$ on the rigid $K$-variety~$(\what{X})_K$, called the \emph{Gelfand-Leray form}.
\begin{thm}[{\cite[9.10]{106}}]\label{ZisS}
Let $X$ be a smooth irreducible $k$-variety of dimension~$m$, $f\colon X\to \mbb{A}^1_k$ a dominant morphism and $\phi$ a volume form on~$X$.
Then
\[S\bigl(\what{X},\tfrac{\phi}{\dif{f}}; T\bigr) = \mbb{L}^{-(m-1)}Z_f(\mbb{L} T)\in \mcl{M}_{X_s}\llb T \rrb.\]
\end{thm}

In \cite{49} Nicaise extended the definition of the volume Poincaré series to special formal $R$-schemes.
If $\mcl{X}$ is stft and $Z$ is a locally closed subset of the special fibre~$\mcl{X}_k$, then the completion $\mcl{Z}$ of~$\mcl{X}$ along $Z$ is a special formal $R$-scheme and $S(\mcl{Z},\omega;T)$ coincides with the image of $S(\mcl{X},\omega;T)$ under the base change morphism $\mcl{M}_{\mcl{X}_k}\llb T\rrb\to \mcl{M}_Z\llb T\rrb$.
This generalization is more than welcome because $\mcl{F}_x$ is canonically isomorphic to the generic fibre of the completion of $X$ along the closed point~$x$ and as such is a special formal scheme.
In particular, the image of $S(\what{X},\tfrac{\phi}{\dif{f}};T)$ in $\mcl{M}_x$ only depends on $\mcl{F}_x$ (cfr.\ \cite[9.8]{49}) and we have
\[Z_{f,x}(\mbb{L} T) = \mbb{L}^{m-1}\, \sum_{d=1}^\infty \Bigl(\int_{\mcl{F}_x(d)} \norme{\tfrac{\phi}{\dif{f}}(d)}\Bigr)\, T^d \in \mcl{M}_x\llb T\rrb.\]


\subsection{Monodromy conjecture for degenerations of Calabi-Yau varieties}

The advantage of $S(\mcl{X},\omega;T)$ over $Z_f(T)$ is that it admits a natural generalization to $K$-varieties.
Let $X$ be a Calabi-Yau variety over $K = k\llpar t\rrpar$, i.e.\ a smooth, proper and geometrically connected variety with trivial canonical sheaf. Assume that $X$ has dimension $m$ and let $\omega$ be volume form on~$X$.
The \emph{motivic zeta function} of the pair~$(X,\omega)$ is defined as
\[Z_{X,\omega}(T) \isDef \mbb{L}^{m} \sum_{d=1}^\infty \Bigl(\int_{X(d)}\norme{\omega(d)}\Bigr)\, T^d.\]
This generalization comes with a variant of the monodromy conjecture, as introduced by Halle and Nicaise in \cite[2.6]{164}.
\begin{déf}\label{intro:monodPropty}
Let $\sigma$ be a topological generator of the Galois group $G(K^s/K)$, where $K^s$ denotes a seperable closure of~$K$.
We say that $X$ satisfies the \emph{global monodromy property} if there is a finite subset $\mcl{S}$ of $\mbb{Z}\oplus \mbb{Z}_{\geq 1}$ such that $Z_{X,\omega}(T)$ belongs to the subring
\[\mcl{M}_{k}\bigl[T, \tfrac{1}{1-\mbb{L}^{a}T^b}\bigr]_{(a,b)\in \mcl{S}}\]
of $\mcl{M}_{k}\llb T\rrb$ and such that the $\tau$th cyclotomic polynomial divides the characteristic polynomial of $\sigma$ on $H^i(X\times_K K^s, \mbb{Q}_\ell)$ for some $i\geq 0$, where $\tau$ denotes the order of $a/b$ in the group $\mbb{Q}/\mbb{Z}$.
\end{déf}
Note that the global monodromy property doesn't depend on the volume form~$\omega$, since for every unit $u\in K^\times$ with valuation $v_K(u)$ we have
\[Z_{X,u\cdot\omega}(T) = Z_{X,\omega}\bigl(\mbb{L}^{-v_K(u)}T\bigr) \in \mcl{M}_k\llb T\rrb.\]

Halle and Nicaise proved in \cite[8.6]{164} that the property holds for abelian varieties.
Their proof relies on a detailed analysis of the behaviour of the Néron model under base change and uses the theory of abelian varieties in an essential way.
It remains a challenge to prove the property in greater generality; we will explain in the last section of this introduction how our work provides a step forward in the case of Calabi-Yau varieties.

\subsection{Motivic nearby cycles and motivic volume}

Guided by previous work on Igusa zeta functions (see \cite[\S 4]{74}), Denef and Loeser considered the opposite of the limit for $T\to\infty$ of $Z_f(T)$ and gave heuristic evidence for it to be considered as a motivic incarnation of the complex of nearby cycles of~$f$ (see \cite[3.5.3]{72}).
From the description of $Z_f(T)$ in \ref{intro:motzetaDL} we get
\[\mcl{S}_f\isDef -\lim_{T\to\infty} Z_f(T) = \sum_{\emptyset\neq J\sset I} (1-\mbb{L})^{\norme{J}-1} [\wtl{E_J^\circ}]\in \mcl{M}_{X_s},\]
which is called the \emph{motivic nearby fibre} of~$f$.
This construction inspired Nicaise and Sebag, who defined in \cite[8.3]{106} the motivic volume of a generically smooth stft $R$-scheme in a similar way.
More details can be found in~\ref{serre:nearbyCycles}.

\begin{déf}
Let $\mcl{X}$ be a generically smooth stft formal $R$-scheme.
The \emph{motivic volume} of $\mcl{X}$ is
\[\Vol(\mcl{X})\isDef -\lim_{T\to\infty} S(\mcl{X},\omega;T)\in \mcl{M}_{\mcl{X}_k},\]
where $\omega$ is any volume form on $\mcl{X}_K$.
Its image in $\mcl{M}_k$ only depends on $\mcl{X}_K$ and not on $\mcl{X}$.
It is denoted by $\Vol(\mcl{X}_K)$ and called the \emph{motivic volume} of~$\mcl{X}_K$.
\end{déf}

Theorem \ref{ZisS} then yields easily the following comparison result.
\begin{thm}[{\cite[9.8]{49}}]\label{SfmotVolume}
Let $f\colon X\to \mbb{A}^1_k$ be a dominant morphism from a smooth $k$-variety~$X$ of dimension~$m$.
We have
\[\Vol(\what{X}) = \mbb{L}^{-(m-1)} \mcl{S}_f \in \mcl{M}_{X_s}\]
and for every closed point $x\in X_s$,
\[\Vol(\mcl{F}_x) = \mbb{L}^{-(m-1)} \mcl{S}_{f,x} \in \mcl{M}_{x}.\]
\end{thm}

In \cite[3.7]{53} Guibert, Loeser and Merle extended $\mcl{S}_f$ to a morphism
\[\mcl{S}_f\colon \mcl{M}_X\to \mcl{M}_{X_s}\]
in such a way that $\mcl{S}_f([X]) = \mcl{S}_f$.
In the process, they defined the \emph{motivic nearby cycles with support} for a dense open subvariety $U$ of~$X$.
It is constructed as the limit of a modified motivic zeta function $Z_{f,U}(T)$ whose coefficients discard arcs too tangent to $X\prive U$.
We will explain in section~\ref{contributions} how we could fit this construction in Nicaise and Sebag's framework.

\subsection{Motivic Serre invariant and the trace formula}

Let $R$ be a complete discrete valuation ring with perfect residue field $k$ and denote by $K$ its fraction field.
The motivic Serre invariant, which is the subject of chapter~\ref{ch:serre}, was first defined in \cite{98} for a smooth separated quasi-compact rigid $K$-variety~$X$.
If $\mcl{X}$ is a weak Néron model for $X$, then $S(X)$ is given by the class of the special fibre of $\mcl{X}$ in the quotient Grothendieck ring $K_0^R(\mrm{Var}_k)/(\mbb{L}-1)$ (see notation in~\ref{serre:K0VarR}).
The independence of $S(X)$ from $\mcl{X}$ is a nontrivial statement which can be proved using the theory of motivic integration (see also \ref{serre:SXoriginel}).

A weak Néron model for a separated rigid $K$-variety $X$ is a smooth stft formal $R$-model of~$X$ whose unramified points are in bijection with the unramified points of~$X$ (see \ref{serre:déf:weakNéron} for a precise definition).
In particular, if $X(K') = \emptyset$ for every unramified extension $K'/K$ then the empty formal scheme is a weak Néron model for~$X$, so that $S(X) = [\emptyset] = 0$.
Hence the nonvanishing of $S(X)$ implies the existence of an unramified point on~$X$.
In particular, $S(X)$ can detect rational points when the residue field $k$ is algebraically closed.

Serre's name was attached to this invariant for its analogy with the $p$-adic invariant $S_p(\cdot)$ defined by Serre in \cite{168} to classify compact $p$-adic manifolds.
When $K$ is a $p$-adic field, the set of rational points of a smooth rigid $K$-variety~$X$ has a structure of $p$-adic manifold and the correspondance
\[Y\mapsto \norme{Y(k)}\]
determines a well-defined morphism
\[K_0^R(\mrm{Var}_k)/(\mbb{L}-1)\to \mbb{Z}/(q-1),\]
where $q$ is the cardinality of $k$.
One shows that $S(X)$ is sent onto $S_p(X(K))$ under this morphism.

In further work, Nicaise and Sebag defined the motivic Serre invariant for larger classes of objects: generically smooth stft formal schemes in \cite{145} and  generically smooth special formal schemes in \cite{49}.
Finally, Nicaise defined $S(X)$ in \cite{112} when $X$ is a smoothly bounded rigid or algebraic $K$-variety.
Typical examples of smoothly bounded algebraic varieties are given by proper smooth varieties.
He also showed how to extend the invariant additively to a morphism of rings $S\colon \mcl{M}_K\to K_0^R(\mrm{Var}_k)/(\mbb{L}-1)$ when $K$ has characteristic zero.
This extension yields an obvious definition of $S(X)$ for a general algebraic variety~$X$, with no smoothness or properness condition, which can still be used to detect unramified points on~$X$.
It also gives a new realization morphism that can be used to distinguish elements in~$\mcl{M}_K$.

The following proposition hints at why we can consider $S(X)$ not only as an unramified-point detector, but as a true measure of the set of unramified points on $X$.
\begin{prop}[{\cite[2.3/5]{NER}}]
Let $\mcl{X}$ be a smooth $R$-scheme over a complete discrete valuation ring.
Then the canonical map
\[\mcl{X}(R^\sh)\to \mcl{X}_k(k^s)\]
is surjective, where $R^\sh$ denotes a strict henselisation of~$R$ and $k^s$ its residue field, which is a separable closure of the residue field~$k$ of~$R$.
\end{prop}

Hence, we see that the unramified points of~$\mcl{X}$ are parameterized by $\mcl{X}_k(k^s)$.
Since $k$ is assumed to be perfect, every closed point $x\in \mcl{X}_k$ induces a $k^s$-point on~$\mcl{X}_k$, so that $[\mcl{X}_k]$ can be seen as a good measure for $\mcl{X}(R^\sh)$.
This can be made more precise in rigid terms.
If $x$ is a closed point of $\mcl{X}_k$, then the inverse image of~$x$ under the specialization morphism $\spFibre\colon (\what{\mcl{X}})_K\to \mcl{X}_k$ is an open rigid unit ball (see \cite[2.1.5]{32}).

Since $S(X)$ measures the set of rational points when $k$ is algebraically closed, we can ask ourselves whether it admits a cohomological interpretation, as does the set of rational points on a variety over a finite field with the Grothendieck-Lefschetz formula.

\begin{thm}[{\cite[Rapport, 3.2]{SGA4.5}}] 
Let $X$ be a variety over $\mbb{F}_q$ and $\ovl{X}$ the base change of $X$ to an algebraic closure of $\mbb{F}_q$.
Let $F\colon \ovl{X}\to\ovl{X}$ be the absolute Frobenius automorphism.
Then for every integer $d\geq 1$ we have
\[\norme{X(\mbb{F}_{q^d})} = \sum_i (-1)^i \Trace\bigl(F^d\bigm| H_c^i(\ovl{X},\mbb{Q}_\ell)\bigr),\]
where $H_c^i(\cdot,\mbb{Q}_\ell)$ denotes $\ell$-adic cohomology with compact support.
\end{thm}

The classical Lefschetz trace formula in topology allows to count the number of fixed points of a continuous map $f\colon X\to X$ by means of the traces of the operator induced by $f$ on the homology spaces of $X$.
Since fixed points of Frobenius correspond to rational points of $X$, the Grothendieck-Lefschetz formula allows similarly to compute rational points from traces of Frobenius.

Nicaise and Sebag obtained a first result in that direction in \cite[5.4]{106}, which Nicaise improved in \cite{49} and \cite{112}. 
We will only state the formula in the algebraic setting.
A smooth and proper algebraic $K$-variety $X$ is called \emph{tame} if there is a regular proper $R$-scheme $\mcl{Y}$ such that $\mcl{Y}_K\isoto X$ and $\mcl{Y}_k$ is a strict normal crossings divisor, whose multiplicities are prime to the characteristic exponent of~$k$.
If $\car k = 0$, then every smooth proper $K$-variety is tame by resolution of singularities.

\begin{thm}[{\cite[6.3]{112}}]
Assume $k$ algebraically closed.
Let $X$ be a smooth and proper tame $K$-variety and $\vphi$ a topological generator of the tame monodromy group $\Gal(K^t/K)$.
Then for every $d\geq 1$ not divisible by~$p$,
\[\chi_\tpl\bigl(S(X\times_K K(d))\bigr) = \sum_i (-1)^i\Trace\bigl(\vphi^d\bigm| H^i(X\times_K K^t,\mbb{Q}_\ell)\bigr).\]
\end{thm}
The extension of the motivic Serre invariant to a morphism on the Grothendieck ring of $K$-varieties has allowed Nicaise to extend his trace formula to arbitrary $K$-varieties when $k$ has characteristic zero.

\begin{cor}
If $k$ has characteristic zero, then for every $K$-variety $X$
\[\chi_\tpl(S(X)) =  \sum_i (-1)^i \Trace\bigl(\vphi\bigm| H^i(X\times_K K^s,\mbb{Q}_\ell)\bigr).\]
\end{cor}

A similar formula for tame generically smooth special formal schemes (see \cite[6.4]{49}) yields a trace formula for the analytic Milnor fibre $\mcl{F}_x$, which exhibits a tight link between $K(d)$-points on $\mcl{F}_x$ and its cohomology.

Finally, one can show that if a (rigid or algebraic) $K$-variety~$X$ admits a gauge form $\omega$, then
\[S(X) = \int_X \norme{\omega}\:\in K_0^R(\mrm{Var}_k)/(\mbb{L}-1).\]
Together with theorem \ref{ZisS}, this yields the remarkable formula
\[Z_{f,x}(T) = \sum_{d=1}^\infty S(\mcl{F}_x(d))\, T^d \:\in K_0(\mrm{Var}_x)/(\mbb{L}-1)\llb T\rrb,\]
which makes the monodromy conjecture appear more convincing than ever.

The following theorem of Nicaise (extending a theorem of A'Campo in \cite{181}) gives another occurence of the fantastic interplay between the objects we defined in this introduction.
\begin{thm}[{\cite[8.10]{49}}]
Let $f\colon X\to \Spec \mbb{C}[t]$ be a dominant morphism from a smooth irreducible $\mbb{C}$-variety.
Let $x$ be a closed point of $f^{-1}(0)$ and denote by $F_x$  (resp.\ $\mcl{F}_x$) the (resp.\ analytic) Milnor fibre of~$f$ at~$x$.
The following are equivalent:
\begin{enumerate}
\item the morphism $f$ is smooth at~$x$;
\item $\mcl{F}_x$ contains a $\mbb{C}\llpar t\rrpar$-rational point;
\item $\mcl{F}_x$ satisfies $S(\mcl{F}_x)\neq 0$;
\item $\mcl{F}_x$ satisfies
\[\sum_i (-1)^i \Trace\bigl(\vphi\bigm| H(\mcl{F}_x\what{\times}_{\mbb{C}\llpar t\rrpar} \what{\mbb{C}\llpar t\rrpar^s},\mbb{Q}_\ell)\bigr)\neq 0,\]
where $\vphi$ denotes any topological generator of the monodromy group $\Gal\bigl(\mbb{C}\llpar t\rrpar^s/\mbb{C}\llpar t\rrpar\bigr)$;
\item $F_x$ satisfies
\[\sum_i (-1)^i \Trace\bigl(M\bigm| H(F_x,\mbb{C})\bigr)\neq 0,\]
where $M$ denotes the monodromy transformation on the graded singular cohomology. 
\end{enumerate}

\end{thm}

\subsection{Contributions}
\label{contributions}
As we explained, Nicaise could extend the motivic Serre invariant to a morphism
\[S\colon \mcl{M}_K\to K_0^R(\mrm{Var}_k)/(\mbb{L}-1)\]
when $K$ has characteristic zero.
The main ingredient of his proof is Bittner's presentation of the Grothendieck ring by classes of proper smooth varieties (see \cite[3.1]{86}).
This presentation relies on resolution of singularities and on the weak factorization theorem (which states that every birational morphism between proper smooth $K$-varieties can be factored as a sequence of blow-ups and blow-downs with smooth centres), two results depending on the assumption $\car K = 0$ in an essential way.
While resolution of singularities is widely thought to be also true in positive characteristic, the status of the weak factorization theorem is less clear.

In chapter~\ref{ch:serre}, we give a new construction of $S(U)$ for a smooth $K$-variety~$U$ under the sole assumption that resolution of singularites holds for $R$-schemes of finite type.
Our approach is logarithmic in spirit: we consider a compactification $X$ of~$U$, and use its boundary to approximate $U^\rig$ by a sequence of quasi-compact rigid open subvarieties $(U_\gamma)$ for which the Serre invariant can be defined in a direct way.
After showing that $S(U_{\gamma})$ stabilizes for $\gamma$ sufficiently large, we can define $S(U)$ as the stabilizing value.
This construction also yields an extension of $S(U)$ to a morphism on $\mcl{M}_K$, which coincides with Nicaise's.
As a first application, we explain how our construction allows to extend the arguments of Esnault and Nicaise in~\cite{120} to give a conditional affirmative answer of Serre's question on the existence of rational fixed points of a $G$-action on the affine space, for $G$ a finite $\ell$-group.
Finally, the last section of the chapter is devoted to motivic nearby cycles. We use there the same logarithmic approach to give a new way of understanding the motivic nearby cycles with support of Guibert, Loeser and Merle as a motivic volume, in a similar way to theorem \ref{SfmotVolume}.

Chapter~\ref{ch:log} focuses on another side of the motivic bestiary: motivic zeta functions.
It expands on results announced in~\cite{186}.
We give there a new explicit expression for the volume Poincaré series $S(\mcl{X},\omega;T)$ when $\mcl{X}$ is a generically smooth log smooth $R$-scheme.
By adding some structure to schemes, logarithmic geometry allows to handle so-called log smooth schemes as if they were smooth.
A key feature of log smooth schemes is that they can be desingularized by means of their combinatorics.
In this process of desingularizing, fake poles are intoduced in the expression of $S(\mcl{X},\omega;T)$, so that our formula yields a much smaller set of candidate poles than would have the corresponding desingularized scheme.
This result paves the way to a proof of the monodromy conjecture for Calabi-Yau varieties that are constructed in the Gross-Siebert programme on mirror symmetry, where log smooth models for those varieties appear naturally (see for example \cite{170}).
Finally, we show in section~\ref{log:sectionNewton} how to recover a formula of Guibert for $Z_f(T)$ when $f$ is a polynomial that is nondegenerate with respect to its Newton polyhedron (see \cite{165}) and section~\ref{log:twoVariables} sheds new light on the (solved) monodromy conjecture for polynomials in two variables.

\counterwithin{subsection}{section}
\counterwithin{num}{section}
\cleardoublepage


%

%
  \graphicspath{{\chaptersdir/formRig/image/}}%
  \chapter{Formal and rigid geometry}\label{ch:formRig}

This very short chapter is devoted to formal and rigid geometry.
It aims to give hands-on knowledge on the subject in order to make some arguments in chapter~\ref{ch:serre} more transparent to the unfamiliar reader.

\section{Formal schemes}

We refer to \cite[\S I.10]{EGA-I} for a thorough introduction to formal schemes.

\begin{num}
A standard operation in ring theory is completion along an ideal.
Let $A$ be a ring and $I$ an ideal of~$A$.
The \emph{$I$-adic completion} of~$A$ is given by the projective limit
\[\what{A}\isDef \varprojlim_n A/I^n.\]
If the canonical morphism $A \to \what{A}$ is an isomorphism, $A$ is said to be \emph{complete with respect to the $I$-adic topology}, or simply \emph{adic}. We call $I$ the \emph{ideal of definition} for the adic ring~$A$.
\end{num}

\begin{num}
Let's take a geometric point of view on this construction.
Each $\Spec A/I^n$ has the same underlying topological space~$Z$, which corresponds to a closed subspace of~$\Spec A$.
The only difference between them lies in their sheaves of sections.
One can actually see those subschemes $\Spec A/I^n$ as infinitesimal neighbourhoods of order~$n$ of~$Z$ into~$\Spec A$.
The logical counterpart to taking the projective limit of the rings $(A/I^n)_n$ would be to take the inductive limit of the schemes $(\Spec A/I^n)_n$.
The resulting object would capture the notion of infinitesimal neighbourhoud of the closed subspace $Z$ inside of~$\Spec A$.
This construction cannot be carried out in the category of schemes, but leads to the notion of formal scheme.
\end{num}

\begin{num}
The \emph{formal spectrum} $\Spf A$ of an adic ring is constructed as follows.
As a topological space it is the closed subspace of $\Spec A$ consisting of all \emph{open} prime ideals of $A$, i.e.\ all prime ideals of $A$ containing the ideal of definition~$I$.
For $f\in A$ we denote by $\mcl{D}(f)$ the subset of $\Spf A$ consisting of all open prime ideals containing $f$.
The sheaf $\mcl{O}_{\Spf{A}}$ is then the sheafification of the presheaf given by $\mcl{O}(\mcl{D}(f)) = \what{A_f}$, where $A_f$ denotes the localization of $A$ along the multiplicative subset $\{f^n\}_{n\geq 1}$, and $\what{A_f}$ is its $I$-adic completion.
The topological space $\Spf A$ endowed with the sheaf $\mcl{O}_{\Spf A}$ is called an \emph{affine formal scheme}.
\end{num}

\begin{num}
A \emph{locally Noetherian formal scheme} is a ringed space locally isomorphic to affine formal schemes~$\Spf A$ with $A$ adic and Noetherian.
We will always assume our formal schemes locally Noetherian and simply call them \emph{formal schemes}.
Moreover, we define unsurprisingly a \emph{morphism of formal schemes} as a morphism of locally ringed spaces.
\end{num}

The following class of rings is the formal counterpart to polynomials rings.

\begin{déf}\label{restSeries}
Let $A$ be an adic ring and $I$ an ideal of definition.
A power series $f = \sum_\nu c_\nu T^\nu \in A\llb T_1,\dots, T_r\rrb$ is called \emph{restricted} if for every $n\geq 1$ there are only finitely many coefficients of~$f$ not belonging to $I^n$.
In other words, the sequence $(c_\nu)\to 0$ as $\norme{\nu}\to \infty$.

The set of restricted formal power series is a subring of $A\llb T_1,\dots, T_r\rrb$, denoted by $A\{T_1,\dots,T_r\}$.
Clearly, if the topology on $A$ is discrete, i.e.\ defined by the zero ideal, then $A\{T_1,\dots,T_r\} = A[T_1,\dots,T_r]$.
\end{déf}

The following proposition shows that the ring of restricted power series is actually the $I$-adic completion of the polynomial ring.

\begin{prop}[{\cite[III \S4.2 prop.~3]{AC1-7}}] \label{restrSeriesIsCompletion}
Let $A$ be an adic ring and $I$ an ideal of definition.
Then $A\{T_1,\dots,T_r\} = \varprojlim_n A[T_1,\dots,T_r]/I^n$.
\end{prop}

\begin{num}
Let $R$ be a complete discrete valuation ring and $\pi$ a uniformizer. It is in particular an adic ring for the ideal generated by~$\pi$, and we can consider its formal spectrum~$\Spf R$.
From now on we will restrict our study to formal schemes over the formal spectrum of~$R$, also called \emph{formal $R$-schemes}, which are the only kind of formal schemes we will encounter in the next chapters.
Note that all the notions that we present below apply much more widely.
\end{num}

\begin{prop}[{\cite[0.10.13.1]{EGA-I}}] \label{finiteTypeCdtn}
A morphism of formal schemes $f\colon \mcl{X}\to \Spf R$ is said to be \emph{topologically of finite type} if $\mcl{X}$ can be covered by finitely many affine formal opens $\mcl{V}_i$ such that $\mcl{O}_\mcl{X}(\mcl{V}_i) = R\{t_1,\dots, t_r\}/J$ for some integer $r\geq 1$ and some ideal~$J$.
\end{prop}

\begin{num}
Let $\mcl{X}\to \Spf R$ be a formal $R$-scheme topologically of finite type. The ringed space $(\mcl{X},\mcl{O}_\mcl{X}/\pi)$ is a scheme of finite type over the residue field~$k$ of~$R$.
It is denoted by~$\mcl{X}_k$ and called the \emph{special fibre} of~$\mcl{X}$.
\end{num}

We motivated the notion of formal scheme as a geometric version of completion.
This naturally raises the question whether a scheme can be completed to a formal scheme.

\begin{déf}
Let $\mcl{X}\to R$ be an $R$-scheme of finite type. Its \emph{$\pi$-adic completion} is given by the formal scheme $\what{\mcl{X}}\isDef (\mcl{X}_k,\what{\mcl{O}}_\mcl{X})$ where
\[\what{\mcl{O}}_\mcl{X}\isDef \varprojlim_n \mcl{O}_\mcl{X}/\pi^n.\]
It is topologically of finite type over~$\Spf R$.
Furthermore, the special fibre of the formal scheme $\wht{\mcl{X}}$ is nothing else but the special fibre $\mcl{X}_k$ of~$\mcl{X}$.
\end{déf}

\begin{eg}
The completion of the affine space $\mbb{A}^r_R$ over~$R$ is the formal spectrum of $R\{T_1,\dots,T_r\}$, as can be seen from proposition~\ref{restrSeriesIsCompletion}.
\end{eg}

Many classical properties of schemes can be defined for formal schemes.
The following proposition will be sufficient for our purposes.

\begin{prop}
Let $\mcl{X}\to R$ be an $R$-scheme of finite type. If $\mcl{X}$ is separated, resp.\ regular, resp.\ flat, resp.\ smooth over~$R$, so is its $\pi$-adic completion.
\end{prop}

\begin{num}\label{hintXeta}
Let $\mcl{X}$ be an $R$-scheme of finite type. When passing to the completion, we a priori lose information on the generic fibre of~$\mcl{X}$.
We will see in the next section that we can actually recover some information on it, by attaching a rigid $K$-variety to~$\mcl{X}$.
\end{num}

\section{Rigid varieties}

We refer to \cite{185} for a thorough introduction and to \cite{NAA} for a comprehensive reference.

\begin{num}
Let $K$ be the fraction field of a complete discrete valuation ring~$R$.
It comes equipped with a non-Archimedean absolute value
\[\norme{\cdot}\colon K\to \mbb{R}_{\geq 0}, \, x \mapsto e^{-v(x)},\]
where $v$ denotes the valuation of~$R$.
The field~$K$ is complete with respect to the topology induced by $\norme{\cdot}$.
\end{num}

\begin{déf}
The \emph{Tate algebra} $T_r = K\{T_1,\dots,T_r\}$ of dimension~$r$ is defined as the $K$-algebra $R\{T_1,\dots,T_r\}\otimes_R K$.

Note that it is \emph{not} a ring of restricted power series as defined in~\ref{restSeries} because $K$ is not an adic ring (the ideals $(\pi^n)$ are unit ideals).
\end{déf}

\begin{déf}
An \emph{affinoid $K$-algebra} is a quotient of some Tate algebra~$T_r$.
Let $A$ be an affinoid $K$-algebra.
We define the ring of \emph{restricted series} in $s$ variables as $A\{T_1,\dots,T_s\}\isDef A\otimes_K K\{T_1,\dots,T_s\}$.
It is clearly an affinoid $K$-algebra (if $T_r\to A$ is a surjective map, then $T_{r+s}\to A\{T_1,\dots,T_s\}$ is surjective as well).
One can show that affinoid algebras are Noetherian (see \cite[3.1.3]{185}).
\end{déf}

\begin{num}
Let $A$ be an affinoid $K$-algebra.
We denote by $\Sp A$ the set of maximal ideals of $A$ and call it the \emph{spectrum} of~$A$.
The couple $(A,\Sp A)$ is called an \emph{affinoid $K$-space}.
A \emph{morphism of affinoid $K$-spaces} is then simply a $K$-algebra morphism $\sigma\colon B\to A$.
Note that such a morphism induces a map $\Sp A\to \Sp B,\, \mfk{m}\mapsto \sigma^{-1}(\mfk{m})$.
If $\mfk{m}$ is a maximal ideal of $A$, then the quotient $A/\mfk{m}$ is a finite extension of~$K$ (see \cite[2.2.12]{185}).
For $f\in A$, we define $f(\mfk{m})$ as the image of~$f$ in~$A/\mfk{m}$.
\end{num}

\begin{num}
We could define a Zariski topology on affinoid $K$-spaces, but the whole point of rigid geometry is to provide a framework to apply analytical methods for non-Archimedean spaces.
Therefore, we will endow these spaces with the canonical topology induced from the field $K$ instead.
Let $X = \Sp A$ be an affinoid $K$-space.
A basis for the \emph{canonical topology} on~$X$ is given by the subsets
\[X(f_1,\dots,f_r) \isDef \Bigl\{x\in X\Bigm| \norme{f_i(x)}\leq 1, \text{ for every } 1\leq i \leq r\Bigr\}\]
for elements $f_i\in A$.
\end{num}

We now set out to define a structural sheaf for affinoid $K$-spaces that suits the canonical topology.

\begin{déf}
A subset $U$ of an affinoid $K$-space $X$ is called an \emph{affinoid subdomain} of~$X$ if there is a morphism of affinoid $K$-spaces $i\colon Y\to X$ such that $i(Y)\sset U$ and if $i$ is universal with respect to this property: any other morphism of affinoid $K$-spaces whose image is contained in~$U$ factors uniquely through~$i$.
\end{déf}
If $i\colon Y\to X$ defines an affinoid subdomain over~$U$, then $i$ is injective and $i(Y)\isoto U$ (see \cite[3.3.10]{185}).

\begin{prop}
Let $X = \Sp A$ be an affinoid $K$-algebra and $f = (f_i)_{1\leq i \leq r}$, $g = (g_j)_{1\leq j\leq s}$ be elements of~$A$.
Then the following subsets of~$X$ are affinoid subdomains:
\begin{itemize}
\item  $X(f) \isDef \bigl\{x\in X\bigm| \norme{f_i(x)}\leq 1, \, \forall i\bigr\}$;
\item  $X(f, g^{-1}) \isDef X(f)\cap \bigl\{x\in X\bigm| \norme{g_j(x)}\geq 1, \, \forall j\bigr\}$;
\item  $X(\tfrac{f_1}{f_0},\dots,\tfrac{f_r}{f_0}) \isDef \bigl\{x\in X\bigm| \norme{f_i(x)}\leq \norme{f_0(x)}, \, \forall i\bigr\}$, if $f_0\in A$ is such that $f_0,f_1,\dots,f_r$ have no common zeroes (or equivalently, they generate the unit ideal of~$A$).
\end{itemize}
\end{prop}

\begin{déf}
Let $X$ be an affinoid $K$-space.
We define the \emph{structural presheaf} of~$X$ as the presheaf $\mcl{O}_X$ which associates to any affinoid subdomain $U = \Sp B$ its corresponding affinoid algebra $\mcl{O}_X(U) = B$.
\end{déf}

The structural presheaf $\mcl{O}_X$ satisfies the sheaf property for finite coverings by affinoid subdomains, but it fails to be a sheaf for the canonical topology.
This is why we will have to restrict the kind of coverings we are allowed to consider.

\begin{déf}
Let $X$ be an affinoid $K$-space.
A subset $U\sset X$ is called an \emph{admissible open} if there are (possibly infinitely many) affinoid subdomains $(U_i)_I$ of $X$ such that $U = \bigcup_I U_i$ and satisfying the following property:
for every morphism of affinoid $K$-spaces $\vphi\colon Z\to X$ such that $\vphi(Z)\sset U$ the covering $(\vphi^{-1}(U_i))_i$ of $Z$ admits a finite refinement by affinoid subdomains.

A covering $(V_j)_j$ of an admissible open $V\sset X$ by means of admissible opens is called \emph{admissible} if for every morphism of affinoid $K$-spaces $\vphi\colon Z\to X$ such that $\vphi(Z)\sset V$ the covering $(\vphi^{-1}(V_j))_j$ of $Z$ admits a finite refinement by affinoid subdomains.

Admissible opens and admissible coverings define what is called a \emph{Grothendieck topology} on~$X$.
The presheaf $\mcl{O}_X$ is a sheaf for this topology.
\end{déf}

\begin{num}
A \emph{rigid $K$-variety} is a locally ringed space locally isomorphic to affinoid $K$-spaces (see \cite[5.3.4]{185} for a more precise definition).
We can glue rigid $K$-varieties and their morphisms as expected (see \cite[5.3.5 and 5.3.6]{185}).
A rigid $K$-variety is \emph{quasi-compact} if it can be covered by finitely many affinoid opens.
\end{num}

\begin{num}
We can adapt many classical properties of schemes to rigid $K$-varieties. We will use these notions freely and refer the reader to more comprehensive treatments such as~\cite{NAA} and \cite{185} for definitions.
\end{num}

\section{Rigid varieties from formal schemes and algebraic varieties}

\begin{num}
As hinted in~\ref{hintXeta}, rigid varieties arise naturally as generic fibres of formal $R$-schemes.
Let $R$ be a complete discrete valuation ring and $\mcl{X} = \Spf A$ an affine formal $R$-scheme of finite type. 
Then $A\otimes_R K$ is an affinoid $K$-space and we define the \emph{generic fibre} of~$\mcl{X}$ as the rigid $K$-variety $\Sp (A\otimes_R K)$.
Gluing yields a generic fibre~$\mcl{X}_K$ for any formal $R$-scheme topologically of finite type~$\mcl{X}$.
This construction is functorial, see details in~\cite[0.2.2]{BERT}.
\end{num}

\begin{déf}
We will say that a formal $R$-scheme topologically of finite type $\mcl{X}$ has \emph{pure relative dimension~$m$} if $\mcl{X}_K$ is equidimensional of dimension~$m$.
It follows from \cite[2.1.8]{48} that in this case, the special fibre $\mcl{X}_k$ has also pure relative dimension~$m$.
The converse holds if $\mcl{X}$ is flat over~$R$.
\end{déf}

\begin{num}
The correspondence $\mcl{X}\mapsto \mcl{X}_K$ even induces an equivalence of categories up to admissible blow-ups of formal $R$-schemes, i.e.\ formal blow-ups whose centre is supported on the special fibre (more details on formal blow-ups can be found in~\cite[2.15]{49}).
More precisely:
\end{num}

\begin{thm}[{\cite[4.1]{FRGI}}]\label{formRig:equivFormModel}
The functor $\mcl{X}\mapsto \mcl{X}_K$ induces an equivalence of categories between
\begin{itemize}
\item the category of flat formal $R$-schemes of topologically finite type, localized by admissible blow-ups;
\item the category of quasi-compact and quasi-separated rigid $K$-varieties.
\end{itemize}
Furthermore, under this correspondence $\mcl{X}_K$ is separated if and only if $\mcl{X}$ is separated and $\mcl{X}_K$ is smooth if $\mcl{X}$ is smooth.
\end{thm}

The following specific result will be extensively used for our constructions in chapter~\ref{ch:serre}.
\begin{prop}[{\cite[4.4]{FRGI}}] \label{formRig:blowUpModel}
Let $\mcl{X}$ be a flat quasi-compact formal $R$-scheme of finite type and $(U_l)_l$ a finite family of quasi-compact open subspaces of~$\mcl{X}_K$.
Then there is an admissible blow-up $h\colon \mcl{X}'\to\mcl{X}$ and a family $(\mcl{U}_l)_l$ of open formal subschemes of $\mcl{X}'$ such that $(\mcl{U}_l)_K = U_l$ for every~$l$.
Furthermore, if $(U_l)$ covers~$\mcl{X}_K$, then $(\mcl{U}_l)$ covers $\mcl{X}'$.
\end{prop}

\begin{num}
We can endow the set of closed points of an algebraic $K$-variety~$X$ with a structure of rigid $K$-variety. The resulting rigid $K$-variety~$X^\rig$ is called \emph{rigidification} of~$X$.
The rigid $K$-variety~$X^\rig$ is not quasi-compact in general, even if $X$ is, unless $X$ is proper (see \ref{compaRig}).
The construction yields a functor $X\mapsto X^\rig$ that preserves fibred products (see \cite[5.1.2]{99}).
\end{num}

\begin{prop}[{\cite[5.2.1]{99}}]
Let $f\colon X\to Y$ be a morphism between algebraic $K$-varieties.
The following properties hold for $f$ if and only if they hold for $f^\rig$: open immersion, smooth, separated, proper.
\end{prop}

\begin{prop}[{\cite[0.3.5]{BERT}}]\label{compaRig}
Let $\mcl{X}$ be an $R$-scheme of finite type.
There is a canonical morphism of rigid $K$-varieties $\alpha\colon (\what{\mcl{X}})_K\to (\mcl{X}_K)^\rig$.
It is an open immersion if $\mcl{X}$ is separated and an isomorphism if $\mcl{X}$ is proper over~$R$.
\end{prop}

\cleardoublepage


%

%
  \graphicspath{{\chaptersdir/serre/image/}}%
  \chapter{The motivic Serre invariant of a variety}\label{ch:serre}


\section{Introduction}
Let $R$ be a complete discrete valuation ring with fraction field $K$ and perfect residue field~$k$.
The motivic Serre invariant of a smooth separated quasi-compact rigid $K$-variety~$X$ was first defined by Loeser and Sebag in \cite[4.5.3]{98}. It gives a measure of the set of unramified points of $X$. Of particular interest is the fact that its nonvanishing implies the existence of an unramified point on~$X$.

This motivic invariant may be defined for a smooth $K$-variety $U$ as soon as $U$ is \emph{bounded}, i.e.\ if there exists a quasi-compact open rigid subvariety $V\sset U^\rig$ that contains all the unramified points of $U$. In that case, the motivic Serre invariant $S(V)$ of $V$ does not depend on the choice of $V$ and one sets $S(U) \isDef S(V)$. The boundedness condition holds for proper varieties, but fails for varieties as common as the affine space.

When $K$ has characteristic zero, the invariant has been extended additively by Nicaise in~\cite[5.4]{112} to arbitrary $K$-varieties.
In order to achieve this, Nicaise developed the theory of weak Néron models of pairs of varieties. The main drawback of his proof is its use of the weak factorization theorem, which is only known to hold when the field~$K$ has characteristic zero.

The goal of this chapter is to give a new construction of the motivic Serre invariant of a smooth $K$-variety $U$ avoiding the use of the weak factorization theorem. Theorem \ref{thmInvariant} shows that we are successful under the assumption that embedded resolution of singularities holds for $R$-schemes of finite type. In particular, we recover the result of Nicaise in equal characteristic zero.

The idea is quite natural. We consider a regular compactification $X$ of $U$ and approximate $U^\rig$ by a sequence of open quasi-compact rigid subvarieties $(U_\gamma)_\gamma$ by removing a tubular neighbourhood of the boundary.
We then show that the motivic Serre invariant of the approximations stabilizes as $\gamma\to\infty$, which allows us to define $S(U)$ as the limit of the sequence.
All this is carried out in section~\ref{motSerre}, with the main theorem stated in~\ref{thmInvariant}.

The motivic Serre invariant has already proved useful in answering Serre's question on the fixed locus of a finite $\ell$-group action on the affine space. Namely, Esnault and Nicaise have shown in \cite{120} that for a finite $\ell$-group $G$ acting on $\mbb{A}^n_K$ with $K$ a Henselian discretely valued field of characteristc zero with algebraically closed residue field of characteristic different from the prime number $\ell$, the action always admits a $K$-rational fixed point.
Their arguments build up from the additive extension of the motivic Serre invariant of Nicaise.
We show at the end of section~\ref{motSerre} how to adapt their arguments to fit our construction.
This yields in particular a conditional affirmative answer to Serre's question in equal characteristic.

Finally, we use in section~\ref{serre:nearbyCycles} a similar approach to give a new interpretation of the motivic nearby cycles with support of Guibert, Loeser and Merle \cite{53}. This application has been suggested by the previous work of Nicaise and Sebag in \cite[9.13]{106}, in which they give a description of the motivic nearby cycles of Denef and Loeser by means of the motivic volume of the so-called Milnor fibre, hence providing a setting similar to the one we work with in section~\ref{motSerre}.

\section{Preliminaries}
\subsection*{Notations and conventions}\label{notations}
\begin{num}
We fix a complete discrete valuation ring $R$ with fraction field $K$ and residue field~$k$. We will always assume that $k$ is perfect. We also fix a uniformizer $\pi$ of $R$ and a $\pi$-adic absolute value $\norme{\cdot}$ on $K$.
Recall that an extension $R\to R'$ is called \emph{unramified} if it has ramification index one (i.e.\ $\pi$ is a uniformizer of~$R'$) and the extension of residue fields $k'/k$ is separable.
If $R^\sh$ is a strict henselization of~$R$, then any unramified extension factors through~$R^\sh$ and any subextension of~$R$ in~$R^\sh$ is unramified.
We denote by $K^\sh$ the fraction field of~$R^\sh$.
For any perfect field extension $k'$ of $k$, we denote by $R\to\mcl{R}(k')$ the unique extension of~$R$ of ramification index one that induces the extension $k'/k$ on the residue fields. We will often write $R' = \mcl{R}(k')$ and $K' = \Frac R'$. 
\end{num}

\begin{num}
A variety over a Noetherian scheme $S$ is a separated $S$-scheme of finite type. A formal $R$-scheme will be called stft if it is separated and topologically of finite type over $R$.
For every $R$-scheme $\mcl{X}$ we denote by $\mcl{X}_k$ and $\mcl{X}_K$ its special and generic fibre, respectively. Likewise, for every stft formal $R$-scheme $\mcl{X}$ we write $\mcl{X}_k$ for its special fibre (which is a $k$-variety) and $\mcl{X}_K$ for its generic fibre (which is a quasi-compact separated rigid $K$-variety). In both cases, we say that  $\mcl{X}$ is \emph{generically smooth} if $\mcl{X}_K$ is smooth over~$K$. If $\mcl{X}$ is an $R$-scheme, we will denote by $\widehat{\mcl{X}}$ its formal $\pi$-adic completion. It is stft if $\mcl{X}$ is separated and of finite type over~$R$.
For every $K$-variety $X$, we denote by~$X^\rig$ its analytification in the category of rigid $K$-varieties.
If $X$ is an algebraic or rigid $K$-variety, then an \emph{unramified point} on~$X$ is a point in $X(K')$ for some finite extension of~$K$ in~$K^{\sh}$.
We will denote by $X(K^\sh)$ the set of unramified points of~$X$.
\end{num}

\begin{num}
If $X$ is a quasi-compact separated rigid $K$-variety, then a \emph{formal $R$-model} of~$X$ is an stft formal $R$-scheme $\mcl{X}$ together with an isomorphism of rigid $K$-varieties $\mcl{X}_K\to X$.
An \emph{admissible blow-up} of $\mcl{X}$ is a formal blow-up whose centre is supported in the special fibre $\mcl{X}_k$. Such an admissible blow-up does not affect the generic fibre $\mcl{X}_K$.
We recall that every quasi-compact separated rigid $K$-variety has a flat formal $R$-model (see~\ref{formRig:equivFormModel}).
\end{num}

\begin{num}
Let $U$ be a $K$-variety. By a \emph{compactification} of $U$ we mean a proper $K$-variety~$X$ endowed with an open immersion $U\to X$ such that $U$ is schematically dense in~$X$. Such a compactification always exists by Nagata's embedding theorem. If $U$ is regular, then we say that $X$ is an snc compactification if $X$ is regular and the boundary $X\prive U$ (with its reduced induced structure) is a divisor with strict normal crossings. Such snc compactifications exist when $K$ has characteristic zero, by Hironaka's embedded resolution of singularities.
\end{num}

\subsection*{Grothendieck rings of varieties}

Grothendieck rings of varieties serve as value rings for motivic integrals and are natural spaces for motivic invariants such as the motivic Serre invariant that we handle in this chapter.

\begin{déf}\label{serre:déf:Groth}
Let $S$ be a Noetherian scheme.
The \emph{Grothendieck ring of $S$-varieties} is defined as an abelian group by the following presentation:
\begin{itemize}
\item generators: isomorphism classes~$[X]$ of $S$-varieties;
\item relations: for any closed subvariety~$Y$ of~$X$, $[X] = [Y] + [X\prive Y]$.
\end{itemize}
We endow this group with a ring structure by putting $[X]\cdot[Y]\isDef [X\times_S Y]$ and denote it by $K_0(\mrm{Var}_S)$.
We also write $\mbb{L} = \mbb{L}_S$ for the class of the affine line $\mbb{A}^1_S$ and consider the localization
\[\mcl{M}_S\isDef K_0(\mrm{Var}_S)[\mbb{L}^{-1}].\]
\end{déf}

\begin{num}\label{serre:K0VarR}
The handling of motivic invariants over discrete valuation rings of mixed characteristic leads to the use of a slightly modified Grothendieck ring of varieties, as explained in \cite[2.4]{32}.
Therefore, we will set $\groth{S}\isDef K_0(\mathrm{Var}_S)$ if $R$ has equal characteristic, and denote by $\groth{S}$ the \emph{modified} Grothendieck ring of $S$-varieties if $R$ has mixed characteristic (see \cite[3.8.1]{36} for a definition; in the latter ring, we identify classes of universally homeomorphic $S$-varieties).
We also set $\mcl{M}^R_S \isDef \groth{S}[\mbb{L}^{-1}]$.
\end{num}

\begin{num}
Let $f\colon X\to S$ be a morphism of Noetherian schemes. It induces two morphisms between Grothendieck rings.
First, a ring morphism
\[f^*\colon K_0(\mrm{Var}_S)\to K_0(\mrm{Var}_X),\, [Y]\mapsto [Y\times_S X]\]
called \emph{base change morphism}; and then a \emph{group} morphism
\[f_!\colon K_0(\mrm{Var}_X)\to K_0(\mrm{Var}_S),\, [U]\mapsto [U],\]
called \emph{forgetful morphism}.
We emphasize that this latter morphism is not a ring morphism, as it sends $1 = [X]$ on~$[X]$.
But note that for any $S$-variety~$Y$ and any $X$-variety~$V$ we have in~$K_0(\mrm{Var}_S)$
\[f_!\bigl([f^*Y]\cdot [V] \bigr) = \bigl[(Y\times_S X)\times_X V\bigr] = [Y\times_S V] = [Y]\cdot [V].\]
In particular, $f_!\bigl(\mbb{L}_X \cdot [V]\bigr) = \mbb{L}_S \cdot [V]$.

These two morphisms factor through localization by~$\mbb{L}^{-1}$ and through the modified Grothendieck rings (see \cite[\S 3.8]{36}).
\end{num}

\begin{num}
In \cite[3.1]{86}, Bittner gave a presentation for $K_0(\mrm{Var}_F)$ in terms of smooth $F$-varieties when $F$ is a field of characteristic zero. Her proof is based on stratification of $F$-varieties into smooth varieties.
Proposition~\ref{regStrat} and corollary~\ref{refineRegStrat} below allow to deduce from Bittner's proof a presentation of $K_0(\mrm{Var}_F)$ in terms of connected regular varieties, with no assumption on the characteristic of~$F$ (this is corollary~\ref{presentationGroth}).
This presentation will be used in~\ref{extensionMK} to extend the motivic Serre invariant additively to arbitrary $K$-varieties.
\end{num}

\begin{déf}\label{déf:regStrat}
Let $F$ be any field and $X$ an $F$-variety. A \emph{regular stratification} of $X$ is a partition of $X$ into finitely many connected regular locally closed subvarieties $(N_\nu)_{\nu\in\mcl{N}}$ of~$X$, called \emph{strata}, such that the schematic closure of each~$N_\nu$ is a union of strata.
\end{déf}

\begin{prop}\label{regStrat}
Let $F$ be any field, $X$ a reduced $F$-variety and let $(V_\alpha)_{\alpha\in A}$ be finitely many closed subvarieties of $X$.
Then there is a regular stratification $(N_\nu)_\nu$ of~$X$ such that each $V_\alpha$ is a union of strata.
\end{prop}
\begin{proof}
We first notice that every irreducible component $X_i$ of~$X$ is a union of strata in any regular stratification.
Indeed, the schematic closure of a stratum~$N_\nu$ containing the generic point of $X_i$ is necessary equal to~$X_i$ because $N_\nu$ is irreducible.
We may therefore assume that neither $X$ nor any irreducible component of~$X$ belongs to $(V_\alpha)$.

We proceed by induction on $\dim X$. If $\dim X = 0$, this is clear.
Assume now $\dim X\geq 1$ and let $Y = X\prive \bigcup_A V_\alpha$.
Denote by $(X_i)_{1\leq i\leq n}$ the irreducible components of~$X$. Thus $Y_i\isDef Y\cap X_i$ are the irreducible components of $Y$.
We set $X_i^\circ = X_i\prive \bigcup_{j\neq i} X_j$, which is open in $X$, and similarly $Y_i^\circ = Y\cap X_i^\circ$.

Each $Y_i^\circ$ admits a connected regular open dense subset $U_i$ and $\dim Y_i^\circ\prive U_i < \dim X$. Hence we can find regular stratifications $(N_{i\mu})_{\mu\in M_i}$ of~$Y_i^\circ\prive U_i$ for every $1\leq i\leq n$.
We have in $X$
\[\ovl{N_{i\mu}} = (\ovl{N_{i\mu}}\cap Y_i^\circ) \sqcup \bigl(\ovl{N_{i\mu}}\cap (X\prive Y_i^\circ)\bigr).\]
Put $X'\isDef X\prive \bigsqcup_i Y_i^\circ$. We have $\bigcup_A V_\alpha\sset X'$ and since $\ovl{N_{i\mu}}\sset X_i$, we also have $\ovl{N_{i\mu}}\cap X' = \ovl{N_{i\mu}}\cap (X\prive Y_i^\circ)$.

We may now apply the induction hypothesis on $X'$ with the family consisting of the closed subvarieties $(V_\alpha)_\alpha$ and $(\ovl{N_{i\mu}}\cap X')_{M_i}$ for every $1\leq i\leq n$ to obtain a regular stratification $(N_\lambda)_\lambda$ of~$X'$.

To sum up, we have
\[X = \bigsqcup_i U_i \sqcup \bigsqcup_i \bigsqcup_{M_i} N_{i\mu} \sqcup \bigsqcup_\lambda N_\lambda.\]

\begin{itemize}
\item Since $X'$ is closed in $X$, each $\ovl{N_\lambda}$ is a union of strata of $(N_\lambda)_\lambda$.
\item Since $\ovl{N_{i\mu}}\cap Y_i^\circ = \ovl{N_{i\mu}}\cap (Y_i^\circ\prive U_i)$, they are unions of strata of $(N_{i\mu})_{M_i}$.
\item Each $\ovl{N_{i\mu}}\cap (X\prive Y_i^\circ)$ is a union of strata of $(N_\lambda)_\lambda$ by construction.
\item $\ovl{U_i} = X_i = Y_i^\circ \sqcup (X'\cap X_i)$ and $X'\cap X_i$ is a union of irreducible components of $X'$, hence a union of strata of $(N_\lambda)_\lambda$.
\item Each $V_\alpha$ is a union of strata of $(N_\lambda)_\lambda$ by construction.\qedhere
\end{itemize}
\end{proof}

\begin{cor}\label{refineRegStrat}
Two regular stratifications of an $F$-variety $X$ admit a common refinement.
\end{cor}
\begin{proof}
If $(N_\nu)_\nu$ and $(N_\lambda)_\lambda$ are two regular stratifications of $X$, we apply \ref{regStrat} on $X$ with the family of closed subvarieties consisting of every $\ovl{N}_\nu$ and $\ovl{N}_\lambda$.
\end{proof}

\begin{cor} \label{presentationGroth}
Let $F$ be any field.
As an abelian group, $K_0(\mathrm{Var}_F)$ is given by the following presentation:
\begin{itemize}
\item {\em generators:} isomorphism classes $[U]$ of connected regular $F$-varieties $U$;
\item {\em relations:} $[U]=[V]+[U\prive V]$ whenever $U$ is a connected regular $F$-variety and $V$ is a connected regular closed subvariety of $U$.
\end{itemize}
\end{cor}

\subsection*{Weak Néron models and the motivic Serre invariant}
\begin{déf}\label{serre:déf:weakNéron}
Let $X$ be a $K$-variety. A \emph{weak Néron model} for $X$ is a pair $(\mcl{U},\iota)$ where
\begin{enumerate}
\item $\mcl{U}$ is a smooth separated $R$-scheme of finite type,
\item $\iota\colon\mcl{U}_K\to X$ is an open immersion,
\item for every finite unramified extension $R'$ of $R$ with fraction field $K'$, the map $\mcl{U}(R') \to X(K')$ induced by $\iota$ is a bijection.
\end{enumerate}
The analogous definitions apply when $X$ is a separated rigid $K$-variety and $\mcl{U}$, $\mcl{X}$ and $\mcl{Y}$ are stft formal $R$-schemes.
\end{déf}

\begin{num}
Beware that, in the algebraic case, our definition of a weak Néron model is less restrictive than the one used in \cite[3.5.1]{NER}, where it is required that $\iota$ is an isomorphism.
It follows easily from the boundedness criterion~\cite[3.5.7]{NER} and~\cite[3.4.2]{NER} that an algebraic $K$-variety $X$ has a weak Néron model if and only if the set of unramified points on $X$ is contained in the $K$-smooth locus~$X_\sm$ of $X$ and bounded in~$X$ in the sense of \cite[1.1.2]{NER} (or equivalently, bounded in~$X_\sm$ by~\cite[1.1.9]{NER}).
Likewise, a separated rigid $K$-variety $Y$ admits a formal weak Néron model if and only if there exists a smooth quasi-compact open rigid subvariety $U$ of~$Y$ such that the map $U(K')\to Y(K')$ is a bijection for every finite unramified extension $K'$ of $K$ (see~\cite[4.8]{112}). Now one deduces easily from the observations in \cite[\S4]{112} that an algebraic $K$-variety $X$ admits a weak Néron model $\mcl{U}$ if and only if its rigid analytification $X^{\rig}$ admits a formal weak Néron model; in that case, $\widehat{\mcl{U}}$ is a formal weak Néron model of~$X^{\rig}$.
\end{num}

\begin{prop}\label{prop-regwner}
If $X$ is a regular and proper $K$-variety, then $X$ has a weak Néron model.
\end{prop}
\begin{proof}
It follows from \cite[17.15.1]{EGA-IV4} that $X$ is smooth over $K$ at each separable point, and thus at each unramified point. The set of unramified points on $X$ is bounded in $X$ because $X$ is proper \cite[1.1.6]{NER}. Thus $X$ has a weak Néron model.
\end{proof}

Weak Néron models are far from being unique, in general. Nevertheless, the theory of motivic integration yields the following remarkable result.

\begin{thm}\label{serre:SXoriginel}
Let $X$ be a separated rigid $K$-variety which admits a weak Néron model $\mcl{U}$.
Then the element
\[S(X) \isDef [\mcl{U}_k] \in \grothModL\]
only depends on $X$, and not on the choice of a weak Néron model. It is called
the \emph{motivic Serre invariant} of $X$.
\end{thm}

\begin{num}\label{SerreAlgébrique}
If $X$ is smooth and quasi-compact, this result is due to Loeser and Sebag \cite[4.5.3]{98}. As was already observed in \cite[5.1]{112}, the general case follows from the proof of~\cite[5.9]{145}.
If $X$ is an algebraic $K$-variety with a weak Néron model $\mcl{U}$, then $\widehat{\mcl{U}}$ is a formal weak Néron model of $X^{\rig}$ so that
\[S(X^{\rig}) = [\mcl{U}_k] \in \grothModL.\] We denote this invariant by $S(X)$ and call it the \emph{motivic Serre invariant} of~$X$.
\end{num}

\subsection*{Technical assumptions}

To prove the main results of this chapter, we need to formulate two technical assumptions, both of which would follow from the existence of embedded resolutions of singularities.

\begin{num}\label{hyp:snc} {\em Existence of snc compactifications.} We will need to assume that, for every regular $K$-variety $U$ and every regular closed subvariety $V$ of $U$, each compactification of $U$ can be dominated by an snc
compactification $X$ of $U$ such that the boundary $X\prive U$ (with its reduced induced structure) is a strict normal crossings divisor that has transversal intersections with the schematic closure of $V$ in $X$. This property is known if $K$ has characteristic zero, by Hironaka's embedded resolution of singularities for $K$-varieties, but it is an open problem if $K$ has positive characteristic and the dimension of $U$ is at least three (recent work of Cossart and Piltant seems to have settled the three-dimensional case, see remark~\ref{rmq:cossart}).
\end{num}

\begin{num}\label{hyp:adapt} {\em Existence of weak Néron models for snc pairs.} Let $X$ be a regular proper $K$-variety and let $\partial$ be a divisor with strict normal crossings on $X$.
We will need to assume that every proper $R$-model of $X$ can be dominated by a weak Néron model $\mcl{X}$ of $X$ such that the schematic closure $\ovl{\partial}$ of $\partial\cap \mcl{X}_K$ in $\mcl{X}$ is a relative strict normal crossings divisor over~$R$; since $\mcl{X}$ is smooth over $R$, this is equivalent to saying that the union of $\ovl{\partial}$ with the special fibre $\mcl{X}_k$ is a divisor with strict normal crossings on $\mcl{X}$.
Note that this implies that every prime component $E$ of $\ovl{\partial}$ is smooth at every point of its special fibre, since $\pi$ is part of a regular system of local parameters on $E$ at such a point and $k$ is assumed to be perfect.
The existence of $\mcl{X}$ is known if $k$ has characteristic zero, by Hironaka's embedded resolution of singularities for $R$-varieties: we can take for $\mcl{X}$ the $R$-smooth locus of a regular proper $R$-model such that the union of the special fibre and the schematic closure of $\partial$ has strict normal crossings.
It is also known if $X$ is smooth over $K$ and $\partial$ is a smooth prime divisor, by the existence of weak Néron models of pairs proven by Nicaise in~\cite[3.15]{112}.
It is not hard to generalize the proof of the latter result to the case where $X$ and $\partial$ are regular instead of smooth, but it seems more difficult to extend it to the case where $\partial$ has several components.
\end{num}

\section{Construction of the motivic Serre invariant of a variety}
\label{motSerre}

In this section we fix a smooth algebraic $K$-variety~$U$.
In general, such a variety does not admit a weak Néron model, which prevents us from defining its motivic Serre invariant in a straightforward manner.
The goal of this section is to provide a construction leading to such a definition.

\subsection*{Main theorem}

\begin{déf}
A \emph{tubular datum} for $U$ is a couple $\approche$ such that $X$ is a normal compactification of $U$, $\{W_{\alpha}\}_{\alpha\in A}$ is a finite admissible cover of~$X^{\rig}$ by quasi-compact open rigid subvarieties, and $f_{\alpha}$ is an analytic function on~$W_{\alpha}$, for each $\alpha$ in $A$.
We require moreover that there exists a Cartier divisor $\partial$ on~$X$ supported on $X\prive U$ such that the restriction to $W_{\alpha}$ of the divisor induced by $\partial$ on~$X^\rig$ is defined by $f_{\alpha}=0$ for every $\alpha$ in $A$.
\end{déf}

\begin{num}
Note that every smooth $K$-variety $U$ admits a tubular datum, since we can start with any compactification $Y$ of $U$ and consider the normalized blow-up of $Y$ at the reduced boundary $Y\prive U$.
\end{num}

\begin{num}\label{constructionUgamma}
Let $\approche$ be a tubular datum for $U$. For every nonnegative integer $\gamma$, we set
\[U_\gamma \isDef \bigcup_{\alpha\in A} \Bigl\{x\in W_\alpha \Bigm| \norme{f_\alpha(x)} \geq \norme{\pi}^\gamma\Bigr\}.\]
Each of the sets
\[W_{\alpha}\Bigl(\frac{\pi^\gamma}{f_\alpha}\Bigr) = \Bigl\{x\in W_{\alpha} \Bigm| \norme{f_\alpha(x)} \geq \norme{\pi}^\gamma\Bigr\}\]
is a quasi-compact open rigid subvariety of $W_{\alpha}$.
It follows that $U_\gamma$ is a quasi-compact open rigid subvariety of the smooth separated rigid variety $U^{\rig}$ and therefore admits a weak Néron model.
It is clear from the definition that the sequence $(U_\gamma)_{\gamma\geq 0}$ is increasing and that
\[U^\rig = \bigcup_{\gamma\geq 0}U_{\gamma}=\bigcup_{\alpha\in A} \Bigl\{x\in W_\alpha \Bigm| \norme{f_\alpha(x)} > 0\Bigr\}.\]
We will call $(U_\gamma)_{\gamma\geq 0}$ the \emph{approximating sequence} associated to $\approche$.
\end{num}

\begin{thm}\label{thmInvariant}
Assume that the properties \eqref{hyp:snc} and \eqref{hyp:adapt} hold (this is known if $k$ has characteristic zero).
Let $U$ be a smooth $K$-variety. Let $\mcl{T} = \approche$ be a tubular datum of~$U$ and let $(U_\gamma)_\gamma$ be the associated approximating sequence.
Then there exists a positive integer $\gamma_0$ such that for every integer $\gamma\geq\gamma_0$ we have $S(U_\gamma) = S(U_{\gamma_0})$ in $\groth{k}/(\mbb{L}-1)$.  Furthermore, this value $S(U_{\gamma_0})$ only depends on $U$ and not on $\mcl{T}$.
\end{thm}

We postpone the proof to \ref{proofMainThm}.

\begin{déf}\label{def:motserre}
With the notations of \ref{thmInvariant}, we define the \emph{motivic Serre invariant}~$S(U)$ of~$U$ as $S(U_{\gamma_0})$.
\end{déf}

\begin{num} If $U(K^\sh) = \emptyset$, then for any tubular datum $\mcl{T}$ we have $U_\gamma(K^\sh) = \emptyset$ for all $\gamma\geq 0$ and $\emptyset$ is a weak Néron model for each $U_\gamma$.
Hence $S(U_\gamma) = 0$ for every~$\gamma$ and $S(U) = 0$.
Thus the nonvanishing of $S(U)$ is a sufficient condition for the existence of an unramified point on $U$.
\end{num}

\begin{eg}
Consider the variety $U = \mbb{A}^1_K = \Spec K[x_0]$ and the tubular datum $\bigl(\mbb{P}^1_K, (W_\alpha,f_\alpha)_{i=1,2}\bigr)$ given by \begin{eqnarray*}
(W_1,f_1)&=& \bigl(\Sp K\{\tfrac{x_0}{x_1}\},1\bigr),
\\(W_2,f_2)&=& \bigl(\Sp K\{\tfrac{x_1}{x_0}\},\tfrac{x_1}{x_0}\bigr),
\end{eqnarray*} where $(x_0:x_1)$ are homogeneous coordinates on~$\mbb{P}^1_K$.
Then for every $\gamma\geq 0$, we have $W_1(\frac{\pi^\gamma}{1}) = W_1$ and
\[W_2\Bigl(\frac{\pi^\gamma}{x_1/x_0}\Bigr) = \Bigl\{x\in W_{2} \Bigm| \norme{\pi^\gamma x_0/x_1} \leq 1\Bigr\},\]
thus $U_\gamma$ is the closed disc $\Sp K\{\pi^\gamma x_0/x_1\}$ of radius $\norme{\pi}^{-\gamma}$ around the origin in~$(\mbb{A}^1_K)^{\rig}$. Since the formal affine line $\Spf R\{\pi^\gamma x_0/x_1\}$ is a weak Néron model for $U_\gamma$, we find that $S(U_\gamma)=1$  in $\groth{k}/(\mbb{L}-1)$ in for every $\gamma\geq 0$. Thus $S(\mbb{A}^1_K)=1$.
\end{eg}

The following proposition confirms that the motivic Serre invariant in definition \ref{def:motserre} agrees with the definition given in~\ref{SerreAlgébrique} when~$U$ admits a weak Néron model.

\begin{prop}\label{prop:defconsist}
Let $U$ be a smooth $K$-variety that admits a weak Néron model~$\mcl{U}$. Then the motivic Serre invariant of~$U$ as defined in~\ref{def:motserre} equals the class $[\mcl{U}_k]$ in $\groth{k}/(\mbb{L}-1)$.
\end{prop}
\begin{proof}
Let $\mcl{T}$ be a tubular datum for~$U$ and let $(U_\gamma)_\gamma$ be the associated approximating sequence.
Since $(\wht{\mcl{U}})_K$ is a quasi-compact subvariety of~$U^\rig$, there exists $\gamma_0\geq 0$ such that
\[(\wht{\mcl{U}})_K\sset U_{\gamma}\sset U^{\rig}\]
for every $\gamma\geq \gamma_0$.
But for every finite unramified extension $K'/K$, we then have
\[U^{\rig}(K') = (\wht{\mcl{U}})_K(K')\sset U_\gamma(K')\sset U^\rig(K'),\]
so that $\wht{\mcl{U}}$ is a weak Néron model of~$U_\gamma$ and $S(U_\gamma)=[\mcl{U}_k]$ for all $\gamma\geq \gamma_0$.
\end{proof}

\subsection*{Greenberg schemes and the motivic measure}
In order to prove theorem \ref{thmInvariant}, we interpret the motivic Serre invariants $S(U_\gamma)$ in terms of motivic measures of cylinders in the Greenberg scheme of a weak Néron model of~$X$. We first recall some important definitions, referring to \cite[2.2]{32} and  \cite{48} for further background.

\begin{num}
Let $\mcl{X}$ be an stft formal $R$-scheme.
For every integer $n\geq 0$, we denote by $\funcGr_n(\mcl{X})$ the Greenberg scheme of level $n$ of $\mcl{X}$.
This is a separated $k$-scheme of finite type and there exists for every perfect field extension $k'$ of $k$ a canonical bijection
\[\funcGr_n(\mcl{X})(k')=\mcl{X}(\mcl{R}(k')/(\pi^{n+1})).\]
For $n=0$, we have $\funcGr_0(\mcl{X})=\mcl{X}_k$.
The Greenberg schemes of different levels fit into a projective system with affine transition morphisms $$\theta^m_n\colon\funcGr_m(\mcl{X})\to \funcGr_n(\mcl{X})$$ for $m\geq n$ which correspond to reduction modulo $\pi^{n+1}$ on the level of $k'$-points.
The projective limit of this system is a $k$-scheme $\funcGr(\mcl{X})$ (of infinite type, in general) endowed with projection morphisms
\[\theta_n\colon\funcGr(\mcl{X})\to \funcGr_n(\mcl{X}).\]
This scheme is called the \emph{Greenberg scheme} of~$\mcl{X}$.
It parameterizes unramified points on~$\mcl{X}$: for every perfect field extension $k'$ of $k$, there is a canonical bijection
\[\funcGr(\mcl{X})(k')=\mcl{X}(\mcl{R}(k')).\]
\end{num}

\begin{num}\label{GrnAd}
The Greenberg scheme of level~$n$ of the formal affine space $\wht{\mbb{A}}^d_R$ is quite easy to describe: it is (non-canonically) isomorphic to $\mbb{A}^{d(n+1)}_k$ (see \cite[2.2.3]{32}).
We will be able to reduce to this case for computations using the following proposition.
\end{num}

\begin{prop}[{\cite[2.2.4]{32}}]\label{serre:GrEtale}
Let $\mcl{X}\to \mcl{Y}$ be an étale morphism of stft formal $R$-schemes.
Then for every $m\geq n \geq 0$ the canonical diagram
\begin{center}
\begin{tikzpicture}
    \matrix (m) [matrix of math nodes, row sep=3em,
    column sep=3em, text height=1.5ex, text depth=0.25ex]
    { \funcGr_m(\mcl{X}) & \funcGr_m(\mcl{Y}) \\
\funcGr_n(\mcl{X}) & \funcGr_n(\mcl{Y}) \\};
    \path[arrow]
    (m-1-1) edge node[auto] {$$} (m-1-2)
            edge node[auto,swap] {$ \theta^m_n $} (m-2-1)
    (m-1-2) edge node[auto] {$ \theta^m_n $} (m-2-2)
(m-2-1) edge node[auto,swap] {$$} (m-2-2);
\draw ($(m-1-1.east)!.07!(m-1-2.west)$)+(0,-.12) coordinate (A);
\draw ($(m-1-1.south)!.25!(m-2-1.north)$)+(.12,0) -| (A);

\end{tikzpicture}
\end{center}
is cartesian.
\end{prop}



\begin{num}
A \emph{cylinder} in $\funcGr(\mcl{X})$ is a subset of the form $\mcl{C}=(\theta_n)^{-1}(C_n)$ for some $n\geq 0$ and some constructible subset $C_n$ of $\funcGr_n(\mcl{X})$.
If $\mcl{X}$ is smooth over $R$ of pure relative dimension $d$, then the morphism $\theta^m_n$ is a locally trivial fibration with fibre $\mbb{A}^{d(m-n)}_k$, for every $m\geq n$.
It follows that the element
\[\mu_{\mcl{X}}(C)\isDef\, [C_n]\,\mbb{L}^{-dn} \in \mcl{M}_{\mcl{X}_k}\] does not depend on the choice of $n$.
It is called the \emph{motivic measure} of the cylinder~$C$.
\end{num}

The following result on cylinders will be needed in~\ref{prop:cylinder}. We prove it here by lack of a proper reference.
\begin{prop}\label{cylindrePtFermés}
A cylinder in $\funcGr(\mcl{X})$ is fully determined by its sets of points over finite extensions of~$k$.
\end{prop}
\begin{proof}
Let $\mcl{C} = \theta_n^{-1}(C_n)$ and $\mcl{D} = \theta_m^{-1}(D_m)$ be two cylinders such that $\mcl{C}(k') = \mcl{D}(k')$ for every finite extension $k'/k$.
We may assume by \cite[4.3.9]{48} that $m=n$, $C_n = \theta_n(\mcl{C})$ and $D_n = \theta_n(\mcl{D})$.

By \cite[10.4.7]{EGA-IV3}, the set of closed points of $\funcGr_n(\mcl{X})$ is very dense in $\funcGr_n(\mcl{X})$. Hence any constructible subset of $\funcGr_n(\mcl{X})$ is fully determined by the set of its closed points (see \cite[10.1.2(c')]{EGA-IV3}).

Now a closed point $x\in C_n$ induces a point $x\in C_n(k')$ for some finite extension of~$k$.
By Greenberg's approximation theorem, there is some $l\geq n$ such that $\theta_n(\funcGr(\mcl{X})) = \theta_n^l(\funcGr_l(\mcl{X}))$. Hence the fibre of $\theta_n^l$ is nonempty over~$x$. 
Since $\funcGr_l(\mcl{X})$ is of finite type over~$k$, it contains a $k''$-point for some finite extension $k''/k'$, and this point lifts to a point in $\mcl{C}(k'')$.
In particular, the set of closed points of~$C_n$ is fully determined by the sets of points over finite extensions of~$k$ of~$\mcl{C}$.
It follows at once from the assumptions that $C_n=D_n$, thus $\mcl{C} = \mcl{D}$.
\end{proof}

\begin{num}
Let $\mcl{I}$ be a coherent ideal sheaf on $\mcl{X}$. If $k'$ is a perfect field extension of $k$ and $\varphi$ is a point in $\mcl{X}(\mcl{R}(k'))$, then we set
$$\ord_\mcl{I}(\varphi) \isDef \min \Bigl\{v\bigl(\varphi^*(f)\bigr) \Bigm| f\in \mcl{I}_{\varphi(s)}\Bigr\}\in \mbb{N}\cup \{+\infty\}$$
where $s$ denotes the closed point of $\Spec \mcl{R}(k')$ and $v$ is the normalized discrete valuation on $\mcl{R}(k')$. In this way, we obtain a function
\[\ord_\mcl{I}\colon \funcGr(\mcl{X}) \to \mbb{N}\cup\{\infty\}\]
and it is easy to see that its fibres over elements of $\mbb{N}$ are cylinders.
If $D$ is the closed subscheme of $\mcl{X}$ defined by $\mcl{I}$, we will also denote this function by $\ord_D$.
\end{num}

\subsection*{Proof of the theorem}\label{subsec-proof}

\begin{prop}\label{prop:cylinder}
Let $\approche$ be a tubular datum for $U$ and let $\mcl{X}$ be a weak Néron model for $X$.
Then for every integer $\gamma\geq 0$, there exists a unique cylinder~$\mcl{C}_\gamma$ in the Greenberg scheme $\funcGr(\wht{\mcl{X}})$ with the following property:
if $k'$ is a finite field extension of $k$, $R'=\mcl{R}(k')$ is the associated unramified extension of~$R$ and $K'$ is the quotient field of~$R'$, then the $k'$-points in $\mathcal{C}_\gamma$ correspond precisely to the $K'$-points of $U_{\gamma}$ under the natural bijections
\[\funcGr(\wht{\mcl{X}})(k')=\wht{\mcl{X}}(R')=X^\rig(K').\]
Moreover,
\[S(U_{\gamma})=\mu(\mcl{C}_\gamma)\] in $\groth{k}/(\mbb{L}-1)$.
\end{prop}
\begin{proof}
The cylinder $\mathcal{C}_\gamma$ is unique by~\ref{cylindrePtFermés}.
To prove its existence and the expression for $S(U_\gamma)$, we can argue as follows.
By \ref{formRig:blowUpModel} and \cite[3.1]{41}, we can find a weak Néron model $\mcl{U}_{\gamma}$ of~$U_\gamma$ dominating $\wht{\mcl{X}}$ and such that the induced morphism of rigid varieties $(\mcl{U}_\gamma)_K\to (\wht{\mcl{X}})_K$ corresponds to the inclusion $U_\gamma\cap (\wht{\mcl{X}})_K\to (\wht{\mcl{X}})_K$.
Then we define $\mathcal{C}_\gamma$ to be the image of the injective morphism $\funcGr(\mcl{U}_\gamma)\to \funcGr(\wht{\mcl{X}})$.
Thus, for every finite extension $k'$ of~$k$, we have
\[\mathcal{C}_\gamma(k')=\funcGr(\mcl{U}_\gamma)(k')=U_\gamma(K').\]
It follows from \cite[2.4.4]{32} that $\mathcal{C}_\gamma$ is a cylinder and that
$$\mu_{\mcl{X}}(\mathcal{C}_\gamma)=\mu_{\mcl{U}_\gamma}(\funcGr(\mcl{U}_\gamma))$$ in $\groth{k}/(\mbb{L}-1)$.
On the other hand, we have
$$ \mu_{\mcl{U}_\gamma}(\funcGr(\mcl{U}_\gamma))=[(\mcl{U}_\gamma)_k]=S(U_\gamma)$$ in $\groth{k}/(\mbb{L}-1)$
because $\mcl{U}_\gamma$ is a weak Néron model for $U_\gamma$.
\end{proof}
Thus, to prove that $S(U_\gamma)$ stabilizes for large $\gamma$, we need to understand how the motivic measure of $\mathcal{C}_\gamma$ depends on $\gamma$. We will first study this question for a particular class of tubular data.

\begin{déf}\label{def:adapted}
A tubular datum $\approche$ for $U$ is said to be \emph{snc} if $X$ is regular and the Cartier divisor $\partial$ supported on $X\prive U$ has strict normal crossings.
A weak Néron model $\mcl{X}$ of~$X$ is called \emph{adapted} to an snc tubular datum $\approche$ if the following conditions are satisfied:
\begin{enumerate}
\item the schematic closure of $\partial\cap \mcl{X}_K$ in $\mcl{X}$ is a relative strict normal crossings divisor over $R$;
\item \label{it:adap2} for each $\alpha$ in $A$ there exists a formal open subscheme $\mcl{W}_{\alpha}$ in $\wht{\mcl{X}}$ such that $$(\mcl{W}_{\alpha})_K  =W_\alpha\cap (\wht{\mcl{X}})_K.$$
\end{enumerate}
\end{déf}
We can find an snc tubular datum for $U$ by our assumption \eqref{hyp:snc}.
The existence of adapted weak Néron models follows from our assumption \eqref{hyp:adapt}: starting from any proper $R$-model of $X$, we can ensure that \eqref{it:adap2} holds by
taking a suitable admissible blow-up (see \ref{formRig:blowUpModel}), and then dominate the resulting model by a weak Néron model as in~\eqref{hyp:adapt}.

\begin{prop}\label{prop:compu}
Let $\approche$ be an snc tubular datum for $U$ and let $\mcl{X}$ be an adapted weak Néron model. Denote by $\mcl{Y}$ the complement in $\mcl{X}$ of the Zariski closure of $X\prive U$.
Then
there exists an integer $\gamma_0\geq 0$ such that
$$S(U_\gamma)=[\mcl{Y}_k]$$ in $\groth{k}/(\mbb{L}-1)$ for all $\gamma\geq \gamma_0$.
\end{prop}
\begin{proof}
We define the cylinders $\mathcal{C}_\gamma$ as in proposition~\ref{prop:cylinder}; it suffices to show that  $$\mu_{\mcl{X}}(\mathcal{C}_\gamma)=[\mcl{Y}_k]$$ in $\groth{k}/(\mbb{L}-1)$ when $\gamma$ is sufficiently large.
Denote by $\partial$ the Cartier divisor on $X$ supported on $X\prive U$ and defined by $f_{\alpha}=0$ on each $W_{\alpha}$.
Let $\ovl{\partial}$ be the schematic closure of $\partial\cap\mcl{X}_K$ in $\mcl{X}$; this is a relative strict normal crossings divisor over $R$, which we write as
\[\ovl{\partial}=N_1D_1+\ldots+N_tD_t.\]
Here $D_1,\ldots,D_t$ are the irreducible components of $\ovl{\partial}$ and $N_1,\ldots,N_t$ are their multiplicities.
For every element $\delta$ of $\mbb{N}^t$, we consider the cylinder
\[\mcl{D}_\delta \isDef \Bigl\{x\in \funcGr(\wht{\mcl{X}})\Bigm|\mrm{ord}_{D_i}(x)=\delta_i\text{ for }i=1,\ldots,t \Bigr\}\]
in $\funcGr(\wht{\mcl{X}})$.

We take formal open subschemes $\mcl{W}_\alpha$ of $\wht{\mcl{X}}$ as in definition \ref{def:adapted}, and we denote by $\mcl{W}_{\alpha,j}$, $j=1,\ldots,c_{\alpha}$ the connected components of $\mcl{W}_{\alpha}$.
For every $j$ we denote by $m_{\alpha,j}\in \mbb{Z}$ the order of $f_{\alpha}$ along the special fibre of $\mcl{W}_{\alpha,j}$.
This means that $f_\alpha = \pi^{m_{\alpha,j}}g_\alpha$ for some section $g_\alpha\in \mcl{O}(\mcl{W}_{\alpha,j})$ defining $\ovl{\partial}$ on~$\mcl{W}_{\alpha,j}$.
If $k'$ is a finite extension of $k$, $R'$ is the associated unramified extension of $R$, $K'$ is the quotient field of $R'$ and $x$ is a point of
\[\funcGr(\mcl{W}_{\alpha,j})(k')\sset W_{\alpha}(K'),\]
then
\[v(f_{\alpha}(x)) = m_{\alpha,j} + \sum_{i=1}^t N_i\, \mrm{ord}_{D_i}(x),\]
where $v$ denotes the $\pi$-adic valuation on $K'$.

We write $A'$ for the set of couples $(\alpha,j)$ with $\alpha$ in $A$ and $1\leq j\leq c_{\alpha}$.
For every nonempty subset $B$ of $A'$, we set
$$\mcl{W}^\circ_{B,s} \isDef \bigcap_{(\alpha,j)\in B}(\mcl{W}_{\alpha,j})_k\,\Big\backslash\, \Bigl(\bigcup_{(\alpha,j)\notin B}(\mcl{W}_{\alpha,j})_k\Bigr).$$
As $B$ runs over the nonempty subsets of $A'$, the sets $\mcl{W}^\circ_{B,s}$ form a partition of~$\mcl{X}_k$ into locally closed subsets.
For every $\delta$ in $\mbb{N}^t$, we set
$$\mathcal{D}_{\delta,B} \isDef \mathcal{D}_\delta\cap (\theta_0)^{-1}(\mcl{W}^\circ_{B,s}).$$
This is still a cylinder in $\funcGr(\wht{\mcl{X}})$.
Note that
$$\mathcal{D}_{\mathbf{0},B}=(\theta_0)^{-1}(\mcl{W}^\circ_{B,s}\cap \mcl{Y}_k),$$
where we write  $\mathbf{0}$ for the zero vector in $\mbb{N}^t$

It is clear from the definition of the cylinder $\mcl{C}_\gamma$ that, for every nonempty subset~$B$ of~$A'$, we can write $$\mathcal{C}_\gamma\cap (\theta_0)^{-1}(\mcl{W}^\circ_{B,s})\sset \funcGr(\wht{\mcl{X}})$$ as a finite disjoint union
$$
\bigcup_{\delta}\mathcal{D}_{\delta,B}
$$ where $\delta$ runs over the set of elements in $\mbb{N}^t$ such that
$$m_{\alpha,j}+N_1\delta_1+\ldots+N_t\delta_t \leq \gamma$$
for at least one element $(\alpha,j)$ in $B$.
We will prove below that the motivic measure of  $\mathcal{D}_{\delta,B}$ is zero in $\groth{k}/(\mbb{L}-1)$ unless $\delta$ is the zero vector $\mathbf{0}$.
Assuming this claim for now, we find that
\[\mu_{\mcl{X}}\Bigl(\mathcal{C}_\gamma\cap (\theta_0)^{-1}(\mcl{W}^\circ_{B,s})\Bigr) = \mu_{\mcl{X}}\Bigl((\theta_0)^{-1}(\mcl{W}^\circ_{B,s}\cap \mcl{Y}_k)\Bigr)\]
in $\groth{k}/(\mbb{L}-1)$ if $\gamma\geq m_{\alpha,j}$ for some $(\alpha,j)$ in $B$, and
\[\mu_{\mcl{X}}\Bigl(\mathcal{C}_\gamma\cap (\theta_0)^{-1}(\mcl{W}^\circ_{B,s})\Bigr)=0\]
in $\groth{k}/(\mbb{L}-1)$ otherwise.
Thus, when $\gamma$ is sufficiently large, the additivity of the motivic measure yields
\[\mu_{\mcl{X}}(\mathcal{C}_{\gamma}) = \sum_{B}\mu_{\mcl{X}}\Bigl((\theta_0)^{-1}(\mcl{W}^\circ_{B,s}\cap \mcl{Y}_k)\Bigr) = \mu_{\mcl{X}}\Bigl((\theta_0)^{-1}(\mcl{Y}_k)\Bigr) =[\mcl{Y}_k]\]
in $\groth{k}/(\mbb{L}-1)$,
where $B$ runs over the nonempty subsets of $A'$.

It remains to prove our claim.
Suppose that $\delta$ is an element in $\mbb{N}^t$ that is different from the zero tuple $\mathbf{0}$.
It suffices to show that  $[\mathcal{D}_\delta]=0$ in $\groth{\mcl{X}_k}/(\mbb{L}-1)$; the result then follows from base change to $\mcl{W}^\circ_{B,s}$.
Since $\ovl{\partial}$ is a relative strict normal crossings divisor, we can cover $\wht{\mcl{X}}$ by formal open subschemes $\mcl{V}$ that admit an étale $R$-morphism
\[h\colon \mcl{V}\to \widehat{\mbb{A}}^d_R = \Spf R\{T_1,\ldots,T_d\}\]
such that the restriction of the divisor $D_1+\ldots+D_t$ to $\mcl{V}$ is defined by the equation $$h^*(T_1\cdots T_r)=0$$
for some $r\geq 0$ less or equal to both $d$ and~$t$.
Using \ref{serre:GrEtale}, we can reduce the computation to the case where $\mcl{X}=\widehat{\mbb{A}}^d_R$ and the component $D_i$ is defined by the equation $T_i=0$, for $i=0,\ldots,r$.
We can also assume without loss of generality that $\delta_i=0$ for every $r<i\leq t$, since otherwise $\mcl{D}_\delta$ would be empty.
For any integrer $n\geq \max_{1\leq i \leq r} \delta_i$, we can write $\mcl{D}_\delta$ as a cylinder over a subvariety $Z$ of~$\funcGr_n(\wht{\mcl{X}})$.
The condition $\mrm{ord}_{D_i}(x)=\delta_i$ translates under the isomorphism
\begin{align*}
\funcGr_n(\wht{\mbb{A}}^d_R)&\isoto \Spec k\Bigl[T_{a,b}\Bigm|a=1,\ldots,d; b=0,\ldots,n\Bigr]\\
&\isoto \funcGr_0(\wht{\mbb{A}}^d_R) \times_k \Spec k\Bigl[T_{a,b}\Bigm|a=1,\ldots,d; b=1,\ldots,n\Bigr]
\end{align*}
as $T_{i,b}=0$ if $b<\delta_i$ and $T_{i,\delta_i}\neq 0$ (see~\ref{GrnAd}).
Since at least one of the coordinates~$\delta_i$ is positive, we can write $Z$ as a product of $\mathbb{G}_{m,\mcl{X}_k}$ with
another variety over $\mcl{X}_k$.
It follows that $[Z]=[\mathcal{D}_\delta]=0$ in $\groth{\mcl{X}_k}/(\mbb{L}-1)$, which is what we needed to prove.
\end{proof}

\begin{num}\emph{Proof of theorem~\ref{thmInvariant}}. \label{proofMainThm}
Let $\approche$ be a tubular datum for~$U$. Then by assumption \eqref{hyp:snc}, we can find a proper birational morphism of $K$-varieties $h\colon X'\to X$ such that $X'$ is regular, $h$ is an isomorphism over~$U$ and the boundary $X'\prive h^{-1}(U)$ is a divisor with strict normal crossings on~$X'$.
Setting $W'_{\alpha}=(h^{\rig})^{-1}(W_{\alpha})$ and $f'_{\alpha}=(h^{\rig})^*f_{\alpha}$ we get an snc tubular datum $\bigl(X',(W'_{\alpha},f'_{\alpha})\bigr)$ which gives rise to the same approximating sequence $U_{\gamma}$. Thus we may assume that $\approche$ is snc.
In this case, it follows from proposition~\ref{prop:compu} that the motivic Serre invariant $S(U_\gamma)$ stabilizes for large $\gamma$.
This proposition also implies that the limit value of $S(U_\gamma)$ only depends
on the pair $(X,U)$, and not on the choice of the pairs $(W_{\alpha},f_{\alpha})$.
On the other hand, the reasoning above shows that the limiting value does not change either if we replace $X$ by another snc compactificatification of~$U$ that dominates~$X$.
Since we can dominate any pair of snc compactifications of $U$ by a common one, we conclude that the limit of the sequence $S(U_{\gamma})$ only depends on $U$, and not on the choice of a tubular datum $\approche$.
\end{num}

\subsection*{Extension to $\mcl{M}_K$}
\begin{num}
In this section, we will assume that properties \eqref{hyp:snc} and \eqref{hyp:adapt} hold, so that the motivic Serre invariant of a smooth $K$-variety $U$ can be defined thanks to \ref{thmInvariant}. Our aim is to define the motivic Serre invariant for arbitrary $K$-varieties $X$; this will be achieved in theorem \ref{extensionMK}.
As a first step, we explain how to compute the motivic Serre invariant of a smooth $K$-variety $U$ from an snc compactification~$X$. Recall that this means that $X$ is a regular compactification of~$U$ and that the boundary $X\prive U$ (with its reduced induced structure) is a strict normal crossings divisor, which we write as $\partial=\sum_{i\in I}E_i.$ \label{nota:E_J}
For every subset $J$ of~$I$, we set \[E_J=\bigcap_{i\in J}E_i.\]
In particular, $E_{\emptyset}=X$.
Note that each of the schematic intersections $E_J$ is again a regular proper $K$-variety, hence admits a weak Néron model by proposition~\ref{prop-regwner}.
\end{num}

\begin{prop}\label{coincidence}
The motivic Serre invariant of $U$ is given by
\[S(U) = \sum_{J\sset I} (-1)^{\norme{J}} S(E_J).\]
\end{prop}
\begin{proof}
By our assumption \eqref{hyp:adapt}, the regular proper $K$-variety $X$ has a weak Néron model $\mcl{X}$ such that the schematic closure $\ovl{\partial}$ of $\partial\cap \mcl{X}_K$ in $\mcl{X}$ is a relative divisor with strict normal crossings over $R$.
In order to apply proposition~\ref{prop:compu}, we need to show that there is a tubular datum for $U$ for which $\mcl{X}$ is adapted.
If the generic fibre of~$\mcl{X}$ is the whole of $X$, it suffices to cover $\mcl{X}$ by opens on which $\ovl{\partial}$ is a principal divisor. The $\pi$-adic completion of this cover yields a cover of~$(\wht{\mcl{X}})_K = X^\rig$.
If it is not, we consider instead a compactification $\ovl{\mcl{X}}$ over~$R$ of the variety obtained by gluing $\mcl{X}$ and $X$ along $U$: we have $\ovl{\mcl{X}}_K = X$ and blowing-up ensures that the schematic closure of~$\partial$ in~$\ovl{\mcl{X}}$ is Cartier.

Proposition~\ref{prop:compu} now implies that
$$S(U)=[(\mcl{X}\prive \ovl{\partial})_k]$$ in $\groth{k}/(\mbb{L}-1)$. For every subset $J$ of $I$, we denote by $\ovl{E}_J$ the schematic closure of $E_J\cap \mcl{X}_K$ in $\mcl{X}$.
We have $\ovl{E}_J(R') = E_J(K') = E_{J,\sm}(K')$ for every finite unramified extension~$R'$ of~$R$ because $\mcl{X}(R') = X(K')$ and $E_J$ is regular.
Furthermore, since $\ovl{E}_J$ is smooth over $R$ at every point of its special fibre, its $R$-smooth locus is a weak Néron model for~$E_J$ and we get in $\groth{k}/(\mbb{L}-1)$
$$\sum_{J\sset I} (-1)^{\norme{J}} S(E_J)=\sum_{J\sset I} (-1)^{\norme{J}} [(\ovl{E}_J)_k]=[(\mcl{X}\prive \ovl{\partial})_k].\qedhere$$
\end{proof}

\begin{num}
We will now extend the definition to \emph{regular} $K$-varieties~$U$ by setting $S(U)\isDef S(U_\sm)$, where $U_\sm$ denotes the $K$-smooth locus of $U$.
Note that this definition is reasonable because the complement of $U_\sm$ in~$U$ does not contain any unramified points (see the proof of proposition~\ref{prop-regwner}).
In particular, any weak Néron model for~$U$ is also a weak Néron model for~$U_\sm$; this \emph{proves} the equality if such a model exists for~$U$.
We can now extend the motivic Serre invariant to arbitrary $K$-varieties thanks to proposition~\ref{coincidence}. More precisely, we have the following result.
\end{num}


\begin{thm}\label{extensionMK}
There exists a unique ring morphism
\[S\colon \mcl{M}_K\to \grothModL\]
such that $S([U]) = S(U)$ for any smooth $K$-variety $U$ and $S([U]) = 0$ if $U$ is a $K$-variety without unramified points. Moreover, if $U$ is a $K$-variety with a weak Néron model $\mcl{U}$, then $S([U])=[\mcl{U}_k]$ in $\grothModL$.
\end{thm}

\begin{proof}
The axioms in the statement imply that $S([U])=S(U_{\mathrm{sm}})$ for every regular $K$-variety $X$ with $K$-smooth locus $U_{\mathrm{sm}}$, since the complement of the smooth locus does not contain any unramified points.
It follows immediately from the presentation given in \ref{presentationGroth} that $S$ is unique if it exists. To prove its existence, we first extend the motivic Serre invariant to a morphism of groups
\begin{equation}\label{eq:Sgroup}
S:K_0(\mathrm{Var}_K)\to \grothModL.
\end{equation}
To this end, it suffices to show that
\begin{equation}\label{eq:addreg}
S(U)=S(V)+S(U\prive V)
\end{equation} whenever $U$ is a connected regular $K$-variety and $V$ is a connected regular closed subvariety of $U$.
We will prove this in three steps.

\smallskip
{\em Step 1: The case where $V$ has no unramified points.}
In this case, $S(V)=0$ and we need to show that $S(U)=S(U\prive V)$. By the definition of the motivic Serre invariant of a regular $K$-variety, we may assume that $U$ is smooth.
Let $X$ be a regular compactification of $U$ such that the closure $\ovl{V}$ of $V$ has transversal intersections with the boundary divisor $\partial=X\prive U$ (we use here our assumption~\eqref{hyp:snc}).
Then $\ovl{V}$ is a regular compactification of~$V$, and it has no unramified points (the same argument as in the proof of \cite[4.6]{112} shows that every weak Néron model of $\ovl{V}$ must have empty special fibre).
We denote by $X'$ the blow-up of $X$ at $\ovl{V}$, by $E\sset X'$ the exceptional divisor and by $E'_i$ the strict transform of each irreductible component $E_i$ of~$\partial$; these are all regular proper $K$-varieties, and \[D = \bigcup_{i\in I} E'_i \cup E\] is a divisor with strict normal crossings on~$X'$.
Applying proposition~\ref{coincidence} to the snc compactification $X'$ of~$U\prive V$, we get, using the notation of~\ref{nota:E_J},
\[S(U\prive V) = \sum_{J\sset I} (-1)^\norme{J} S(E'_J) - \sum_{J\sset I} (-1)^\norme{J} S(E'_J \cap E).\]
By induction on the dimension of~$U$, we can assume that $S(E'_J) = S(E'_J\cap E) + S(E'_J\prive E)$ because $E$ does not contain any unramified points.
Hence we obtain
\[S(U\prive V) = \sum_{J\sset I} (-1)^\norme{J} S(E'_J\prive E) = \sum_{J\sset I} (-1)^\norme{J} S(E_J) = S(U)\]
since $X'\prive E\to X\prive \ovl{V}$ is an isomorphism and $\ovl{V}$ has no unramified points.

Note that this result implies that $S(U)=S(U\prive W)$ for every (not necessarily regular) closed subvariety $W$ of $U$ without unramified points, by induction on a regular stratification of~$W$ (see definition \ref{déf:regStrat}).
This more general equality will be used in the next steps to reduce to the case where $U$ and~$V$ are smooth.

\smallskip
{\em Step 2: The blow-up relation.} Let $Y$ be any regular $K$-variety with a weak Néron model, let $Z$ be a regular closed subvariety of $Y$ and denote by $Y'\to Y$ the blow-up of $Y$ at $Z$, with exceptional divisor $E$. We will show that $$S(Y)-S(Z)=S(Y')-S(E).$$ This relation will be used in Step 3 to finish the proof.
Removing from $Y$ the union of the nonsmooth loci of $Y$ and $Z$, we can assume that both $Y$ and $Z$ are smooth by Step~1. This case was proven in \cite[5.3]{112}.

\smallskip
{\em Step 3: General case.} Now we can prove the general case of \eqref{eq:addreg}.
Removing from $U$ and $V$ the union of the nonsmooth loci of $U$ and $V$, we can assume that both $U$ and~$V$ are smooth by Step 1.
We again take an snc compactification~$X$ of~$U$ such that the closure $\ovl{V}$ of $V$ has transversal intersections with the boundary divisor $\partial=X\prive U$, and we denote by $X'$ the blow-up of $X$ at $\ovl{V}$, with exceptional divisor $E$. Observe that $\ovl{V}$ is an snc compactification of $V$.
We write $\partial=\sum_{i\in I}E_i$ and, for every subset $J$ of $I$, we set $$E_J^\circ=(\bigcap_{i\in J}E_i)\prive (\bigcup_{i\notin J}E_j).$$
We use once more proposition~\ref{coincidence} to compute $S(U)$, $S(U\prive V)$ and $S(V)$.
Assume by induction that \eqref{eq:addreg} holds for every regular variety $U'$ of dimension strict smaller than~$U$ and every closed subvariety $V'$ of~$U'$ (it clearly holds if $U'$ has dimension $0$).
Then we can write
\begin{align*}
S(U) &= S(X) + \sum_{\emptyset\neq J\sset I} (-1)^\norme{J} S(E_J)\\
&= S(X) - \sum_{\emptyset\neq J\sset I} S(E_J^\circ)\\
&= S(X) - \sum_{\emptyset\neq J\sset I} S(E_J^\circ \cap \ovl{V}) - \sum_{\emptyset\neq J\sset I} S(E_J^\circ\prive\ovl{V})\\
&= S(X) +S(V) - S(\ovl{V}) -\sum_{\emptyset\neq J\sset I} S(E_J^\circ\prive \ovl{V}),
\end{align*}
and similarly
\[ S(U\prive V) = S(X')-S(E)-\sum_{\emptyset\neq J\sset I}S(E_J^\circ\prive \ovl{V})\]
because $X'\prive E\to X\prive \ovl{V}$ is an isomorphism. Now the result follows from Step 2.

\smallskip It follows easily from \ref{coincidence} that the map \eqref{eq:Sgroup} is a ring morphism, because weak Néron models are compatible with products (see \cite[4.10]{112}).
It remains to show that, if $U$ is a $K$-variety with a weak Néron model $\mcl{U}$, then $S([U])=[\mcl{U}_k]$ in $\grothModL$.
The existence of a weak Néron model implies that all the unramified points of $U$ are contained in the $K$-smooth locus $U_{\mathrm{sm}}$ of $U$; since $\mcl{U}$ is also a weak Néron model of $U_{\mathrm{sm}}$, we have
$$S([U])=S(U_{\mathrm{sm}})=[\mcl{U}_k]$$ by proposition~\ref{prop:defconsist}.
\end{proof}

\begin{rmq}\label{rmq:cossart}
Cossart and Piltant seem to have recently proved resolution of singularities for three-dimensional quasi-excellent schemes in~\cite[1.1]{182}.
This result together with previous work of Cossart, Jannsen and Saito (\cite[0.3]{183}) yields embedded resolution in three-dimensional schemes.
Our construction consequently produces a motivic Serre invariant for $K$-surfaces, with no assumption on the characteristic of~$K$.
\end{rmq}

\subsection*{Application: Serre's question}

\begin{num}\label{serre:SerreQuestion}
Let $F$ be a field and $G$ a finite $\ell$-group, with $\ell$ a prime number invertible in $F$. In~\cite{151}, Serre asks whether every $G$-action on the affine space $\mbb{A}_F^n$ admits a rational fixed point.
Esnault and Nicaise gave a positive answer, among other cases, when $F$ has characteristic zero and is the fraction field~$K$ of a complete discrete valuation ring with algebraically closed residue field of characteristic different from~$\ell$.
In fact, they even proved that for any $K$-variety~$U$, the nonvanishing of the rational volume of~$U$ (see definition~\ref{ratVol}) modulo~$\ell$ implies the existence of a rational fixed point for the $G$-action on~$U$ (see \cite[5.18]{120}).

In this section, we will assume that conditions \eqref{hyp:snc} and \eqref{hyp:adapt} hold and show how to adapt the arguments of~\cite[5.18]{120} to our extension of the motivic Serre invariant.
This yields in particular a conditional proof of Serre's question with no restriction on the characteristic of~$K$.
\end{num}

\begin{déf}\label{ratVol}
For any $K$-variety~$U$, we define the \emph{rational volume} of $U$ as
\[s(U) \isDef \chi_c(S(U))\in \mbb{Z},\]
where $\chi_c$ denotes the Euler characteristic with proper support of a $k$-variety (see \cite[2.3]{120}).
\end{déf}

\begin{déf}
Let $U$ be a smooth $K$-variety and $G$ a finite group acting on $U$. An \emph{equivariant compactification} of~$U$ is a compactification~$X$ of~$U$ on which the action of~$G$ extends.
A tubular datum $\approche$ is said to be \emph{equivariant} if $X$ is an equivariant compactification of~$U$ and if each rigid variety~$U_\gamma$ of the associated approximating sequence is $G$-invariant.
\end{déf}

\begin{prop}\label{equivApproach}
Let $U$ be a smooth $K$-variety. There exists an equivariant tubular datum for~$U$.
\end{prop}
\begin{proof}
Let $X'$ be an equivariant compactification of~$U$. Blowing-up $X'$ in $X'\prive U$ and normalizing yields a normal equivariant compactification~$X$ of $U$ whose boundary $\partial$ is a Cartier divisor.
By \cite[4.3]{120}, there exists a $G$-model~$\mcl{X}$ of~$X$. Covering $\mcl{X}$ by open subsets $\mcl{V}_\alpha$ on which $\ovl{\partial}$ is defined by some~$f_\alpha$ yields a tubular datum $\bigl(X,((\wht{\mcl{V}}_\alpha)_K,f_\alpha)\bigr)$.
We need to show that every $U_\gamma$ is $G$-invariant.

Let $x\in (\wht{\mcl{V}}_\alpha)_K$ be such that $\norme{f_\alpha(x)}\geq \norme{\pi}^\gamma$ and denote by $\phi_g$ the automorphism of~$X^\rig$ corresponding to the element $g\in G$.
If $\phi_g(x)\in (\wht{\mcl{V}}_\beta)_K$, then we have $\phi_g^*f_\alpha = uf_\beta$ for some unit~$u$ of norm one, because $\ovl{\partial}$ is $G$-invariant.
In particular,
\[
\norme{f_{\beta}(\phi_g(x))} = \norme{\phi_g^*f_\alpha(\phi_g(x))} =\norme{f_\alpha(x)}
\geq \norme{\pi}^\gamma,
\]
from which we conclude that $\phi_g(x)\in U_\gamma$.
\end{proof}

\begin{thm}
Let $\ell$ be a prime number and $G$ a finite $\ell$-group. Let $K$ be a complete discretely valued field whose residue field $k$ is perfect and of characteristic different from~$\ell$.
Let $U$ be a $K$-variety with good $G$-action.
Then \[s(U)\equiv s(U^G) \mod \ell.\]
In particular, if $\ell$ does not divide $s(U)$, then $U^G(K^\sh)\neq\emptyset$.
\end{thm}
\begin{proof}
Since every side of the congruence is additive with respect to $G$-invariant closed immersions, we can assume that $U$ is smooth (one can deduce a presentation for the equivariant Grothendieck ring as in \ref{extensionMK} by stratifying every $K$-variety with good $G$-action into $G$-invariant regular varieties).

Let $\approche$ be an equivariant tubular datum for~$U$. Since the rigidification functor preserves fibred products and immersions, we see from the construction in~\cite[3.1]{133} that 
\[(U^G)^\rig = (U^\rig)^G.\]
Let $(U_\gamma)$ be the associated approximating sequence.
Since each $U_\gamma$ is $G$-invariant, we can consider their subvariety of fixed points~$(U_\gamma)^G$ and we have
\[(U_\gamma)^G = (U^\rig)^G\cap U_\gamma = (U^G)^\rig\cap U_\gamma = (U^G)_\gamma,\]
where $(U^G)_\gamma$ is an element of the approximating sequence associated to the equivariant tubular datum $\bigl(\ovl{U^G},(W_\alpha\cap (U^G)^\rig,f_\alpha)\bigr)$ of the smooth $K$-variety~$U^G$ (see \cite[3.4]{133}). 

Choose $\gamma$ sufficiently large such that $S(U) = S(U_{\gamma})$ and $S(U^G) = S(U_\gamma^G)$ and let $\mfk{U}$ be a weak $G$-Néron model for~$U_{\gamma}$.
By \cite[4.6]{120}, $\mfk{U}^G$ is a weak Néron model for $U_{\gamma}^G$. We then have modulo $\ell$
\[s(U) = \chi_c\bigl(S(U_{\gamma})\bigr) = \chi_c(\mfk{U}_s) = \chi_c(\mfk{U}_s^G) = \chi_c\bigl(S(U_{\gamma}^G)\bigr) = s(U^G).\qedhere\]
\end{proof}

\section{Motivic nearby cycles}
\label{serre:nearbyCycles}

Throughout this section, we will assume that $\car k = 0$.
\subsection*{Motivic nearby cycles with support}

\begin{num}
For any $k$-scheme of finite type $S$, we denote by $\mcl{M}^{\mbb{G}_m}_{S\times \mbb{G}_m}$ the equivariant Grothendieck ring of $S$-varieties with
monomial $\mbb{G}_m$-action defined in \cite[2.3]{53}. It comes equipped with a forgetful morphism to $\mcl{M}_{S}$ which is the composition of the morphism $$\mcl{M}^{\mbb{G}_m}_{S\times \mbb{G}_m}\to \mcl{M}_{S\times \mbb{G}_m}$$ that forgets the $\mbb{G}_m$-action and the morphism $$\mcl{M}_{S\times \mbb{G}_m}\to \mcl{M}_{S}$$ induced by base change to  $S\times \{1\}\isoto S$. Modulo the canonical isomorphism $\mcl{M}^{\mbb{G}_m}_{S\times \mbb{G}_m}\cong \mcl{M}_S^{\hat{\mu}}$ from \cite[2.6]{53} this amounts to forgetting the action of the profinite group scheme of roots of unity $\hat{\mu}$.
\end{num}

\begin{num}
Let $$g\colon Y\to \mbb{A}^1_k=\Spec k[t]$$ be a dominant morphism from a smooth $k$-variety $Y$ of pure dimension~$m$.
We denote by $S_g\in\equiGroth{Y}$ the \emph{motivic nearby cycles} of $g$ as defined in \cite[3.6.2]{53}, where $Y_k$ denotes the zero locus of $g$.
\end{num}

\begin{num} In \cite[3.9]{53} Guibert, Loeser and Merle extend the invariant $S_g$ to a morphism
\[\mcl{S}_g\colon \mcl{M}_Y\to \equiGroth{Y},\]
in such a way that
$\mcl{S}_g\bigl([Y]\bigr) = S_g$.
For this, they define an intermediate object $S_{g,U}$, the \emph{motivic nearby cycles of $g$ with support in $U$}, where $U$ is a dense open subset of~$Y$, as follows.
Denote by $D$ the closed subset $Y\prive U$ of $Y$ with its reduced induced structure.
For every nonnegative integer $n$ we denote by $\mcl{L}_{n}(Y)$ the $n$-th jet scheme of $Y$, parameterizing points on~$Y$ with coordinates in $k[t]/(t^{n+1})$. For every positive integers $d$ and $\gamma$ we consider the reduced subvariety
\[ \mcl{Y}_d^{\gamma d}(g,U) = \Bigl\{\varphi \in \mcl{L}_{\gamma d}(Y) \Bigm| \ord g(\varphi) = d, \; \ord_{D} \varphi\leq \gamma d\Bigr\}\] of $\mcl{L}_{\gamma d}(Y)$,
equipped with the morphism $\mcl{Y}_d^{\gamma d}(g,U)\to\mbb{G}_m$ which sends a jet $\varphi$ to the angular component of $g(\varphi)$. This morphism is piecewise monomial with respect to the action of $\mbb{G}_m$ defined by $a\cdot \varphi(t) = \varphi(at)$, so that we may consider the class of $\mcl{Y}_d^{\gamma d}(g,U)$ in $\equiGroth{Y}$. Introducing
\[Z^\gamma_{g,U} (T) \isDef \sum_{d\geq 1} \bigl[\mcl{Y}_d^{\gamma d}(g,U)\bigr]\, \mbb{L}^{-\gamma md}\, T^d \in\equiGroth{Y}\llbracket T\rrbracket\]
one shows that the limit $-\lim_{T\to\infty}Z^\gamma_{g,U} (T)$ is well-defined in $\equiGroth{Y}$ and independent of $\gamma$ for $\gamma$ sufficiently large. This value is used as the definition of $S_{g,U}$.
See \cite[3.8]{53} for more details. Then it is shown in \cite[3.9]{53} that there exists a unique $\mathcal{M}_k$-linear group morphism  \[\mcl{S}_g\colon \mcl{M}_Y\to \equiGroth{Y}\] such that $\mcl{S}_g([U])$ is the image of $S_{g\circ h,U}$ under the forgetful morphism $\equiGroth{Z}\to \equiGroth{Y}$ for every smooth $k$-variety $Z$, every proper morphism $h\colon Z\to Y$ and every dense open subvariety $U$ of $Z$. Note that uniqueness follows immediately from the fact that the classes of such $Y$-varieties $U$ generate the Grothendieck group $K_0(\mathrm{Var}_Y)$, by resolution of singularities.
\end{num}

The aim of this section is to give another construction of the object $S_{g,U}$ using tubular data to approach $U$ by quasi-compact rigid subvarieties, as in the previous section.

\subsection*{The motivic volume of a quasi-compact rigid variety}
\begin{num}
We set $R = k \llbracket t\rrbracket$ and $K=k\llpar t \rrpar$. For every positive integer $d$ we write $R(d) = R[T]/(T^d-t)$ and $K(d) = \Frac R(d)$.
If $\mcl{X}$ is an stft formal $R$-scheme, we denote by $\mcl{X}(d)$ the stft formal $R(d)$-scheme $\mcl{X}\times_R R(d)$, and similarly for a rigid $K$-variety.
\end{num}




\begin{num}\label{serre:déf:motInt}
We will need a slight refinement of the motivic volume of a rigid $K$-variety constructed in \cite{106} and \cite{49}. First, we recall the definitions given there.
Let $\mcl{X}$ be an stft formal $R$-scheme of pure relative dimension~$m$. Suppose that $\mcl{X}_K$ is smooth over $K$ and that it admits a volume form $\omega$, that is, a nowhere vanishing differential form of maximal degree.
Let $\mcl{Y}\to \mcl{X}$ be a Néron smoothening, i.e.\ $\mcl{Y}$ is a weak Néron model for the rigid variety~$\mcl{X}_K$.
We will again denote by $\omega$ the pullback of $\omega$ to $\mcl{Y}_K$.
For every connected component $C$ of $\mcl{Y}_k$, we denote by $\ord_{C}\omega$ the unique integer $a$ such that $\pi^{-a}\omega$ extends to a relative volume form on an open neighbourhood of the generic point of~$C$ in~$\mcl{Y}$.
Then the change of variables formula for motivic integrals implies that the expression
\[\int_{\mcl{X}}\norme{\omega}\isDef\sum_{C\in \pi_0(\mcl{Y}_k)}[C]\,\mbb{L}^{-\ord_C\omega}\quad \in \mcl{M}_{\mcl{X}_k}\]
only depends on $\mcl{X}$ and $\omega$, and not on the choice of the Néron smoothening (see \cite[2.6.1]{32}), where $\pi_0(\mcl{Y}_k)$ denotes the set of connected components of~$\mcl{Y}_k$.
We call this invariant the \emph{motivic integral} of $\omega$ on $\mcl{X}$.
\end{num}

\begin{num}
Let $\mcl{X}$ be a generically smooth stft formal $R$-scheme and $\omega$ a volume form on~$\mcl{X}_K$.
One can show that the \emph{volume Poincaré series}
\begin{equation}\label{eq:genser}
S(\mcl{X},\omega;T) = \sum_{d\geq 1}\Bigl(\int_{\mcl{X}(d)}\norme{\omega(d)}\Bigr)T^d \in \mcl{M}_{\mcl{X}_k}\llb T\rrb,
\end{equation}
where $\omega(d)$ denotes the inverse image of~$\omega$ in~$X(d)$, is actually an element of the ring
\[\mcl{M}_{\mcl{X}_k}\Bigl[\frac{\mbb{L}^aT^b}{1-\mbb{L}^aT^b}\Bigr]_{(a,b)\in S}\]
for some finite subset $S$ of~$\mbb{Z}\oplus\mbb{Z}_{\geq 1}$.
This ring is equipped with an $\mcl{M}_{\mcl{X}_k}$-linear morphism
\[\lim_{T\to\infty}\colon \mcl{M}_{\mcl{X}_k}\Bigl[\frac{\mbb{L}^aT^b}{1-\mbb{L}^aT^b}\Bigr]_{(a,b)\in S} \to \mcl{M}_{\mcl{X}_k}\]
that sends any product $\prod_{(a,b)\in I} \frac{\mbb{L}^aT^b}{1-\mbb{L}^aT^b}$ onto $(-1)^\norme{I}\, [\mcl{X}_k]$ (see~\cite[8.1]{106}).
The image of $S(\mcl{X},\omega;T)$ under this morphism is called the \emph{motivic volume} of~$\mcl{X}$ and it will be denoted by $\Vol(\mcl{X})$ (the notation $S(\mcl{X};\widehat{K^s})$ was used in \cite{106} and \cite{49}).
It does not depend on~$\omega$ and can be computed explicitly on a log-resolution for the pair $(\mcl{X},\mcl{X}_k)$ (cfr.\ \cite[8.2]{106}).
By an additivity argument we can still define the motivic volume without the assumption that $\mcl{X}_K$ admits a volume form; see the remark following \cite[8.3]{106} and \cite[\S7.4]{49}.
\end{num}

\begin{prop}\label{indepModel}
Let $\mcl{X}$ be an stft formal $R$-scheme and let $V$ be a separated smooth quasi-compact rigid $K$-variety endowed with a morphism $f\colon V\to \mcl{X}_K$. Let $\mcl{V}$ be a formal model of $V$ and let $\mcl{V}\to \mcl{X}$ be a morphism of formal $R$-schemes that extends the morphism $f$.
\begin{enumerate}
\item Let $\omega$ be a volume form on $V$.
Then the image of the integral
\[\int_{\mcl{V}} \lvert\omega\rvert\]
under the forgetful morphism $\mcl{M}_{\mcl{V}_k}\to \mcl{M}_{\mcl{X}_k}$ depends only on $V$, $\mcl{X}$ and~$f$, and not on the chosen model $\mcl{V}$.
\item The image of the motivic volume $\mathrm{Vol}(\mcl{V})$ in $\mcl{M}_{\mcl{X}_k}$ depends only on $V$, $\mcl{X}$ and $f$, and not on~$\mcl{V}$.
\end{enumerate}
\end{prop}

\begin{proof}
Let $\mcl{V}\to\mcl{X}$ and $\mcl{V}' \to\mcl{X}$ be two such models of the morphism $f$. By \cite[2.5.13(2)]{32}, there exists a weak Néron model $\mcl{W}$ of $V$ which dominates $\mcl{V}$ and $\mcl{V}'$.
By definition of the motivic integral of a volume form, we have equalities
$$\int_\mcl{V} \norme{\omega} = \int_\mcl{W} \norme{\omega} =\int_{\mcl{V}' } \norme{\omega}$$ in  $\mcl{M}_{\mcl{X}_k}$.
This proves the first part of the statement; the second part now immediately follows from the first, since the definition of~$\Vol(\mcl{V})$ is based on such integrals.
\end{proof}

\begin{déf}
With the notations of proposition~\ref{indepModel}, we denote by $\mathrm{Vol}_{\mcl{X}}(V)$ the image of $\mathrm{Vol}(\mcl{V})$ in $\mcl{M}_{\mcl{X}_k}$ and we call this invariant the \emph{motivic volume of $V$ over $\mcl{X}$}.
\end{déf}

\subsection*{The motivic volume of an algebraic variety}
\begin{num}
We fix a separated $R$-scheme of finite type $Y$, a smooth connected $K$-variety $U$ of dimension~$m$ and a morphism of $K$-varieties
$h\colon U\to Y_K$, where we set $Y_K=Y\times_R K$.
We denote by $\widehat{Y}$ the $t$-adic completion of $Y$.
The generic fibre~$(\widehat{Y})_K$ is a quasi-compact open rigid subvariety of $(Y_K)^{\rig}$, and the morphism~$h$ induces a morphism of rigid $K$-varieties $h^{\rig}\colon U^{\rig}\to (Y_K)^{\rig}$.

We choose a tubular datum $\approche$ for the smooth $K$-variety $U$ and we consider the associated approximating sequence $(U_\gamma)_{\gamma\geq 0}$ of quasi-compact open rigid subvarieties of $U^{\rig}$.
For every $\gamma\geq 0$, we set \[U'_\gamma=U_\gamma\cap (h^{\rig})^{-1}(\widehat{Y}_K).\]
One can show just as in the case of schemes (see \cite[6.6.4]{EGA-I}) that this is again a quasi-compact open rigid subvariety of $U^{\rig}$.
\end{num}
\begin{thm}\label{thm-motvol}
There exists an integer $\gamma_0\geq 0$ such that, for every $\gamma\geq \gamma_0$, we have
$$
\mathrm{Vol}_{\widehat{Y}}(U'_\gamma) = \mathrm{Vol}_{\widehat{Y}}(U'_{\gamma_0})
$$ in $\mathcal{M}_{Y_k}$.
\end{thm}

\begin{num}\label{serre:NéronSmoothXd}
Before diving into the proof, we recall some facts from \cite[\S 4]{106} about Néron smoothenings of extensions $\mcl{X}(d) = \mcl{X}\times_R R(d)$ of a regular stft formal $R$-scheme~$\mcl{X}$ whose specicial fibre~$\mcl{X}_k = \sum_I N_i\, E_i$ is a strict normal crossings divisor.
For every nonempty subset $J$ of~$I$, we set \label{serre:EJtilde}
\[E_J = \bigcap_{j\in J}E_j \quad \text{and} \quad E_J^\circ=E_J\prive \Bigl(\bigcup_{i\notin J}E_i\Bigr).\]
We cover $E_J^\circ$ by open subsets $\Spf A$ on which $\pi = u\prod_J x_i^{N_i}$, where $u$ is a unit on~$A$ and the $x_i$ are part of a regular system of parameters.
Then the $R$-formal schemes $\Spf A[T]/(uT^{m_J} - 1)$ glue to a Galois cover $\wtl{E^\circ_J}$ of~$E_J^\circ$, where $m_J$ is the greatest common divisor of the multiplicities $N_i$, for $i\in J$.

One can compute explicitly a Néron smoothening for $\mcl{X}(d)$ when $d$ is not $\mcl{X}_k$-linear, that is, if whenever $d$ can be written as an $\mbb{N}$-linear combination
\[d = \sum_I \delta_i N_i,\] then all but one $\delta_i$ is zero. 
The smoothening is given in this case by the $R(d)$-smooth locus of the normalization of~$\mcl{X}(d)$, and its special fibre is canonically isomorphic to
\[\bigsqcup_{N_i|d} \wtl{E_i^\circ}.\]
Furthermore, \cite[5.17]{107} ensures that one can always assume that $d$ is not $\mcl{X}_k$-linear after a finite sequence of blow-ups along some $E_J$. This will allow us to carry over our computations to the general case.
\end{num}

\begin{proof}[Proof of theorem~\ref{thm-motvol}]
Arguing as in the proof of theorem \ref{thmInvariant}, we see that it suffices to prove the result if $\approche$ is an snc tubular datum. In this case, we can make explicit computations on a suitable $R$-model of $U$.

Assume that $\approche$ is an snc tubular datum for $U$ and let $\partial$
be the corresponding Cartier divisor on $X$.
By Nagata's embedding theorem we can factorize $h$ as
\[U\to V\to Y_K,\]
where $U$ is a dense open of~$V$ and $V$ is proper over $Y_K$.
Dominating $X$ by another $snc$ compactification of $V$ (hence of $U$) and pulling back the approach datum to this new compactification, we can assume that the morphism $h:U\to Y_K$ extends to a proper $K$-morphism $\overline{h}:V\to Y_K$ on some open subvariety $V$ of $X$ containing $U$.
Then the morphism $\ovl{h}^\rig$ is proper as well, hence quasi-compact. We denote by $V_0$ the quasi-compact open rigid subvariety $(\ovl{h}^{\rig})^{-1}(\widehat{Y}_K)$ of~$V^\rig$.
Let $\mcl{X}$ be a proper $R$-model of $X$ and denote by $\ovl{\partial}$  the schematic closure of $\partial$ in~$\mcl{X}$.
By \ref{formRig:blowUpModel}, we can modify this model by means of an admissible blow-up in such a way that it satisfies condition 2. of definition~\ref{def:adapted} and that $V_0$ is the generic fibre of an open formal subscheme~$\mcl{V}$ of the $t$-adic completion~$\widehat{\mcl{X}}$ of~$\mcl{X}$.
We can even assume by \cite[4.1(c)]{FRGI} that the morphism $\mcl{V}_K\to \widehat{Y}_K$ induced by $\ovl{h}$ extends to a morphism of formal $R$-schemes $\mcl{V}\to \widehat{Y}$.
By resolution of singularities, we can moreover arrange that $\mcl{X}$ is regular and that $\ovl{\partial}+\mcl{X}_k$ is a divisor with strict normal crossings on $\mcl{X}$.
We write $$\mcl{X}_k=\sum_{i \in I}N_i\, E_i$$
and claim that
\begin{equation}\label{eq:motvol0}
\mathrm{Vol}_{\wht{Y}}(U'_\gamma) =\sum_{\emptyset \neq J\sset I}(1-\mbb{L})^{\norme{J}-1}\bigl[\wtl{E^\circ_J}\cap (\mcl{V}_k\prive \ovl{\partial}_k)\bigr]
\end{equation}
in $\mathcal{M}_{Y_k}$ whenever $\gamma$ is sufficiently large.

By an additivity argument one can assume that $U'_\gamma$ admits a volume form $\omega$, so that we can compute the motivic volume as a limit of a generating series of the form \eqref{eq:genser}.
We are going to compute each coefficient of the series separately.

Let $d\geq 1$ be an integer and assume that it is not $\mcl{X}_k$-linear. We denote by~$\mcl{Z}$ the Néron smoothening of $\wht{\mcl{X}}(d)$ given by the $R(d)$-smooth locus of the normalization of~$\wht{\mcl{X}}(d)$.
The same proof as in~\ref{prop:cylinder} shows that there is a cylinder~$\mcl{C}_\gamma^d$ of~$\funcGr(\mcl{Z})$ whose $k'$-points correspond precisely to the $K'$-points of~$U'_\gamma$.
It follows from (the proof of) \cite[2.4.4]{32} that for any weak Néron model $\mcl{U}_\gamma^d$ of $U'_\gamma(d)$
\[\int_{\mcl{U}_\gamma^d} \norme{\omega(d)} = \int_{\mcl{C}_\gamma^d} \norme{\omega(d)} \:\in \mcl{M}_{\mcl{X}_k}.\]

Keeping the notations of~\ref{prop:compu}, we partition~$\mcl{X}_k$ into locally closed subsets $\mcl{W}_{B,s}^\circ$ and define for every $\delta\in\mbb{N}^t$ the cylinder $\mcl{D}_{\delta}$ as
\[\Bigl\{x\in \funcGr(\mcl{Z})\Bigm|\mrm{ord}_{D_i}(x)=\delta_i\text{ for }i=1,\ldots,t \Bigr\},\]
where we keep writing $D_i$ for the schematic closure of the $i$th component of the divisor $\partial\times_K K(d)$ in~$\mcl{Z}$.
We also set
\[\mcl{D}_{\delta,B} \isDef \mcl{D}_\delta\cap (\theta_0)^{-1}\bigl(\mcl{W}^\circ_{B,s} \cap \mcl{V}_k\bigr).\]
Note that we are a bit sloppy with notation here: by writing $\mcl{W}^\circ_{B,s} \cap \mcl{V}_k$ we really mean its inverse image in the special fibre of~$\mcl{Z}$.
Then
\[\mcl{C}_\gamma^d \cap (\theta_0)^{-1}(\mcl{W}^\circ_{B,s}) = \bigsqcup_\delta \mcl{D}_{\delta,B}\]
where $\delta$ runs over the set of elements in $\mbb{N}^t$ such that
\[d m_{\alpha,j}+N_1\delta_1+\ldots+N_t\delta_t \leq d \gamma\]
for at least one element $(\alpha,j)$ in~$B$.

We are going to compute the motivic measure of each $\mcl{D}_{\delta}$, from which we will deduce the measure of $\mcl{D}_{\delta,B}$ by base change.
Arguing similarly as in~\ref{prop:compu}, we cover $\mcl{Z}$ by formal open subschemes on which we have an étale morphism to the formal affine space~$\wht{\mbb{A}}^d_R$.
For any $\delta\in\mbb{N}^t$, let $K = K(\delta)$ be the subset of $\{1,\dots,t\}$ defined by $\delta_i>0$ if and only if $i\in K$.
The argument given in~\ref{prop:compu} shows that the class of the subvariety~$Z$ of~$\funcGr_n(\mcl{Z})$ over which $\mcl{D}_\delta$ is a cylinder is given by
\[[Z] = [D_K^\circ]\, \mbb{L}^{n(d-\norme{K})}\, \mbb{L}^{\sum_{i\in K} n-\delta_i} (\mbb{L}-1)^{\norme{K}},\]
from which we get
\[\mu_\mcl{X}(\mcl{D}_\delta) = [D_K^\circ]\, (\mbb{L}-1)^{\norme{K}}\, \mbb{L}^{-\sum_{K}\delta_i} \in \mcl{M}_{\mcl{X}_k}.\]
If we set
\[\gamma_B\isDef \max_B (\gamma-m_{\alpha,j}),\]
then we have
\[\mu_\mcl{X}\Bigl(\mcl{C}_\gamma^d\cap (\theta_0)^{-1}(\mcl{W}^\circ_{B,s})\Bigr) = \sum_{\substack{\delta\in\mbb{N}^t\\\sum N_i\delta_i\leq d\gamma_B}} \!\!\!\!\!\!\bigl[D_K^\circ \cap \mcl{W}^\circ_{B,s} \cap \mcl{V}_k\bigr]\, (\mbb{L}-1)^{\norme{K}}\,\mbb{L}^{-\sum_{K}\delta_i}\]
and one gets from \cite[6.13]{106} the following expression for $\int_{\mcl{C}_\gamma^d} \norme{\omega(d)}$ in $\mcl{M}_{\mcl{X}_k}\llb T\rrb$:
\[\sum_B \sum_{K\sset\{1,\dots,t\}} \sum_{\substack{i\in I\\ N_i|d}}\, \bigl[D_K^\circ \cap \wtl{E}_i^\circ \cap \mcl{W}^\circ_{B,s} \cap \mcl{V}_k\bigr]\, (\mbb{L}-1)^{\norme{K}}\!\!\!\!\!\! \sum_{\substack{\delta\in\mbb{N}^K\\\sum N_i\delta_i\leq d\gamma_B}} \!\!\!\!\!\!\mbb{L}^{-(\sum\delta_i) -\tfrac{d}{N_i} \ord_{E_i}\omega} \]
Applying \ref{invBlowupEJ} and resorting to its notation, we can get rid of our assumption that $d$ is not $\mcl{X}_k$-linear thanks to \cite[5.17]{107} and deduce that
\[ \sum_{d\geq 1}\Bigl(\int_{\mcl{U}_\gamma^d}\norme{\omega(d)}\Bigr)T^d = \sum_B \sum_{d\geq 1} i_B^*S(\gamma_B,d)\, T^d\]
in~$\mcl{M}_{\mcl{X}_k}\llb T\rrb$,
where $i_B$ denotes the inclusion $\mcl{W}^\circ_{B,s} \cap \mcl{V}_k\to \mcl{X}_k$.
Finally, we can assume by \ref{limitVolSeries} that all $\gamma_B$ are equal for $\gamma$ sufficiently large, and we eventually get
\[\Vol_{\wht{Y}}(U'_\gamma) = \sum_{\emptyset \neq J\sset I}(1-\mbb{L})^{\norme{J}-1}\bigl[\wtl{E^\circ_J}\cap (\mcl{V}_k\prive \ovl{\partial}_k)\bigr] \in \mcl{M}_{Y_k}.\]
Note that even though we worked throughout over $\mcl{M}_{\mcl{X}_k}$, the right-hand side factors through $\mcl{M}_{\mcl{V}_k}$, hence we can consider its image in~$\mcl{M}_{Y_k}$. 
\end{proof}

\begin{lem}\label{invBlowupEJ}
For any $\wtl{\gamma}\geq 0$, the expression $S(\wtl{\gamma},d)$ defined as
\[\Bigl(\!\!\sum_{ K\sset\{1,\dots,t\}} \!\!\!\![D_K^\circ]\, (\mbb{L}-1)^{\norme{K}} \!\!\!\!\!\!\!\!\sum_{\substack{\delta\in\mbb{N}^K\\\sum N_i\delta_i\leq d\tld{\gamma}}} \!\!\!\!\!\!\mbb{L}^{-\sum\delta_i} \Bigr) \cdot
\Bigl( \sum_{\emptyset\neq J\sset I} [\wtl{E_J^\circ}] \,(\mbb{L}-1)^{\norme{J}-1}\!\!\!\!\!\! \sum_{\substack{\epsilon\in\mbb{N}^J\\ \sum N_i\epsilon_i =d}}\!\!\! \mbb{L}^{-\sum \epsilon_i\mu_i}\Bigl)\]
is invariant in $\mcl{M}_{\mcl{X}_k}\llb T\rrb$ under blow-up of a stratum $E_J$, where $\mu_i = \ord_{E_i}\omega$.
\end{lem}
\begin{proof}
This is immediate from \cite[7.6]{106} and the fact that
\[ [D_K^\circ]\, [\wtl{E}_J^{\prime \circ}] = [D_K^{\prime\circ} \cap \wtl{E}_J^{\prime\circ}]\]
in $\mcl{M}_{\mcl{X}_k}$, where $(\cdot)'$ denotes strict transforms.
\end{proof}

\begin{lem}\label{limitVolSeries}
For $\wtl{\gamma}$ sufficiently large we have
\[-\lim_{T\to \infty} \sum_{d\geq 1} S(\wtl{\gamma},d)\, T^d
= \sum_{\emptyset \neq J\sset I}(1-\mbb{L})^{\norme{J}-1}\bigl[\wtl{E^\circ_J} \prive \ovl{\partial}_k\bigr] \in \mcl{M}_{\mcl{X}_k}\]
\end{lem}
\begin{proof}
This follows from the same argument as in~\cite[3.8]{53}.
\end{proof}

\begin{déf}
We call the limit value
$$\mathrm{Vol}_{\widehat{Y}}(U'_{\gamma_0}) \in \mathcal{M}_{Y_k}$$ in theorem \ref{thm-motvol} the \emph{motivic volume} of~$U$ over $Y$, and we denote it by $\mathrm{Vol}_{Y}(U)$.
\end{déf}

\begin{num}
We can also extend $\Vol_{Y}(\cdot)$ to a ring morphism
\[K_0(\mrm{Var}_{Y_K})\to K_0(\mrm{Var}_{Y_k})\]
in a similar way to~\ref{extensionMK}.
This yields a definition of motivic volume for any $Y_K$-variety, not necessarily smooth over~$K$.
Indeed, the expression~\eqref{eq:motvol0} shows that if $U$ is a $K$-variety embedded in a smooth $K$-variety~$Z$ proper over~$Y_K$ and such that $Z\prive U$ is the support of a divisor with strict normal crossings $\sum_L N_l\,D_l$, then
\[\Vol_{Y}(U) = \sum_{K\sset L} \Vol_{Y}(D_K),\]
where $D_K = \bigcap_{k\in K} D_k$.
Let $V$ be a smooth connected closed subvariety of~$U$.
Then considering an snc compactification of~$Z$ whose boundary intersects $\ovl{V}$ transversally and blowing up along~$\ovl{V}$, one can show as in step~3 of~\ref{extensionMK} that 
\[\Vol_{Y}(U) = \Vol_{Y}(V) + \Vol_{Y}(U\prive V).\]
\end{num}

\subsection*{A non-Archimedean interpretation of the motivic nearby cycles with supports}
\begin{num}
Finally, we consider a morphism of smooth connected $k$-varieties $h\colon U\to Y$ and a dominant morphism
$$g\colon Y\to \mbb{A}^1_k=\Spec k[t].$$ We set $Y_R=Y\times_{k[t]}R$ and $U_K=U\times_{k[t]}K$. The following theorem compares the motivic volume of $U_K$ over $Y_R$ with the motivic nearby cycles of $g$ with supports in $U$.
\end{num}

\begin{thm}
We have
\begin{equation}\label{eq:motvol}
\mathrm{Vol}_{Y_R}(U_K) = \mbb{L}^{-(m-1)}S_{g}([U])
\end{equation} in $\mathcal{M}_{Y_k}$, where $m$ is the dimension of~$U$.
\end{thm}
\begin{proof}
The equality is proven by explicit computation on a suitable compactification of the morphism $U\prive U_k\to Y$.
Note that neither side of the equality \eqref{eq:motvol} changes if we remove $U_k$ from $U$, since this does not affect $U_K$ and $S_g([U_k])=0$. Thus we can assume that $U_k$ is empty.
By Nagata's embedding theorem and resolution of singularities, we can find a  compactification $\ovl{h}\colon Z\to Y$ of the morphism $h\colon U\to Y$ such that $Z$ is proper over $Y$ and smooth over $k$ and $Z\prive U$ is a strict normal crossings divisor on $Z$. Likewise, we can compactify the morphism $g\circ \ovl{h}\colon Z\to \mbb{A}^1_k$ to a proper morphism of $k$-schemes $Z'\to \mbb{A}^1_k$ with $Z'$ smooth over $k$ and $Z'\prive U$ a strict normal crossings divisor on $Z'$.
Then $X=Z'\times_{k[t]}K$ is an snc compactification of $U_K$. Setting $\mcl{X}=Z'\times_{k[t]}R$ and taking for $\mcl{V}$ the $t$-adic completion of $Z$, we see by comparing \eqref{eq:motvol0} with the formula in \cite[3.8.1]{53} that
\[\Vol_{Y_R}(U_K) = \mbb{L}^{-(m-1)} S_{g}([U]).\]
\end{proof}

\begin{num}
It is possible to recover the $\hat{\mu}$-structure on $S_{g}([U])$ (i.e., to recover $S_{g}([U])$ as an element of
$\mcl{M}_{Y_k}^{\hat{\mu}}\cong \equiGroth{Y}$)
in a similar way, using the Galois action on the rigid varieties $\mcl{X}_K\times_K K(d)$ to define a $\hat{\mu}$-action on the coefficients of the generating series \eqref{eq:genser}. This construction is developed in ongoing work of A.~Hartmann.
\end{num}

\cleardoublepage

%

%
  \graphicspath{{\chaptersdir/prelim_log/image/}}%
  \addtocontents{toc}{\protect\newpage}
\chapter{Preliminaries on log geometry}\label{ch:prelim_log}

We give here a crash course on log geometry.
Our favourite introduction are lecture notes by Illusie and Ogus \cite{135}, which can be found with some google-fu.
Two comprehensive references are Ogus's book draft \cite{137.2014} and Gabber and Ramero's notes \cite{148}, both available online.

Log geometry is tightly linked with compactifications and degenerations. It revolves around the notion of log structure (see~\ref{logstr}), which gives a way to encode information about the boundary of a compactification. 
Its power undoubtly comes from the simplicity of the definition, which fits seamlessly in the theory of schemes.
The best illustration of this is given by log smoothness.
One can mimick the theory of infinitesimal liftings via an adapted notion of thickening, which leads to a logarithmic counterpart of the classical notion of smooth and étale morphisms (see \cite[\S 3]{139}).
The analogy allows to translate classical theorems where the sheaf~$\Omega$ of differentials is replaced by the sheaf~$\Omega^\logr$ of log differentials (allowing simple poles). In particular, the sheaf of log differentials of log smooth schemes is locally free, and this allows us to treat our beloved sncd schemes as if they were smooth (see \ref{prelim:sncdLogSmooth}).

The whole theory rests on monoids, this will be the topic of the first four sections. Log schemes are introduced in section~\ref{prelim:log_schemes}.

All monoids considered are assumed to be commutative.

\section{Monoids}

Basically, a monoid is the same thing as a commutative group, except that we do not require every element to have an inverse.
A morphism of monoids is a map that respects internal laws and neutral elements.
Monoids and their morphisms form a category, denoted by $\catMnd$, and this category is complete and cocomplete.
We illustrate this by describing a few limits and colimits.

\begin{num}
The \emph{product} of two monoids $P_1$ and~$P_2$ is given by the set-theoretic product $P_1\times P_2$ endowed with componentwise multiplication.

The \emph{fibred product} of two morphisms $v_1\colon P_1\to Q$ and $v_2\colon P_2\to Q$ is constructed as the submonoid
\[\bigl\{(x,y)\bigm| u_1(x) = u_2(x)\bigr\}\sset P_1\times P_2\]
and denoted by $P_1\times_Q P_2$.
\end{num}

\begin{num}
The \emph{amalgamated sum} of two morphisms $u_1\colon P\to Q_1$ and $u_2\colon P\to Q_2$ is defined as the quotient of the product $Q_1\times Q_2$ by the congruence relation generated by
\[S = \Bigl\{ \Bigl( \bigl(u_1(p),1\bigr), \bigl(1,u_2(p)\bigr)\Bigr) \Bigm| p\in P\Bigr\}.\]
It is denoted by $Q_1\oplus_P Q_2$.
If either $P$, $Q_1$ or $Q_2$ is a group, then this congruence relation can be described as the equivalence relation
\[(q_1,q_2)\sim (q_1',q_2') \iff \text{ there are } a,b\in P \text{ such that } \begin{cases}q_1u_1(b) = q_1'u_1(a)\\ q_2 u_2(a) = q_2' u_2(b).\end{cases}\]
The amalgamated sum $Q_1\oplus_P \{1\}$ is called the \emph{cokernel} of $u_1$ and denoted by $\Cok u_1$ or $Q_1/P$.

If $P = \{1\}$, then $Q_1\oplus_{\{1\}} Q_2 = Q_1\times Q_2$.
We will therefore write $Q_1\oplus Q_2$ for $Q_1\times Q_2$ when we consider it as a colimit with inclusion morphisms instead of a limit with projection morphisms.
\end{num}

\begin{num}
Let $M$ be a monoid.
We denote by $M^\times$ the submonoid of invertible elements of~$M$, which is a group, and by $M^\sharp$ the quotient $M/M^\times = M\oplus_{M^\times}\{1\}$.
A monoid $M$ is called \emph{sharp} if $M^\times = \{1\}$ and, unsurprisingly, $M^\sharp$ is always sharp.
\end{num}

\begin{déf}
A morphism of monoids $\phi\colon P\to Q$ is called \emph{local} if $\phi^{-1}(Q^\times) = P^\times$ (the inclusion $\sups$ always holds).
This is equivalent to $\phi(P\prive P^\times)\sset Q\prive Q^\times$ (see also \ref{maxIdeal}).
We denote by $\Hom_{\loc}(P,Q)$ the set of local morphisms between the monoids $P$ and~$Q$.
It canonically carries a structure of monoid.

For a monoid $M$, we denote by $M^\vee$ its \emph{dual monoid}, defined by the monoid of morphims $\Hom(M,\mbb{N})$.
We will also denote by $M^{\vee,\loc}$ its submonoid $\Hom_\loc(M,\mbb{N})$.
\end{déf}

We introduce now some fundamental properties of monoids.

\begin{déf}
A monoid $M$ is called \emph{finitely generated} if there are finitely many elements $x_1,\dots,x_n\in M$ such that every $x\in M$ can be written as a product $x = \prod_{i=1}^r x_i^{n_i}$ with $n_i\geq 0$.
\end{déf}

\begin{num}\label{localization}
Let $M$ be a monoid and $S\sset M$ a submonoid of~$M$.
We define the \emph{localization} $S^{-1}M$ as the quotient of $M\times S$ by the equivalence relation
\[(x,s)\sim (y,t) \iff axt = ays \text{ for some } a\in S.\]
We have a canonical morphism $M\to S^{-1}M, x\mapsto (x,1)$.
When $S = M$, we obtain the \emph{groupification} of~$M$, denoted by $M^\gp$.
It is an abelian group.
\end{num}

\begin{rmq}\label{prelim:gpPreservesColim}
The groupification yields a functor $(\mhyp)^\gp\colon \catMnd\to \catAb$ to the category of abelian groups, which is a left adjoint for the forgetful functor $\catAb\to \catMnd$.
This means that for every abelian group $G$, we have canonical bijections
\[\Hom_\catAb(M^\gp, G)\isoto \Hom_\catMnd(M,G).\]
In particular, it follows from abstract nonsense that
\[(Q_1\oplus_P Q_2)^\gp = Q_1^\gp\oplus_{P^\gp} Q_2^\gp,\]
for any pair of morphisms of monoids $u_1\colon P\to Q_1$ and $u_2\colon P\to Q_2$.
\end{rmq}

Just as localization morphisms of rings need not be injective, $M$ doesn't always embed into~$M^\gp$.
For example, the groupification of the monoid $\mbb{N}$ endowed with the law $x\cdot y\isDef \max(x,y)$ is trivial.
Pointed monoids also give examples of such monoids (see \ref{pointedMonoid}).
This motivates the following definition.

\begin{déf}
A monoid $M$ is \emph{integral} if the canonical morphism $M\to M^\gp$ is injective.
It is called \emph{fine} if it is moreover finitely generated.
\end{déf}

Finally, the property of saturation below is reminiscent of the notion of normal ring.

\begin{déf}
An integral monoid $M$ is called \emph{saturated} if for every $x\in M^\gp$, $x^\alpha\in M$ for $\alpha\geq 1$ implies that $x\in M$.
\end{déf}

\begin{rmq}
Note that if $M$ is sharp and saturated, then $M^\gp$ is a torsion-free abelian group.
\end{rmq}

\begin{déf}
If $M$ is any monoid, its \emph{integral closure} $M^\intg$ is the image of~$M$ in~$M^\gp$.
Moreover, if $M$ is an integral monoid, we denote by $M^\sat$ the \emph{saturation} of~$M$, which is the submonoid of~$M^\gp$ consisting of all elements $x$ such that $x^\alpha\in M$ for some $\alpha\geq 1$.
\end{déf}

\begin{num}\label{intgMorph}
Even when $Q_1, Q_2$ and~$P$ are integral monoids, the amalgamated sum $Q_1\oplus_P Q_2$ need not be integral.
This motivates the following definition.
A morphism of integral monoids $P\to Q_1$ is called \emph{integral} if $Q_1\oplus_P Q_2$ is integral for every integral monoid $Q_2$.
We will see in \ref{log:defEpi} that any morphism $\mbb{N}\to P$ is integral.

The same kind of problem arises for saturated monoids. Replacing every occurence of ``integral'' by ``saturated'' yields the corresponding notion of \emph{saturated morphism}.
\end{num}

The monoids we will be mostly interested in are fine and saturated monoids, abbreviated by \emph{fs monoids}.

\section{Ideals and pointed monoids}

\begin{déf}\label{déf:ideal}
An \emph{ideal} of a monoid~$M$ is a subset~$I$ of~$M$ such that $ax\in I$ for every $x\in I$ and $a\in M$.
An ideal $I$ is \emph{prime} if $I\neq M$ and $xy\in I$ implies that either $x$ or $y$ belong to $I$, or equivalently, $I$ is prime if $M\prive I$ is a submonoid of~$M$.
\end{déf}

\begin{déf}
Submonoids of~$M$ arising as complements of prime ideals are called \emph{faces}.
More explicitly, a submonoid $F$ of $M$ is a face if $xy\in F$ with $x,y\in M$ implies that both $x$ and $y$ belong to $F$.
The terminology comes from the analogy with convex polyhedral cones, see \ref{facesCone}.
\end{déf}

\begin{num}\label{maxIdeal}
We easily see that unions of (resp.\ prime) ideals are (resp.\ prime) ideals.
In particular every monoid $M$ admits a unique maximal ideal $M^+\isDef M\prive M^\times$ and a unique minimal ideal $\emptyset$.
\end{num}

\begin{déf}
Let $\mfk{p}$ be a prime ideal of a monoid~$M$.
The \emph{height} of~$\mfk{p}$ is the maximum length~$n$ of a chain of prime ideals
\[\mfk{p} = \mfk{p}_n\supsn \cdots\supsn \mfk{p}_0 = \emptyset.\]
It is denoted by $\hgt \mfk{p}$.
We also define the \emph{dimension} of~$M$ as the height of its maximal ideal~$M^+$.
\end{déf}

\begin{num}\label{locM_p}
Since the complement of a prime ideal $\mfk{p}$ is a submonoid of~$M$, we can consider the localization
\[M_\mfk{p}\isDef (M\prive \mfk{p})^{-1}M\]
(see \ref{localization}).
When $M$ is fine, we have
\[\hgt \mfk{p} = \dim M_\mfk{p} = \rk (M_\mfk{p}^\sharp)^\gp.\]
(see \cite[I.1.4.7]{137.2014} for the last equality). 
Note that $M^\gp = M_\emptyset$ and that $M_\mfk{p}^\sharp\isDef (M_\mfk{p})^\sharp$ equals the quotient of monoids $M/(M\prive\mfk{p})$.

This construction provides the basis for the definition of the spectrum of a monoid, see section~\ref{prelim:monoidalSpaces}.
\end{num}

\subsection*{Pointed monoids}
In order to construct quotients of monoids by ideals (see \ref{quotientM/I}), we will need to formally add a zero element.
This leads to the notion of \emph{pointed monoid}.

\begin{déf}
Let $M$ be a monoid. An \emph{$M$-module} is a set $S$ endowed with a map $M\times S\to S, (x,s)\mapsto x\cdot s$ such that $x\cdot (y\cdot s) = (xy)\cdot s$ and $1\cdot s = s$ for every $x,y\in M$ and $s\in S$.
\end{déf}

\begin{déf}\label{pointedMonoid}
A \emph{pointed monoid} is a monoid $M$ endowed with a morphism of $M$-modules $0_M\colon 0\to M$, where $0$ denotes the trivial $M$-module.
We will denote by $\catMnd^\bullet$ the category of pointed monoids.
The forgetful functor $\catMnd^\bullet\to \catMnd$ admits a left adjoint
\[\catMnd\to \catMnd^\bullet,\quad M\mapsto M^\bullet\]
where $M^\bullet$ is the module $M\oplus 0$ endowed with the inclusion map $0_M\colon 0\to M\oplus 0$.
Observe that $M^\bullet$ can be seen as a monoid in a canonical way.
We define similarly the notion of \emph{pointed module} of a pointed monoid and for any $M^\bullet$-module~$S$, we define $S^\bullet$ as the pointed module $S\oplus 0$.
\end{déf}

\begin{eg}
Let $A$ be a ring and $I$ an ideal. The multiplicative monoid $(A,\cdot)$ is naturally a pointed monoid with $0$ as pointed element, and $I$ is a pointed $A$-module.
\end{eg}

\begin{num}\label{quotientM/I}
Let $M$ be a monoid and $I\sset M$ an ideal. We define the quotient $M/I$ as the cokernel (as pointed $M$-modules) of the induced morphism $I^\bullet\to M^\bullet$. It has the structure of a pointed monoid and we have a canonical morphism of monoids $M\to M/I$.
\end{num}

\begin{num}
Let $A$ be a ring. The forgetful functor $A\mhyp \catAlg\to \catMnd^\bullet$ admits a left adjoint
\[\catMnd^\bullet\to A\mhyp \catAlg, \quad (M,0_M)\mapsto A\langle M\rangle\isDef A[M]/(e_0),\]
where $e_0$ denotes the image of $0_M$ in $A[M]$.
Observe that $A\langle M^\bullet\rangle = A[M]$.
\end{num}

\begin{prop}\label{isoMonoidAlgebras}
Let $\mfk{p}$ be a prime ideal of a monoid $M$ and $A$ a ring. Then the canonical morphism of $A[M]$-algebras
\[A[M\prive\mfk{p}]\to A\langle M/\mfk{p}\rangle\]
is an isomorphism.
\end{prop}
\begin{proof}
We see $A[M\prive\mfk{p}]$ as an $A[M]$-algebra by
\[A[M]\to A[M\prive\mfk{p}],\quad e_m\mapsto \begin{cases} e_m & m\notin \mfk{p}\\ 0 & m\in \mfk{p}\end{cases}, \]
which is well-defined since $M\prive\mfk{p}$ is a submonoid of $M$.
But this is exactly the definition of $A\langle M/\mfk{p}\rangle$.
\end{proof}

\section{Valuations and divisors}
\label{prelim_log:sec:valEtDiv}

\begin{déf}\label{déf:valuation}
Let $M$ be an fs monoid and $\mfk{p}$ a prime ideal of height~1.
Then $M_\mfk{p}^\sharp$ is a sharp fs monoid of dimension 1, hence canonically isomorphic to~$\mbb{N}$ (see \cite[3.4.16]{148}).
The induced map $v_\mfk{p}\colon M\to\mbb{N}$ is called the \emph{valuation} at~$\mfk{p}$.
It extends to a map $M^\gp\to \mbb{Z}$ which we keep writing $v_\mfk{p}$.
\end{déf}

\begin{déf}
Let $M$ be an integral monoid. A \emph{fractional ideal} of $M$ is a nonempty $M$-submodule $I$ of $M^\gp$ such that $xI\sset M$ for some $x\in M$.
\end{déf}

\begin{déf}
Let $I$ be a fractional ideal of an integral monoid $M$. We set
\[I^{-1}\isDef \bigl\{x\in M^\gp \bigm| xI\sset M\bigr\} \text{ and } I^* = (I^{-1})^{-1}.\]
The fractional ideal $I$ is said to be \emph{reflexive} if $I^* = I$.
Note that $I^*$ (and even~$I^{-1}$) is always a reflexive fractional ideal of~$M$.
We denote by $\Divr(M)$ the set of reflexive fractional ideals of $M$ and we call its elements \emph{divisors} of $M$.
\end{déf}

\begin{num}
We endow $\Divr(M)$ with a structure of commutative monoid by setting
\[I\odot J \isDef (IJ)^*.\]

For every morphism $\phi\colon P\to Q$ of integral monoids, $IQ \isDef \phi^\gp(I)Q$ is a fractional ideal of $Q$ and $\phi$ induces a morphism of monoids
\[\Divr(\phi)\colon \Divr(P)\to \Divr(Q), \quad I\mapsto (IQ)^*.\]
\end{num}

\begin{prop}[{\cite[3.4.25 and 3.4.32]{148}}]\label{IsoDivM}
Let $M$ be an fs monoid and denote by $D$ the set of its height-one prime ideals. Then $\Divr(M)$ is a group and we have an isomorphism
\[\mbb{Z}^{\oplus D}\toisoto \Divr(M), \quad (n_{\mfk{p}})_D\mapsto \bigcap_D \mfk{m}_{M_\mfk{p}}^{n_{\mfk{p}}},\]
where the intersection takes place in~$M^\gp$ via the canonical inclusions $M_\mfk{p}\to M^\gp$.
Furthermore, under this isomorphism, the fractional ideal generated by an element $a\in M^\gp$ coincides with the tuple $(v_\mfk{p}(a))_D$ (notation of~\ref{déf:valuation}).
\end{prop}

\begin{num}
Let $I$ be a fractional ideal of an fs monoid $M$. If $I =\bigcap_D \mfk{m}_{M_\mfk{p}}^{n_{\mfk{p}}}$ we will denote by $v_{\mfk{p}}(I)$ the integer $n_\mfk{p}$. 
\end{num}

\section{Cones}

The study of sharp fs monoids is essentially equivalent to the study of strictly convex rational polyhedral cones.
This point of view will prove especially useful for computations (see for instance \ref{log:poles} and section \ref{log:sectionNewton}).

\begin{déf}
Let $L$ be a free abelian group of finite rank and $\sigma\sset L_\mbb{R} \isDef L\otimes_\mbb{Z}\mbb{R}$ a convex cone.
We say that $\sigma$ is
\begin{itemize}
\item \emph{polyhedral} if it admits a finite number of generators,
\item \emph{rational} if it is generated by elements of~$L$,
\item \emph{strictly convex} if it does not contain any nonzero linear subspace.
\end{itemize}
\end{déf}

\begin{num}
Let $\sigma$ be a rational polyhedral cone in $L_\mbb{R}$.
The \emph{dual cone} is defined by 
\[\sigma^\vee\isDef \bigl\{u\in L_\mbb{R}^\vee \bigm| u(x)\geq 0 \text{ for all } x\in \sigma\bigr\},\]
where $L_\mbb{R}^\vee$ denotes the dual vector space of~$L_\mbb{R}$.
It is a rational polyhedal cone as well (see \cite[1.3]{166}).

\label{prelim:déf:faceCone}
A \emph{face} of $\sigma$ is a subcone of the form
\[\tau = \bigl\{x\in\sigma \bigm| u(x) = 0 \text{ for some } u\in L_\mbb{R}^\vee\bigr\}.\]
Faces of rational polyhedral cones are rational polyhedral cones (\emph{loc.\ cit.}).

\label{déf:dimCone}
The \emph{dimension} of the cone $\sigma$, denoted by $\dim \sigma$, is the dimension of the $\mbb{R}$-vector space it generates.
We have $\dim \sigma^\vee = \rank L - \dim \sigma^\times$, where $\sigma^\times$ denotes the largest linear subspace contained in $\sigma$ (see \cite[3.3.12(ii)]{148}).
\end{num}


Let's see how this relates to monoids.

\begin{num}\label{prelim:mndCones}
Let $M$ be a fine monoid and denote by
\[\sigma(M) \isDef \bigl\{{\textstyle \sum}_i (\lambda_i\otimes  x_i) \bigm| \lambda_i\in \mbb{R}_{\geq 0} \text{ and } x_i\in M\bigr\}\]
the cone generated by $M$ in $(M^\gp)_\mbb{R}$.
This is a rational polyhedral cone.

By \cite[I.2.3.6]{137.2014}, the canonical morphism of monoids $M\to \sigma(M)$ induces a bijection between the faces of~$M$ and the faces of the cone $\sigma(M)$.\label{facesCone}
In particular, $\sigma(M)$ is strictly convex if and only if $M$ is sharp.
We also have (see formulas in \ref{déf:dimCone} and \ref{locM_p})
\[\dim \sigma(M)^\vee = \dim M^\gp - \dim M^\times = \dim M.\]
Finally, we have $\sigma(M)\cap M^\gp = M^\sat$ (cfr.\ \cite[I.2.3.5(4)]{137.2014}).
\end{num}

\begin{num}\label{prelim:valuations}
Let $F$ be a facet of~$M$ (i.e.\ a face of codimension one). Then $\mfk{p} = M\prive F$ is a height-one prime ideal and we have a valuation
\[v_\mfk{p}\colon M\to M_\mfk{p}^\sharp\isoto \mbb{N}.\]
But since $M_\mfk{p}^\sharp\isoto M/F$, we see that the kernel of $v_\mfk{p}$ is precisely~$F$.
Hence $v_\mfk{p}$ can be identified with the inward normal direction of the facet $\sigma(F)$ of the cone $\sigma(M)$.

In particular, if $\mfk{p}_1,\dots,\mfk{p}_n$ are the height-one prime ideals of $M$, then the dual cone $\sigma(M)^\vee$ is generated over $\mbb{R}_{\geq0}$ by the valuations $v_{\mfk{p}_i}$ (see \cite[3.3.11]{148}).
Moreover,
\[\sigma(M^{\vee,\loc}) = \bigl\{{\textstyle \sum}_i \lambda_i v_{\mfk{p}_i}\bigm| \lambda_i>0 \text{ for every }1\leq i \leq n\bigr\}.\]
We say that $\sigma(M^{\vee,\loc})$ is \emph{strictly positively spanned} by the $v_{\mfk{p}_i}$.
\end{num}

\section{Monoidal spaces and fans}
\label{prelim:monoidalSpaces}
\begin{num}\label{déf:MndSp}
A \emph{monoidal space} is a topological space $T$ endowed with a sheaf $\mcl{M}_T$ of monoids.
A \emph{morphism of monoidal spaces} $f\colon T'\to T$ consists of a continuous map $f$ and a map of sheaves $h\colon f^{-1}\mcl{M}_T\to \mcl{M}_{T'}$ such that for every $t\in T'$ the map of monoids $h_t\colon \mcl{M}_{T,f(t)}\to \mcl{M}_{T',t}$ is local.
We will denote by $\catMndSp$ the corresponding category.

A monoidal space $(T,\mcl{M}_T)$ is called \emph{sharp} if $\mcl{M}_T$ is sharp, i.e.\ $\mcl{M}_T(U)$ is a sharp monoid for every open $U$ of~$T$.
We can associate to every monoidal space $(T,\mcl{M}_T)$ a sharp monoidal space $T^\sharp \isDef (T,\mcl{M}_T^\sharp)$ where $\mcl{M}_T^\sharp$ is the sheaf associated to the presheaf
\[U\mapsto \mcl{M}_T(U)^\sharp.\]
\end{num}

A fundamental example of monoidal space is given by the spectrum of a monoid.
We describe this construction below.

\begin{num}
For $f\in P$ we set $P_f\isDef S^{-1}P$, where $S = \{f^n\mid n\geq 0\}$.
The sets
\[D(f)\isDef \bigl\{\mfk{p}\in \Spec P\mid f\notin \mfk{p}\bigr\}\]
generate a topology on the set $\Spec P$ of prime ideals of~$P$.
If we define $M_P(D(f))\isDef P_f$, then $(\Spec P, M_P)$ is a monoidal space, called the \emph{spectrum} of~$P$ and simply denoted by~$\Spec P$.
It represents the functor
\[\catMndSp\to \catSet,\quad T\mapsto \Hom_\catMnd\bigl(P, \mcl{M}_T(T)\bigr),\] which means that
\[\Hom_\catMndSp(T,\Spec P) \isoto \Hom_\catMnd\bigl(P, \mcl{M}_T(T)\bigr).\]
We easily see that the stalk $M_{P,\mfk{p}}$ of the sheaf $M_P$ at~$\mfk{p}$ is nothing else than the localization $P_\mfk{p}$.

The monoidal space $(\Spec P)^\sharp$ is called the \emph{sharp spectrum} of~$P$.
\end{num}

To conclude this section, we introduce the notion of fan, which is to spectra of monoids what schemes are to spectra of rings.
Fans are a fundamental tool in the study of log regular schemes (cfr.\ \ref{FanInducedByLogReg}). We will use them extensively in chapter~\ref{ch:log}.

\begin{déf}
A \emph{fan} is a sharp monoidal space $(F,M_F)$ that can be covered by open subsets isomorphic to sharp spectra of monoids.
We denote by $\catFan$ the full subcategory of sharp monoidal spaces whose objects are fans.

A fan is called
\begin{itemize}
\item \emph{locally finite} (resp.\ \emph{finite}) if it can be covered by (resp.\ finitely many) spectra of finitely generated monoids.
\item \emph{saturated} if it can be covered by spectra of saturated monoids.
\item \emph{locally fs} (resp.\ \emph{fs}) if it can be covered by (resp.\ finitely many) spectra of fs monoids.
\end{itemize} 
\end{déf}

This notion of fan is essentially equivalent to the classical notion of a fan of polyhedral cones as defined in \cite{166}.
If $\Delta$ is a classical fan in an abelian group~$L$, then a cone $\tau$ of $\Delta$ corresponds to a point $t$ of an fs fan~$F$ and $M_{F,t}^\vee$ should be seen as $\tau\cap L$.
Furthermore, an inclusion of faces $\sigma\sset \tau$ corresponds to a specialization morphism $M_{F,t}\to M_{F,s}$ where $t$ and $s$ are the points respectively associated to $\tau$ and~$\sigma$.
We invite the interested reader to consult \cite[9.5]{143} for more details on this correspondence.

Just as in the classical case, we have a notion of subdivision of a fan, which takes a particularly elegant form.

\begin{num}
Let $F$ be a fan and $P$ a monoid. The \emph{set of $P$-points} of~$F$ is
\[F(P)\isDef \Hom_{\catFan}\bigl((\Spec P)^\sharp, F\bigr).\]
Let $\vphi\in F(P)$ and let $U$ be an open subset of $F$ which contains $\vphi(P^+)$.
Since $P^+$ is a maximal ideal, the inverse image $\vphi^{-1}(U)$ is the whole of $(\Spec P)^\sharp$.
In particular, the induced morphism of sheaves factors through the stalk $M_{F,\vphi(M^+)}\to P^\sharp$ and this morphism is local (by definition~\ref{déf:MndSp}).
Hence we have shown that
\[ F(P) = \bigsqcup_{t\in F} \Hom_{\catMnd,\loc} (M_{F,t}, P^\sharp).\]
\end{num}

\begin{déf}\label{déf:subdivision}
Let $F$ be a fan. A \emph{subdivision} of $F$ is a morphism of fans $\vphi\colon F'\to F$ such that
\begin{enumerate}
\item ${M}_{F,\varphi(t)}^\gp\to {M}_{F',t}^\gp$ is surjective for every $t\in F'$,
\item the canonical map $F'(\mbb{N})\to F(\mbb{N})$ is bijective.
\end{enumerate}
\end{déf}

\begin{num}
Subdivisions of fans allow to desingularize them (in the sense given by proposition~\ref{prelim:desingFan} below).
This result lies at the heart of desingularization of log smooth schemes over a discrete valuation ring to sncd schemes (see \ref{prelim:inverseSubdiv} and \ref{log:equivReg}).
For a fan~$F$, we will denote by $F_\reg$ the set $\{t\in F\mid M_{F,t}\isoto\mbb{N}^{r(t)}\}$, where $r(t)$ denotes the dimension of~$M_{F,t}$.
\end{num}

\begin{prop}[{\cite[3.6.31]{148}}] \label{prelim:desingFan}
Let $F$ be a locally fs fan. There is a subdivision $\vphi\colon F'\to F$ such that $F'_\reg = F'$ and that restricts to an isomorphism $\vphi^{-1}(F_\reg)\toisoto F_\reg$.
\end{prop}

\section{Log structures and charts}
\label{prelim:log_schemes}
\begin{num}\label{logstr}
Let $X$ be a scheme.
A \emph{pre--log structure} on $X$ consists of a sheaf of monoids~$\mcl{M}_X$ on~$X$ together with a morphism of sheaves of monoids $\alpha\colon \mcl{M}_X\to (\mcl{O}_X, \cdot)$ to the multiplicative monoid of $\mcl{O}_X$.
A pre--log structure is called a \emph{log structure} if moreover $\alpha^{-1}(\mcl{O}_X^\times)\to \mcl{O}_X^\times$ is an isomorphism.
This is equivalent to saying that $\alpha_x\colon \mcl{M}_{X,x}\to \mcl{O}_{X,x}$ is local and induces an isomorphism $\mcl{M}_{X,x}^\times\toisoto \mcl{O}_{X,x}^\times$ for every point $x\in X$.

Every pre--log structure $\alpha\colon \mcl{M}_X\to \mcl{O}_X$ induces a log structure $\mcl{M}_X^a$, called the \emph{associated log structure}.
It is given on the stalks by
\begin{equation}\label{eq:Ma_stalks}
\mcl{M}_{X,x}^a = \mcl{M}_{X,x}\oplus_{\alpha_x^{-1}(\mcl{O}_{X,x}^\times)}\mcl{O}_{X,x}^\times \to \mcl{O}_{X,x}.
\end{equation}
\end{num}

Probably the best way to motivate this definition is to see how it fits naturally in the definition of log differentials (see \ref{defOmegaLog}).

\begin{rmq}
Some authors endow log schemes with the étale topology.
We will stick throughout with the Zariski topology, which is better suited to deal with fans of log regular schemes (see~\ref{FanInducedByLogReg}).
\end{rmq}

\begin{eg}
Every scheme admits a trivial log structure defined by $\mcl{M}_X = \mcl{O}_X^\times$.
It is the structure associated to the trivial sheaf of monoids.
We denote by~$X^\circ$ the log scheme obtained by endowing $X$ with its trivial log structure.
We have a canonical morphism of log schemes $X\to X^\circ$.
\end{eg}

\begin{eg}
Let $D$ be a divisor on a scheme~$X$ and denote by~$U$ its complement. The sheaf of monoids
\[\mcl{M}_X(V)\isDef \bigl\{f\in\mcl{O}_X(V) \bigm| \rest{f}{U\cap V} \text{ is invertible}\bigr\}\]
is a log structure, called the \emph{divisorial log structure} associated to~$D$.
\end{eg}

\begin{déf}
A \emph{morphism of log schemes} $(f,\phi)\colon X\to Y$ is a morphism~$f$ between the underlying schemes together with a morphism of sheaves $\phi\colon f^{-1}\mcl{M}_Y\to \mcl{M}_X$ such that $\alpha_X\circ\phi = f^\sharp\circ f^{-1}(\alpha_Y)$:
\begin{center}
\begin{tikzpicture}
    \matrix (m) [matrix of math nodes, row sep=3em,
    column sep=3em, text height=1.5ex, text depth=0.25ex]
    { f^{-1}\mcl{M}_Y & \mcl{M}_X \\
f^{-1}\mcl{O}_Y & \mcl{O}_X\\};
    \path[arrow]
    (m-1-1) edge node[auto] {$\phi $} (m-1-2)
            edge node[auto,swap] {$f^{-1}\alpha_Y$}  (m-2-1)
    (m-1-2) edge node[auto] {$\alpha_X$}(m-2-2)
(m-2-1) edge node[auto] {$f^\sharp$} (m-2-2);

\end{tikzpicture}\end{center}
In particular it is also a morphism of monoidal spaces.
\end{déf}

\begin{num}
If $(f,\phi)\colon X\to Y$ is a morphism of log schemes, then the morphism
\[\alpha_X\circ \phi\colon f^{-1}\mcl{M}_Y\to\mcl{M}_X\to \mcl{O}_X\]
gives a pre--log structure on $X$ and $\phi$ extends to a morphism $\phi^a\colon (f^{-1}\mcl{M}_Y)^a\to \mcl{M}_X$.
We say that $f$ is \emph{strict} if $\phi^a$ is an isomorphism.
\end{num}

\begin{num}
Let $X$ be a log scheme.
A \emph{chart} for $X$ is a strict morphism $c\colon X\to \Spec \mbb{Z}[P]$ for some monoid $P$, where we endow $\Spec \mbb{Z}[P]$ with the standard log structure induced by the canonical morphism $P\to \mbb{Z}[P]$.
It is called \emph{coherent} (resp.\ fine, resp.\ fs) if $P$ is finitely generated (resp.\ fine, resp.\ fs).
Similarly, a log scheme $X$ is called \emph{coherent} (resp.\ fine, resp.\ fs) if it locally admits a coherent (resp.\ fine, resp.\ fs) chart.
To give a chart is equivalent to giving a morphism $\beta\colon P\to \mcl{M}_X(X)$ of monoids inducing an isomorphism $P_X^a\toisoto \mcl{M}_X$, where $P_X$ denotes the constant sheaf on~$X$.

We will say that a chart $P\to \mcl{M}_X(X)$ is \emph{local at a point}~$x\in X$ if the induced morphism $P\to\mcl{M}_{X,x}$ is local.
\end{num}

\begin{num}\label{chartForMorphism}
Likewise a \emph{chart for a morphism} $f\colon X\to Y$ of log schemes consists of a morphism $u\colon P\to Q$ of monoids and charts $c_Q\colon X\to \Spec \mbb{Z}[Q]$ and $c_P\colon Y\to \Spec \mbb{Z}[P]$ fitting into a commutative diagram:
\begin{equation}\label{dia:chartForMorphism}
\begin{tikzpicture}[baseline=(current  bounding  box.center)]
    \matrix (m) [matrix of math nodes, row sep=3em,
    column sep=3em, text height=1.5ex, text depth=0.25ex]
    { X & \Spec \mbb{Z}[Q] \\
Y & \Spec \mbb{Z}[P]\\};
    \path[arrow]
    (m-1-1) edge node[auto] {$c_Q $} (m-1-2)
            edge node[auto,swap] {$f$}  (m-2-1)
    (m-1-2) edge node[auto] {$\Spec \mbb{Z}[u]$}(m-2-2)
(m-2-1) edge node[auto] {$c_P$} (m-2-2);

\end{tikzpicture}
\end{equation}
Such a chart is said to be \emph{fs} if both $P$ and~$Q$ are fs.
\end{num}

\begin{num}
If $\beta\colon P\to \mcl{M}_X(X)$ is a chart, then we have at each point $x\in X$ a diagram of amalgamated sums (see \eqref{eq:Ma_stalks})
\begin{center}
\begin{tikzpicture}
    \matrix (m) [matrix of math nodes, row sep=3em,
    column sep=3em, text height=1.5ex, text depth=0.25ex, inner ysep=6pt]
    {  \beta_x^{-1}(\mcl{O}_{X,x}^\times)& \mcl{O}_{X,x}^\times & \{1\}\\
P & P^a  & P^a/\mcl{O}_{X,x}^\times\\};
    \path[arrow]
    (m-1-1) edge node[auto] {$ $} (m-1-2)
            edge node[auto,swap] {$$}  (m-2-1)
    (m-1-2) edge node[auto] {$$}(m-2-2)
edge node[auto] {$$}(m-1-3)
(m-2-1) edge node[auto] {$$} (m-2-2)
(m-1-3) edge (m-2-3)
(m-2-2) edge (m-2-3);
\draw ($(m-2-2.center)!.3!(m-2-1.west)$)+(0,.2) coordinate (A);
\draw ($(m-2-2.north)!.3!(m-1-2.south)$)+(-.2,0) -| (A);

\draw ($(m-2-3.center)!.3!(m-2-2.west)$)+(0,.2) coordinate (A);
\draw ($(m-2-3.north)!.3!(m-1-3.south)$)+(-.2,0) -| (A);
\end{tikzpicture}
\end{center}
from which we deduce
\begin{equation}
\label{eq:PaSharp}
P/\beta_x^{-1}(\mcl{O}_{X,x}^\times) \isoto P^a/\mcl{O}_{X,x}^\times \isoto \mcl{M}_{X,x}^\sharp.
\end{equation}
This also shows that if $\beta_x\colon P\to \mcl{M}_{X,x}$ is local, then $P^\sharp\isoto \mcl{M}_{X,x}^\sharp$.
Thus we can see charts as an easy way to describe the most informative part of the log structure.
\end{num}

\begin{prop}
Let $X$ be an fs log scheme and $x\in X$ a point.
Let $P = \mcl{M}_{X,x}^\sharp$.
Then there is a neighbourhood $U$ of $x$ and a chart $\beta\colon P\to \mcl{M}_X(U)$ such that the induced morphism $P\to \mcl{M}_{X,x}\to \mcl{M}_{X,x}^\sharp$ is the identity.
\end{prop}
\begin{proof}
Since $X$ is fs, $\mcl{M}_{X,x}^\sharp$ is fs and sharp.
Consequently $(\mcl{M}_{X,x}^\sharp)^\gp = \mcl{M}_{X,x}^\gp/\mcl{M}_{X,x}^\times$ is torsion-free and we can find a section $s\colon \mcl{M}_{X,x}^\gp/\mcl{M}_{X,x}^\times\to \mcl{M}_{X,x}^\gp$.
If $x\in \mcl{M}_{X,x}$ then $s([x]) - x\in \mcl{M}_{X,x}^\times$, hence $s([x])\in \mcl{M}_{X,x}$.
So we see that $s$ restricts to a section $\mcl{M}_{X,x}^\sharp\to \mcl{M}_{X,x}$.
Furthermore, $P^a = P\oplus \mcl{M}_{X,x}^\times\isoto \mcl{M}_{X,x}$.
Finally, proposition \cite[7.1.18(iv.c)]{148} allows us to extend this chart to a neighbourhood~$U$ of~$x$. 
\end{proof}


We present two examples of charts that will be extensively used in chapter~\ref{ch:log}.
\begin{eg}\label{prelim:Rdagger}
Let $R$ be a discrete valuation ring and $\pi$ a uniformizer.
We denote by $\logs{S}$ the scheme $S = \Spec R$ endowed with the divisorial log structure induced by its closed point.
Clearly, $\mcl{M}_{\logs{S}}(\logs{S}) = R\prive \{0\}$ and the morphism $\mbb{N}\to R\prive\{0\},\, 1\mapsto \pi$ is a chart for~$\logs{S}$.
We have indeed $\mbb{N}^a = \mbb{N}\oplus R^\times\isoto R\prive \{0\}$.
\end{eg}

\begin{eg}\label{prelim:logSncd}
We resume the notation of~\ref{prelim:Rdagger}.
Let $\mcl{X}\to S$ be an sncd $S$-scheme, i.e.\ a regular scheme of finite type over~$S$ whose special fibre $\mcl{X}_k  = \sum_{i\in I} N_i\, E_i$ is a divisor with strict normal crossings.
We endow $\mcl{X}$ with the divisorial log structure induced by~$\mcl{X}_k$.
For $J\sset I$ we write
\[E_J = \bigcap_J E_j \quad \text{ and } \quad E_J^\circ = E_J\prive \bigcup_{i\notin J} E_i,\]
where $E_\emptyset = \mcl{X}$.
Around each point $x\in E_J$ we can find an affine open~$U$ in~$\mcl{X}$ on which $\pi = u\prod_J x_j^{N_j}$ for $u$ a unit and where $V(x_j) = E_j$.
Then the morphism of monoids
\[\mbb{N}^{J}\to \mcl{M}_\mcl{X}(U), \,e_j\mapsto x_j\]
is a chart on~$U$.
\end{eg}


\begin{num}
We can easily see that the morphism $\mcl{X}\to S$ extends to a morphism of log schemes $\mcl{X}\to\logs{S}$.
We keep the notations of \ref{prelim:logSncd}.
If $u\neq 1$, we cannot find a morphism $\mbb{N}\to \mbb{N}^J$ making \eqref{dia:chartForMorphism} commutative.
But we can check that the morphism $\mbb{N}^J\oplus\mbb{Z}\to \mcl{M}_\mcl{X}(U)$ which maps the generator~$1$ of~$\mbb{Z}$ onto~$u$ is still a chart for the log structure on~$U$.
Hence defining
\[\mbb{N}\to \mbb{N}^J\oplus\mbb{Z},\quad 1\mapsto \bigl((N_j)_J,1\bigr)\]
yields a chart for the morphism $U\to\logs{S}$ of log schemes.

The lemma below shows that we have to slightly modify the morphism $\mbb{N}^J\to \mcl{M}_\mcl{X}(U)$ in order to get a chart of the form $\mbb{N}\to \mbb{N}^J$ for $U\to\logs{S}$.
\end{num}

\begin{lem}\label{SharpifiedChart}
Let $\phi\colon P\to Q$ be an fs chart for a morphism $X\to Y$ of log schemes.
Assume moreover that $P$ is sharp, $\phi$ is local and injective and the order of the torsion subgroup of $(Q^\sharp)^\gp/P^\gp$ is invertible on~$X$.
Then there is an étale cover $\wtl{X}\to X$ on which the morphism $P\to Q^\sharp$ extends to a chart for $\wtl{X}\to Y$.
\end{lem}
\begin{proof}
Since $\phi$ is injective and both $P$ and $Q$ are integral, $P^\gp\to Q^\gp$ is injective as well.
Furthermore, $P^\gp\to (Q^\sharp)^\gp$ is also injective because $\phi$ is local and $P$ is sharp.

Let $(\bar{e}_i)_{1\leq i \leq r}$ be a basis for $(Q^\sharp)^\gp$ adapted for $P^\gp$, i.e.\ we have integers $d_i\geq 1$ such that $(\bar{e}_i^{d_i})_{1\leq i\leq k}$ is a basis for $P^\gp$, for some $k\leq r$.
We denote by $x_i$ the element of $P^\gp$ corresponding to $\bar{e}_i^{d_i}$.
We lift $(\bar{e}_i)$ to elements $(e_i)$ of~$Q^\gp$, hence defining a section $s\colon (Q^\sharp)^\gp\to Q^\gp$.
Then we can find units $u_i\in Q^\times$ such that $\phi(x_i) = u_ie_i^{d_i}$ for each $1\leq i \leq k$.

Let $\beta\colon Q\to \mcl{M}_X(X)$ be the chart morphism.
We define $\wtl{X}$ as the étale cover obtained by adding a $d_i$-th root of $\beta(u_i)$ for each $d_i$, $1\leq i \leq k$ (they are invertible on $X$ by assumption).
Finally, the chart
\[Q^\sharp\to \mcl{M}_X(X), \, x = \prod_{i=1}^r \bar{e}_i^{n_i}\mapsto \beta(s(x))\, \prod_{i=1}^k \beta(u_i)^{n_i/d_i}\]
is well-defined, because $s(x)$ belongs to~$Q$ for $x\in Q^\sharp$ and we easily see that the diagram
\begin{center}
\begin{tikzpicture}
    \matrix (m) [matrix of math nodes, row sep=3em,
    column sep=3em, text height=1.5ex, text depth=0.25ex]
    { P & \mcl{M}_Y(Y) \\
Q^\sharp & \mcl{M}_X(X)\\};
    \path[arrow]
    (m-1-1) edge node[auto] {$ $} (m-1-2)
            edge node[auto,swap] {$$}  (m-2-1)
    (m-1-2) edge node[auto] {$$}(m-2-2)
(m-2-1) edge node[auto] {$$} (m-2-2);

\end{tikzpicture}
\end{center}
commutes.
\end{proof}

\subsection*{Fibred products of log schemes}

\begin{num}\label{logFibredProduct}
Let $f_1\colon X_1\to Y$ and $f_2\colon X_2\to Y$ be morphisms of log schemes.
Their \emph{fibred product} exists in the category of log schemes and is obtained by endowing the usual fibred product $X_1\times_Y X_2$ with the log structure associated to $p_1^{-1}\mcl{M}_{X_1}\oplus_{p_Y^{-1}\mcl{M}_Y} p_2^{-1}\mcl{M}_{X_2}$, where $p_1, p_2$ and $p_Y$ are the obvious projections.
Consequently, if $u_1\colon P\to Q_1$ and $u_2\colon P\to Q_2$ are charts for the morphisms $f_1$ and $f_2$ respectively, then the induced morphism $X_1\times_Y X_2\to \Spec \mbb{Z}[Q_1\oplus_P Q_2]$ is a chart as well.

We also immediately see that if $f\colon X\to Y$ is a morphism of log schemes, then the diagram
\begin{center}
\begin{tikzpicture}
    \matrix (m) [matrix of math nodes, row sep=3em,
    column sep=3em, text height=1.5ex, text depth=0.25ex]
    { X & Y \\
X^\circ & Y^\circ \\};
    \path[arrow]
    (m-1-1) edge node[auto] {$ f $} (m-1-2)
            edge node[auto,swap] {$  $} (m-2-1)
    (m-1-2) edge node[auto] {$  $} (m-2-2)
(m-2-1) edge node[auto,swap] {$ f $} (m-2-2);

\end{tikzpicture}\end{center}
is a fibred product of log schemes if and only if $f$ is strict.
\end{num}

\begin{num}
Interestingly, the category of fine (resp.\ fs) log schemes also admits fibred products.
If $f_1$ and $f_2$ are morphisms of fine log schemes and $u_1$, $u_2$ are fine charts as in \ref{logFibredProduct}, then
\[X_1\times_Y^\intg X_2\isDef (X_1\times_Y X_2)\times_{\mbb{Z}[Q_1\oplus_P Q_2]} \mbb{Z}[Q_1\oplus_P^\intg Q_2]\]
is the fibred product of~$X_1$ and~$X_2$ in the category of fine log schemes, where $\oplus^\intg$ indicates that we take the integral closure of the amalgamated sum.
It follows immediately that
\[X_1\times_Y^\intg X_2\to \Spec \mbb{Z}[Q_1\oplus_P^\intg Q_2]\]
is a chart.
When $u_1$ and $u_2$ are fs charts for morphisms $f_1$ and $f_2$ of fs log schemes, we have a similar expression for $X_1\times_Y^\fs X_2$, where $\oplus^\intg$ is replaced by $\oplus^\sat$, the saturation of the amalgamated sum.
Beware that those fibred products don't only have different log structures, they also have different underlying schemes.
\end{num}

\section{Log smoothness and log regularity}

Log differentials differ from classical differentials in that they may have at worst simple poles.
They are tightly related to the notion of log smoothness.
We formalize the notion as follows.

\begin{déf}\label{defOmegaLog}
Let $f\colon X\to S$ be a morphism of log schemes and let $\Omega_{X/S}$ be the classical sheaf of differentials of~$X^\circ/S^\circ$.
The \emph{module of log differentials}~$\Omega_{X/S}^{\logr}$ of $X/S$ is an $\mcl{O}_X$-module constructed as the quotient of $\Omega_{X/S}\oplus (\mcl{O}_X\otimes_{\mbb{Z}} \mcl{M}_X^\gp)$ by the submodule generated by
\begin{itemize}
\item $\bigl(\dif{\alpha_X(a)},0\bigr) - \bigl(0,\alpha_X(a)\otimes a\bigr)$, for $a\in\mcl{M}_X$,
\item $\bigl(0,1\otimes \phi(b)\bigr)$, for $b\in f^{-1}(\mcl{M}_S)$.
\end{itemize}
\end{déf}

\begin{num}
We have a canonical morphism of monoids
\[\dlog\colon \mcl{M}_X^\gp\to \Omega^\logr_{X/S},\, a\mapsto \bigl[(0,1\otimes a)\bigr]\]
that satisfies
\[\alpha_X(a)\cdot \dlog(a) = \dif{\alpha_X(a)}\]
for every $a\in\mcl{M}_X$.

If $X$ and $S$ are coherent log schemes (resp.\ and $f$ is locally of finite type), then $\Omega^\logr_{X/S}$ is quasi-coherent (resp.\ finitely generated) by \cite[IV.1.2.9]{137.2014}. 
Furthermore, we easily see that if the log structure on~$X$ is trivial, then $\Omega_{X/S}^\logr = \Omega_{X/S}$.
\end{num}

\begin{num} \label{prelim:computeOmegaFromChart}
If $P\to Q$ is a chart for the morphim $X\to S$, then $\Omega^\logr_{X/S}$ can be constructed as the quotient of
\[\Omega_{X/S}\oplus (\mcl{O}_X\otimes_{\mbb{Z}} Q^\gp)\]
by the same relations as in~\ref{defOmegaLog} for $a\in Q$ and $b\in P$ (see~\cite[1.7]{139}). 
\end{num}

\begin{prop}[{\cite[IV.1.3.1]{137.2014}}]\label{prelim:OmegaLogCartesian}
Consider the following cartesian diagram in the category of log schemes
\begin{center}
\begin{tikzpicture}
    \matrix (m) [matrix of math nodes, row sep=3em,
    column sep=3em, text height=1.5ex, text depth=0.25ex]
    { X' & S' \\
X & S \\};
    \path[arrow]
    (m-1-1) edge node[auto] {$ g $} (m-1-2)
            edge node[auto,swap] {$ p $} (m-2-1)
    (m-1-2) edge node[auto] {$ s $} (m-2-2)
(m-2-1) edge node[auto,swap] {$ f $} (m-2-2);
\draw ($(m-1-1.east)!.2!(m-1-2.west)$)+(0,-.12) coordinate (A);
\draw ($(m-1-1.south)!.2!(m-2-1.north)$)+(.12,0) -| (A);

\end{tikzpicture}\end{center}
Then the canonical morphism $p^*\Omega_{X/S}^\logr\to \Omega_{X'/S'}^\logr$ is an isomorphism.
In particular, if $X\to S$ is a strict morphism of log schemes, then $\Omega_{X/S}^\logr = \Omega_{X/S}$ (see \ref{logFibredProduct}).
\end{prop}
This result is even true in the category of fine (resp.\ fs) log schemes (\emph{loc.\ cit.}).

\begin{prop}[{\cite[IV.1.2.5]{137.2014}}]\label{OmegaLogPQ} 
Let $u\colon P\to Q$ be a morphism of monoids and $f\colon \Spec \mbb{Z}[Q]\to \Spec \mbb{Z}[P]$ the induced morphism of log schemes.
Then the canonical morphism
\[\mcl{O}_{\mbb{Z}[Q]}\otimes_\mbb{Z} Q^\gp/P^\gp \to \Omega_{\mbb{Z}[Q]/\mbb{Z}[P]}^\logr,\quad 1\otimes [a]\mapsto \dlog a\]
is an isomorphism.
\end{prop}

Log smoothness (and log étaleness) can be defined by means of log thickenings, much in the way that smoothness is (see for example \cite[\S 3]{139}). 
We will use here a more convenient definition due to Kato (cfr.\ \cite[3.5]{139}).

\begin{déf}\label{KatoCriterion}
A morphism $f\colon X\to Y$ of log schemes is \emph{log smooth} (resp.\ \emph{log étale}) if there is a chart $P\to Q$ for~$f$ étale-locally on $X$ and $Y$ such that
\begin{enumerate}
\item the kernel and the torsion part of the cokernel (resp.\ the kernel and the cokernel) of $P^\gp\to Q^\gp$ are finite groups whose order is invertible on $X$,
\item the induced morphism of schemes $h\colon X\to Y\times_{\mbb{Z}[P]} \mbb{Z}[Q]$ is smooth (resp.\ étale) in the classical sense.
\end{enumerate}
Note that the morphism $h$ is strict and that we can equivalently require $h$ étale in the definition of log smoothness (see \cite[3.6]{139}).

A log smooth morphism is locally of finite presentation (see \cite[3.3]{139}).
\end{déf}

\begin{eg}\label{prelim:sncdLogSmooth}
An sncd $S$-scheme $\mcl{X}$ endowed with the log structure defined in~\ref{prelim:logSncd} is log smooth over~$\logs{S}$ if the multiplicities $(N_i)_i$ of the components of the special fibre are all prime to the characteristic exponent of the residue field~$k$ of~$R$.
We resume the notation of~\ref{prelim:logSncd}.
Under this assumption, lemma~\ref{SharpifiedChart} ensures that we can find a chart
\[\mbb{N}\to \mbb{N}^J\]
for $\mcl{X}\to\logs{S}$ étale-locally on~$\mcl{X}$ that satisfies the first condition of definition~\ref{KatoCriterion}.
It remains to show that the induced morphism
\[h\colon \mcl{X}\to \Spec \bigl(R\otimes_{\mbb{Z}[\mbb{N}]}\mbb{Z}[\mbb{N}^J]\bigr)\defIs Y\]
is smooth.
The images $x_i$ of the basis vectors of $\mbb{N}^J$ are part of a regular system of parameters at each point $x\in\mcl{X}$.
Renumbering $(x_i)_J$ as $(x_1,\dots,x_r)$ and completing them to a regular system of parameters~$(x_i)_{1\leq i \leq n}$, we get a basis $(\dif{x_i})_{1\leq i \leq n}$ for $\Omega_{\mcl{X}/S}$, hence providing an étale morphism $\mcl{X}\to \mbb{A}^n_R$
which factors through an étale morphism
\[\mcl{X}\to \Spec \frac{R[x_1,\dots,x_n]}{(\pi-x_1^{N_1}\cdots x_r^{N_r})} = \mbb{A}^{n-r}_Y.\]
But $h$ is the composition of $\mcl{X}\to \mbb{A}^{n-r}_Y\to Y$, thus it is smooth.
\end{eg}

\begin{prop}[{\cite[IV.3.1.2]{137.2014}}]\label{prelim:logSmoothSmooth} 
Log smooth and log étale morphisms are stable under composition and base change, and a strict morphism of log schemes is log smooth (resp.\ log étale) if and only if it is smooth (resp.\ étale).
\end{prop}

\begin{déf}\label{intgMorphSchemes}
A morphism of integral log schemes $f\colon X\to Y$ is called \emph{integral} if for every $x\in X$ the morphism of monoids $\mcl{M}_{Y,f(x)}\to \mcl{M}_{X,x}$ is integral.
\end{déf}

\begin{prop}[{\cite[II.3.8 and II.3.15]{135}}] \label{prelim:OmegaLogLocFree}
Let $f\colon X\to Y$ be a log smooth morphism and denote by $\dim_x f$ the relative dimension of~$f$ at a point~$x$.
The sheaf $\Omega_{X/Y}^\logr$ is locally free of finite rank.

If furthermore $X$ and $Y$ are fine, then we have at every point $x\in X$
\[\rk \Omega_{X/Y,x}^\logr \leq \dim_x f,\]
with equality if $f$ is integral.
\end{prop}

\begin{prop}[{\cite[IV.3.2.4]{137.2014}}]\label{OmegaLogÉtale} 
Let $f\colon X\to Y$ and $g\colon Y\to S$ be morphisms of log schemes.
If $f$ is log étale, then the canonical morphism
\[f^*\Omega_{Y/S}^\logr\to \Omega_{X/S}^\logr\]
is an isomorphism.
\end{prop}

\begin{cor}
Let $X\to Y$ be a log smooth morphism of log schemes and $P\to Q$ a chart as in \ref{KatoCriterion}.
Then
\[\Omega_{X/Y}^\logr\isoto\mcl{O}_X\otimes_{\mbb{Z}} h^*(Q^\gp/P^\gp).\]
\end{cor}
\begin{proof}
This is immediate from \ref{OmegaLogPQ}, \ref{prelim:OmegaLogCartesian} and \ref{OmegaLogÉtale}.
\end{proof}

\begin{prop}[{\cite[IV.2.3.1 and 3.2.3]{137.2014}}]\label{prelim:OmegaLogExact} 
Let $f\colon X\to Y$ and $g\colon Y\to S$ be morphisms of log schemes.
The following sequence is exact:
\[f^*\Omega^\logr_{Y/S}\to \Omega^\logr_{X/S}\to \Omega^\logr_{X/Y}\to 0.\]
Furthermore, if $f$ is log smooth, then it is injective and locally split.
Conversely, if $g\circ f$ is log smooth and the sequence is injective and locally split, then $f$ is log smooth.
\end{prop}

\subsection*{Log regularity}
We now collect some results on log regular schemes; the most important construction is the fan associated to a log regular scheme in~\ref{FanInducedByLogReg}. The main reference for this part is \cite{143}.

\begin{déf}
A locally Noetherian and locally fs log scheme $X$ is said to be \emph{log regular} if for every $x\in X$
\begin{enumerate}
\item the ring $\mcl{O}_{X,x}/\mcl{M}_{X,x}^+\mcl{O}_{X,x}$ is regular;
\item $\dim \mcl{O}_{X,x} = \dim \mcl{O}_{X,x}/\mcl{M}_{X,x}^+\mcl{O}_{X,x} + \dim \mcl{M}_{X,x}$.
\end{enumerate}
\end{déf}

If the log structure of~$X$ is trivial at~$x$, then log regularity at~$x$ boils down to regularity.

\begin{eg}
Let $\logs{S}$ be the spectrum of a discrete valuation ring~$R$ endowed with the log structure described in \ref{prelim:Rdagger}, and $s$ its closed point.
Then $\mcl{M}^+_{\logs{S},s}\mcl{O}_{\logs{S},s} =\pi R$ and $\dim \mcl{M}_{\logs{S},s} = \dim \mbb{N} = 1$, from which we easily conclude that $\logs{S}$ is a log regular scheme.
\end{eg}

\begin{prop}[{\cite[6.2]{143}}]\label{prelim:logRegFlat}
Let $x$ be a point on a locally Noetherian log scheme~$X$ and $P\to \mcl{M}_X(X)$ an fs chart local at~$x$, with $P$ sharp.
Assume moreover that $\mcl{O}_{X,x}$ contains a field~$k$.
Then $X$ is log regular at~$x$ if and only if $\mcl{O}_{X,x}/\mcl{M}^+_{X,x}\mcl{O}_{X,x}$ is a regular local ring and the induced map $k[P]\to \mcl{O}_{X,x}$ is flat.
\end{prop}

\begin{prop}[{\cite[8.2]{143}}]\label{prelim:logLisseSurLogReg}
Let $X\to Y$ be a log smooth morphism between locally fs log schemes.
If $Y$ is log regular, so is $X$.
\end{prop}

\begin{prop}[{\cite[7.3]{143}}]\label{pOXisPrime}
Let $x$ be a point on a log regular scheme~$X$, and $\mfk{p}$ a prime ideal of $\mcl{M}_{X,x}$.
Then $\mfk{p}\mcl{O}_{X,x}$ is a prime ideal of the same height as~$\mfk{p}$.
\end{prop}

\begin{num}\label{FanInducedByLogReg}
Let $X$ be a log regular scheme. We set
\[F = F(X) \isDef \bigl\{x\in X\bigm| \mcl{M}_{X,x}^+\mcl{O}_{X,x} = \mfk{m}_x\bigr\}.\]
If $i$ denotes the inclusion of $F$ in  $X$, the sharp monoidal space $\bigl(F,(i^{-1}\mcl{M}_X)^\sharp\bigr)$ has the structure of a locally fs fan and is called the \emph{fan associated to $X$}.
It is finite if $X$ is quasi-compact.
More precisely we have
\end{num}

\begin{prop}[{\cite[10.1]{143}}]\label{chartForFan}
Let $X$ be a log regular scheme and $P\to \mcl{M}_X(X)$ an fs chart. Let $x\in X$ and denote by $\mfk{p}$ the inverse image of $\mcl{M}_{X,x}^+$ in~$P$.
Then there is an affine open neighbourhood $V$ of~$x$ such that
\[F(X)\cap V \isoto (\Spec P_\mfk{p})^\sharp.\]
\end{prop}

\begin{eg}
The fan associated to the spectrum of a discrete valuation ring~$\logs{S}$ is isomorphic to $(\Spec \mbb{N})^\sharp$ as we can see from the chart given in~\ref{prelim:Rdagger}.
\end{eg}

Fans of log regular schemes will prove immensely useful in chapter~\ref{ch:log}.

\begin{num}
We have a morphism of sharp monoidal spaces
\[\pi\colon X^\sharp\to F, \quad x\mapsto \mcl{M}_{X,x}^+\mcl{O}_{X,x}.\]
It satisfies $\pi\circ i = \id_F$.

The construction of $F(X)$ is functorial: if $f\colon X\to Y$ is a morphism between log regular schemes, then we get a morphism of fans $F(X)\to F(Y)$ such that the following diagram commutes
\begin{equation}\label{prelim:diag_fanMorph}
\begin{tikzpicture}[baseline=(current  bounding  box.center)]
    \matrix (m) [matrix of math nodes, row sep=3em,
    column sep=3em, text height=1.5ex, text depth=0.25ex]
    { X & Y \\
F(X) & F(Y) \\};
    \path[arrow]
    (m-1-1) edge node[auto] {$f$} (m-1-2)
            edge node[auto,swap] {$\pi_X$}  (m-2-1)
    (m-1-2) edge node[auto] {$\pi_Y$}(m-2-2)
(m-2-1) edge node[auto] {$F(f)$} (m-2-2);

\end{tikzpicture}\end{equation}
where the bottom arrow is given by
\[x\mapsto \pi_Y\circ f\circ i_X(x).\]
This latter morphism is induced by any chart $P\to Q$ for~$f$.
\end{num}

\begin{num}
For every $t\in F$, we set $U(t) = \pi^{-1}(t)$, which is a locally closed subscheme of~$X$.
The family $(U(t))_{t\in F}$ form the \emph{log stratification} of~$X$, and we will denote by $E(t)$ the schematic closure of $U(t)$ endowed with its reduced structure.
It is equal to the closure of $i(t)$ in~$X$.
We can compute the log stratification of~$X$ via charts:
\end{num}

\begin{prop}\label{prelim:calculEt}
Let $X$ be a log regular scheme of fan~$F$ and $\beta\colon P\to \mcl{M}_X(X)$ an fs chart.
For each prime ideal $\mfk{p}$ of $P$ we set
\[T(\mfk{p}) = \bigl\{t\in F\bigm| \beta^{-1}\mcl{M}_{X,i(t)}^+ = \mfk{p}\bigr\}.\]
Then $\bigcup_{T(\mfk{p})} E(t) = X \times_{\mbb{Z}[P]} \mbb{Z}\langle P/\mfk{p}\rangle$.
In particular $\bigcup_{T(P^+)} U(t) = X\times_{\mbb{Z}[P]} \mbb{Z}\langle P/P^+\rangle$.
\end{prop}
\begin{proof}
Since for every $x\in X$ the morphism $\mcl{M}_{X,x}\to \mcl{O}_{X,x}$ is local we have
\[\mfk{p}\mcl{O}_{X,x}\sset \mfk{m}_x \iff \mfk{p}\mcl{M}_{X,x}\sset \mcl{M}_{X,x}^+\]
In this case it follows from \eqref{eq:PaSharp} that $\mfk{p}\mcl{M}_{X,x}$ is a prime ideal of~$\mcl{M}_{X,x}$ whose inverse image in~$P$ is~$\mfk{p}$.
Let $t\in F$ be the point corresponding to the prime ideal $\mfk{p}\mcl{O}_{X,x}$ (see \ref{pOXisPrime}).
Then $t\in T(\mfk{p})$ and $x\in E(t)$.
\end{proof}

Recall that we have a notion of subdivision of a fan, see~\ref{déf:subdivision}.

\begin{prop}\label{prelim:inverseSubdiv}
Let $X$ be a log regular scheme and $F$ its fan. Let $\vphi\colon F'\to F$ be a subdivision of $F$ with $F'$ locally fs.
Then there exists a unique log étale morphism $\vphi^*X\to X$ up to isomorphism such that $F'\to F$ is isomorphic to $F(\vphi^*X)\to F$.
In particular, $\vphi^*X$ is log regular.
\end{prop}
\begin{proof}
See \cite[7.6.14 and 7.6.18]{148} for a complete proof. The log scheme $\vphi^*(X)$ is constructed locally as follows. 

Let $x\in X$ and $P\to \mcl{M}_X(V)$ be an fs chart in an affine neighbourhood~$V$ of~$x$ such that $F\cap V = (\Spec P)^\sharp$.
Let $t\in F'$ with $\vphi(t) = \pi(x)$ and denote by~$Q$ the fibred product of monoids
\[M_{F',t}\times_{M_{F',t}^\gp} P^\gp.\]
Then $\vphi^*X$ is obtained by gluing all the $V\times_{\mbb{Z}[P]}\mbb{Z}[Q]$.

Note that $(\Spec M_{F',t})^\sharp\isoto (\Spec Q)^\sharp$ by \cite[7.6.12]{148} and $Q^\gp = P^\gp$.
\end{proof}

\cleardoublepage


%

%
  \graphicspath{{\chaptersdir/log/image/}}%
  \chapter{Computing zeta functions on log smooth models}\label{ch:log}

\section{Introduction}
Let $f$ be a polynomial in $k[x_1,\dots,x_n]$, where $k$ is a field of characteristic zero.
The motivic zeta function $Z_f(T)$ of~$f$ was first introduced by Denef and Loeser in \cite{74} in analogy with the $p$-adic setting, where Igusa local zeta functions of polynomials had been long studied for their use in understanding the asymptotic behaviour of the number of congruence solutions of $f = 0$ modulo powers of~$p$.
They naturally stated a motivic version of the monodromy conjecture holding sway in the $p$-adic world.
This conjecture relates the poles of $Z_f$ to eigenvalues of the monodromy action on the complex of nearby cycles on $f^{-1}(0)$ (see ibid.\ \S2.4).
Denef and Loeser proved an explicit formula for~$Z_f$ involving a resolution of singularities for~$f$ and yielding naturally a set of candidate poles for $Z_f$ (see theorem~\ref{intro:motzetaDL}).
Yet, a major issue in proving the monodromy conjecture is to determine which of the apparent poles of $Z_f$ are actual poles; most of the proofs in special cases rely on a detailed study of the geometry of the exceptional locus of resolutions of singularities for~$f$.

In \cite{106} Nicaise and Sebag introduced the volume Poincaré series $S(\mcl{X},\omega;T)$ of a pair $(\mcl{X},\omega)$ consisting of a generically smooth stft (for separated and topologically of finite type) formal scheme $\mcl{X}$ over the formal spectrum of a complete discrete valuation ring~$R$ and a volume form on its generic fibre~$\mcl{X}_K$.
They then showed how to express $Z_f$ as such a series, yielding a new interpretation of the function.
This new way of considering $Z_f$ allowed Halle and Nicaise in \cite{146} to associate a motivic zeta function $Z_X$ to a Calabi-Yau variety $X$ over the fraction field of~$R$ and state an analogue of the monodromy conjecture, which they could prove when $X$ is a tamely ramified abelian variety (see definition~\ref{intro:monodPropty}).

In this chapter we will mostly stick to the algebraic setting and show how to compute $S(\mcl{X},\omega;T)$ when $\mcl{X}$ is a generically smooth log smooth $S$-scheme.
By adding a suitable structure to schemes, logarithmic geometry allows to handle so-called log smooth schemes as if they were smooth.
In particular the sheaf of log differentials is locally free and we see how a volume form nicely induces a divisor on~$\mcl{X}$ by considering this sheaf.
Another feature of log smooth schemes is that their fans, as defined by Kato in \cite{143}, can be used to exhibit a desingularization of the scheme.
In this process of desingularizing, fake poles are intoduced in the expression of $S(\mcl{X},\omega;T)$, so that our formula, depending directly on the fan, reduces substantially the set of candidate poles.
This paves the way to the study of motivic zeta functions for degenerations of Calabi-Yau varieties that appear in the Gross-Siebert programme on mirror symmetry, and for which log smooth models are naturally constructed (see for example \cite{170}).

The main results of this chapter were announced in~\cite{186}.
In section~\ref{sec:Nmonoids}, we study properties of $\mbb{N}$-monoids more in depth.
Section~\ref{sec:logRegOverDVR} applies those results to the structure of the log stratification of a log regular scheme under subdivisions of fans. We conclude by a proof of a formula for $S(\mcl{X},\omega;T)$ with $\mcl{X}$ log smooth in section~\ref{sec:formulaForPoincaré}.
Finally, we show in section~\ref{log:sectionNewton} how to use our formula to recover a result of Guibert (see~\cite{165}) expressing the motivic zeta function of a polynomial that is nondegenerate with respect to its Newton polyhedron and section~\ref{log:twoVariables} shows how our formula applies to yield precisely the right poles for~$Z_f(T)$ when $f$ is a polynomial in two variables.

In this chapter we work throughout with Zariski sites. Every log scheme considered will be assumed locally Noetherian and fine and saturated, abbreviated by~\emph{fs}.

\section{$\mbb{N}$-monoids}
\label{sec:Nmonoids}

We will make use of two more advanced properties of monoids in our proofs.

\begin{num}\label{val-exact}
An integral monoid $Q$ is called \emph{valuative} if for any $x\in Q^\gp$, either $x\in Q$ or $x^{-1}\in Q$.
A submonoid $M$ of an integral monoid $Q$ is \emph{exact} if $M^\gp\cap Q = M$.
In particular we see that any valuative submonoid~$M$ of~$Q$ is exact if $Q^\times\sset M$.
\end{num}

\begin{déf}
Let $Q$ be an integral monoid. For a submonoid $M$ of $Q$, we set
\[M^{\sat,Q} \isDef \bigl\{x\in Q^\sat\bigm| \exists \alpha\geq 1 \colon x^\alpha \in M\bigr\}.\]
We easily see that $M^{\sat,Q}$ is a face of $Q^\sat$ if $M$ is a face of $Q$.
\end{déf}

\begin{prop}\label{facesQsat}
Let $Q$ be an integral monoid.
The correspondence $F\mapsto F^{\sat,Q}$ induces a bijection between the faces of $Q$ and the faces of $Q^\sat$.
In particular $(Q^\sat)^\times \cap Q = Q^\times.$
\end{prop}
\begin{proof}
We first show $F^{\sat,Q}\cap Q = F$.
Let $x\in F^{\sat,Q}\cap Q$. Then $x^\alpha\in F$ for some $\alpha\geq 1$. But since $x\in Q$, this implies $x\in F$.

Let $G$ be a face of $Q^\sat$. We show $(G\cap Q)^{\sat,Q} = G$. Let $x\in G$ and $\alpha\geq 1$ such that $x^\alpha\in Q$. Then $x^\alpha\in G\cap Q$ and by definition $x\in (G\cap Q)^{\sat,Q}$.
\end{proof}

\begin{num}\label{log:defEpi}
Let $\phi\colon \mbb{N}\to P$ be an $\mbb{N}$-monoid and $d$ a positive integer. We set $e_\pi\isDef\phi(1)$ and $P(d)\isDef P\oplus_\mbb{N} \mbb{N}_d$, where $\mbb{N}\to\mbb{N}_d$ denotes multiplication by $d$.
If $P$ is integral, then $\phi$ is an integral morphism since $\mbb{N}$ is valuative (see \cite[I.4.5.3(5)]{137.2014}), so that $P(d)$ is integral as well.\label{NmorphismIntegral}
Since the internal law on $\mbb{N}$ is additive, we will often write the laws of $P$ and $P(d)$ additively as well.
\end{num}

\begin{prop}\label{PiIsN}
Let $\mbb{N}\to P$ be a local morphism of fs monoids. Then the submonoid 
\[\Pi = \bigl\{x\bigm| \alpha x = k[e_\pi] \text{ for some } \alpha\geq 1, k\geq 0\bigr\}\]
of $P^\sharp$ is isomorphic to $\mbb{N}$ and $\Pi^\sat = \Pi^{\sat,P^\sharp}= \Pi$.
\end{prop}
\begin{proof}
We first show that $\Pi$ is valuative.
By construction of the groupification, every $x\in \Pi^\gp$ can be written as $x = x_1-x_2$ with $x_1,x_2\in \Pi$.
Let $\alpha_i\geq 1$ be such that $\alpha_ix_i = k_i[e_\pi]$ for some $k_i\geq 0$, $i=1,2$.
Then
\[(\alpha_1\alpha_2) x = (\alpha_2k_1-\alpha_1k_2) [e_\pi],\]
so that either $x$ or $-x$ belongs to~$\Pi$.

Since $\Pi$ is valuative, it is an exact submonoid of $P^\sharp$, hence finitely generated (see \cite[I.2.1.17(2)]{137.2014}).
We conclude by \cite[I.2.4.2]{137.2014} that $\Pi\isoto \mbb{N}$ since $\Pi$ is nonzero by localness of $\mbb{N}\to P$.
The second assertion follows from the fact that $P$ is saturated.
\end{proof}

\begin{déf}
In the setting of \ref{PiIsN}, we will denote by $\bar{e}_\pi$ the generator of $\Pi$ and call the integer $m\geq 1$ such that $[e_\pi] = m\bar{e}_\pi$ the \emph{root index} of $\mbb{N}\to P$.
\end{déf}

\begin{rmq}\label{ePiFreeFactor}
The root index of a local morphism $\mbb{N}\to P$ is 1 if and only if $[e_\pi]$ is nonzero and generates a free factor of the free group $(P^\sharp)^\gp$.
More generally, the root index corresponds to the order of the torsion part of $(P^\sharp)^\gp/\mbb{Z}$.
\end{rmq}

\begin{prop}\label{unitsInPd}
Let $\mbb{N}\to P$ be a local morphism. Then the canonical morphism $P\to P(d)$ restricts to an isomorphism $P^\times\toisoto P(d)^\times$.
\end{prop}
\begin{proof}
First note that $(x,0)\in P(d)^\times$ if and only if $x\in P^\times$. In particular $(0,1)\in P(d)^\times$ if and only if $e_\pi\in P^\times$.
Now if $(v,t)\in P(d)^\times$, then both $(v,0)$ and $(0,t)$ belong to $P(d)^\times$ and since $\phi$ is local, $e_\pi\notin P^\times$, which implies $t=0$.
\end{proof}

\begin{prop}\label{PtoQsatSharp}
Let $\mbb{N}\to P$ be a local morphism of root index $m$ between fs monoids. If $\gcd(d,m) = 1$, then the morphism $P\to P(d)^\sat$ restricts to an isomorphism $P^\times\toisoto P(d)^{\sat,\times}$.
\end{prop}
\begin{proof}
By \ref{unitsInPd}, $P^\times \isoto P(d)^\times$ so that we only have to show $P(d)^\times\isoto P(d)^{\sat,\times}$.
Let $x\in P(d)^{\sat,\times}$ and $\alpha\geq 1$ such that $\alpha x\in P(d)^\times$ (see \ref{facesQsat}).
We may write $x$ as $(v,t)$ with $0\leq t< d$ and $v \in P^\gp$.
Since $\alpha x\in P(d)^{\times}$, $\alpha t = qd$ for some $q\geq 0$ and
\[\alpha x= (\alpha v + qe_\pi,0)\]
with $\alpha v + qe_\pi\in P^\times$.
We then have in $(P^{\sharp})^\gp$ the equality $\alpha [-v] = q [e_\pi]$. But then $[-v]\in \Pi$ so that $[-v] = s\bar{e}_\pi$ for some $s\geq 0$ and $\alpha s = mq$. Since $\gcd(d,m) = 1$ and $\alpha$ also divides $qd$, $\alpha$ is a divisor of $q$. But since $q<\alpha$ (because $t<d$) we necessarily have $s = 0 = q$.
We eventually get $v\in P^\times$ and $t = 0$, i.e.\ $x\in P(d)^\times$.
\end{proof}

\begin{prop}\label{maxIdealOfPd}
Let $\mbb{N}\to P$ be a local morphism of root index $m$ between fs monoids. If $d$ divides $m$, then the maximal ideal of $P(d)^\sat$ is generated by~$P^+$.
\end{prop}
\begin{proof}
Let $x = (v,t)\in P(d)^{\sat,+}$ and $\alpha\geq 1$ such that $\alpha x\in P(d)^+$. If we write $\alpha t = qd+s$ with $0\leq s<d$ we have
\[\alpha x = (\alpha v,\alpha t) = (\alpha v + qe_\pi,s)\] with $\alpha v + qe_\pi\in P$.
We are going to show that $x$ can be written as the sum of an element of $P^+$ and a unit of $P(d)^{\sat}$.
Let $\bar{e}\in P^+$ be such that $[\bar{e}] = \bar{e}_\pi\in P^\sharp$. Then there is $u\in P^\times$ such that $e_\pi = m\bar{e} + u$.
We have in $P^\gp$
\begin{align*}
\alpha (v+ \tfrac{m}{d} t\bar{e}) &= \alpha v + \tfrac{m}{d}\bar{e} (qd+s)\\
&= \alpha v + q e_\pi - qu +\tfrac{m}{d}s\bar{e},
\end{align*}
which belongs to $P^+$ because either $\alpha v+ qe_\pi\in P^+$ or $s\geq 1$ by \ref{unitsInPd}.
Since $P$ is saturated, we get $v + \tfrac{m}{d} t\bar{e}\in P^+$ and it only remains to show that $x - (v+ \tfrac{m}{d} t\bar{e},0)$ is a unit of~$P(d)^\sat$.
But
\begin{align*}
d\,\bigl( x - (v+ \tfrac{m}{d} t\bar{e},0) \bigr) &= d\,(-\tfrac{m}{d} t\bar{e},t)\\
&= (-mt \bar{e},td)\\
&= (-mt \bar{e} + te_\pi,0)\\
&= (tu,0),
\end{align*}
is a unit in $P(d)$. Hence $x - (v+ \tfrac{m}{d} t\bar{e},0)\in P(d)^{\sat,\times}$.
\end{proof}

The following proposition will be useful for computations (see \ref{étaleCover} and \ref{Xtau1}).

\begin{prop}\label{saturation}
Let $P$ be an fs monoid and
\[\mbb{N}\to P\oplus F,\quad 1\mapsto (e,\omega)\]
a morphism, where $F$ is a free group. Let $d$ be an integer such that there is some $\bar{e}\in P$ with $d\bar{e} = e$ and assume that there are elements $(\omega_i)_I$ for~$F$ such that $(\omega,(\omega_i))$ is a basis of $F$.
Then the morphism \[\phi\colon P\oplus F\to P\oplus F,\, \bigl(x, a_0\omega + \sum a_i\omega_i\bigr)\mapsto \bigl(x, da_0\omega + \sum a_i\omega_i\bigr)\]
factors through an isomorphism
\[(P\oplus F) \oplus_{\mbb{N}}^\sat \mbb{N}_d \isoto P\oplus F.\]
\end{prop}
\begin{proof}
Let $Q= (P\oplus F) \oplus_{\mbb{N}} \mbb{N}_d $.
The morphism $\phi$ together with
\[\mbb{N}\to P\oplus F,\quad 1\mapsto (\bar{e},\omega)\]
induces a unique morphism $Q\to P\oplus F$.
Since the latter is saturated, we have a unique factorisation through $Q^\sat$.

On the other hand, the morphism
\[P\oplus F\to Q^\sat,\, \bigl(y,a_0\omega + \sum a_i\omega_i\bigr)\mapsto \bigl((y-a_0\bar{e},\sum a_i\omega_i),a_0\bigr),\]
is well defined since 
\begin{align*}
d\, \bigl((y-a_0\bar{e},\sum a_i\omega_i),a_0\bigr) &= \bigl((dy-a_0e,d\sum a_i\omega_i),da_0\bigr)\\
&= \bigl((dy,a_0\omega + d\sum a_i\omega_i),0\bigr) \in Q.
\end{align*}
We may now easily check that those morphisms are each other inverse.
\end{proof}

\begin{cor}\label{rootPdsat}
Let $\mbb{N}\to P$ be a local morphism of root index~$m$ with $P$ fs. If $d$ divides~$m$, then $\mbb{N}_d\to P(d)^\sat$ has root index~$\tfrac{m}{d}$.
\end{cor}
\begin{proof}
Let $\bar{e}_\pi\in P^\sharp$ be such that $m[e_\pi] = \bar{e}_\pi$.
We can complete $(\bar{e}_\pi)$ to a basis of $(P^\sharp)^\gp$.
In particular we get a section $P^\sharp\to P$ and $\mbb{N}\to P$ factors as
\[\mbb{N}\to P^\sharp \oplus P^\times, \quad 1\mapsto ([e_\pi],u),\]
for some unit $u\in P^\times$.
The statement is now easy to deduce from the proof of~\ref{saturation}.
\end{proof}

\section{Log regular schemes over a discrete valuation ring}
\label{sec:logRegOverDVR}

\begin{num}
Let $R$ be a discrete valuation ring and fix a uniformizer $\pi$.
We denote by~$K$ its fraction field and by~$k$ its residue field.
Recall that we denote by $\logs{S}$ the scheme $S = \Spec R$ endowed with the log structure induced by its closed point, as described in \ref{prelim:Rdagger}.

Let $d$ be a positive integer.
We set $R(d) \isDef R[T]/(T^d-\pi)$ which is a totally ramified extension of $R$ of degree $d$, and $S(d) = \Spec R(d)$.
The chart $\mbb{N}\to\mbb{N}_d$ of $R\to R(d)$ induces an isomorphism
\begin{equation}\label{RdTenseur}
R(d)\isoto R\otimes_{\mbb{Z}[\mbb{N}]} \mbb{Z}[\mbb{N}_d].\end{equation}
\end{num}

\begin{num}
Let $\mcl{X}$ be a log regular scheme over $\logs{S}$ and denote its fan by $F$.
By functoriality, we have a morphism of fans
\[F\to F(\logs{S}) \isoto (\Spec \mbb{N})^\sharp\]
that gives a structure of $\mbb{N}$-fan to~$F$.
We will denote by $F_\eta$ and $F_s$ its fibres over the points $\emptyset$ and $\mbb{N}^+$ respectively.
It then follows from \eqref{prelim:diag_fanMorph} that $\pi^{-1}_\mcl{X}(F_\eta) = \mcl{X}_K$ and $\pi^{-1}_\mcl{X}(F_s) = \mcl{X}_k$.

If we denote by $e_\pi$ the image of 1 under the induced morphism $\mbb{N}\to M_F(F)$ (compare with \ref{log:defEpi}), then
$F_\eta = \{t\in F\mid e_\pi \in M_{F,t}^\times\}$ and $F_s = \{t\in F\mid e_\pi \in M_{F,t}^+\}.$
\end{num}

\begin{déf}
We say that a subdivison $\vphi\colon F'\to F$ of $\mbb{N}$-fans \emph{preserves the horizontal part} if $\vphi^{-1}(F_\eta) \isoto F'_\eta$.
\end{déf}

\begin{rmq}
If $\mcl{X}\to \logs{S}$ is a log regular $\logs{S}$-scheme and $\vphi\colon F'\to F$ is a subdivision of its fan~$F$ preserving the horizontal part, then we have an isomorphism of generic fibres $(\vphi^*\mcl{X})_K\isoto \mcl{X}_K$, as we can see from the proof of~\ref{prelim:inverseSubdiv}.
\end{rmq}

\begin{num}\label{XdFibredProd}
Let $\mcl{X}$ be a log regular scheme over $\logs{S}$ and $d$ a positive integer.
We will denote by $\mcl{X}(d)$ and $\mcl{X}(d)^\fs$ the fibred product $\mcl{X}\times_\logs{S} \logs{S(d)}$ in the category of log schemes and fs log schemes respectively.
By \ref{RdTenseur}, we have
\[\mcl{X}(d) = \mcl{X}\times_{\mbb{Z}[P]} \mbb{Z}[P(d)] \text{ and } \mcl{X}(d)^\fs = \mcl{X}\times_{\mbb{Z}[P]} \mbb{Z}[P(d)^\sat]\]
for any chart $\mbb{N}\to P$ for $\mcl{X}\to \logs{S}$. 
\end{num}

\begin{num}
Let $t\in F_s$ and $m = m(t)$ be the root index of $\mbb{N}\to M_{F,t}$.
For every divisor $d$ of $m$, we will denote by $U_d(t)$ the inverse image of $U(t)$ in $\mcl{X}(d)^\fs$.
We also set $\wtl{U}(t) \isDef U_m(t)$.
\end{num}

\begin{prop}\label{calculUdt}
Let $\mcl{X}\to\logs{S}$ be a log regular scheme over~$\logs{S}$ and consider a chart $\mbb{N}\to P$ with $P$ fs.
We keep the notations of \ref{prelim:calculEt}.
Then every point of $T(P^+)$ belongs to $F_s$ and has the same root index~$m$, and for every~$d$ dividing~$m$ we have
\[\bigcup_{T(P^+)} U_d(t)= \mcl{X}\times_{\mbb{Z}[P]} \mbb{Z}\langle P(d)^{\sat}/P(d)^{\sat,+}\rangle.\]
\end{prop}
\begin{proof}
For any $t\in T(P^+)$ we have a diagram of monoids
\begin{center}
\begin{tikzpicture}
    \matrix (m) [matrix of math nodes, row sep=3em,
    column sep=3em, text height=1.5ex, text depth=0.25ex]
    { \mbb{N} & R \prive\{0\}\\
P & \mcl{M}_{\mcl{X},i(t)}\\};
    \path[arrow]
    (m-1-1) edge node[auto] {$ $} (m-1-2)
            edge node[auto,swap] {$$}  (m-2-1)
    (m-1-2) edge node[auto] {$$}(m-2-2)
(m-2-1) edge node[auto] {$\beta$} (m-2-2);

\end{tikzpicture}\end{center}
Since the top and right arrows are local morphisms, $\beta^{-1}(\mcl{M}_{\mcl{X},i(t)}^+) = P^+$ implies that $e_\pi\in P^+$. The statement over the root indices is clear since $P^\sharp = M_{F,t}$ for every $t\in T(P^+)$.
By \ref{prelim:calculEt} and \ref{XdFibredProd},
\[\bigcup_{T(P^+)} U_d(t) = \mcl{X}\times_{\mbb{Z}[P]} \mbb{Z}\langle P(d)^\sat /I\rangle,\]
where $I$ is the ideal of $P(d)^\sat$ generated by $P^+$, which is $P(d)^{\sat,+}$ by \ref{maxIdealOfPd}.
\end{proof}

The following lemma will help us in proving proposition \ref{FibresOnUt}, which is the main result of this section.

\begin{lem}\label{integralPushout}
Consider the following amalgamated sum of monoids with $P,Q_1$ and $Q_2$ integral:
\begin{center}
\begin{tikzpicture}
    \matrix (m) [matrix of math nodes, row sep=3em,
    column sep=3em, text height=1.5ex, text depth=0.25ex, inner ysep=6pt]
    {  P& Q_1 \\
Q_2 & Q_1\oplus_P Q_2 \\};
    \path[arrow]
    (m-1-1) edge node[auto] {$ $} (m-1-2)
            edge node[auto,swap] {$$}  (m-2-1)
    (m-1-2) edge node[auto] {$$}(m-2-2)
(m-2-1) edge node[auto] {$$} (m-2-2);
\draw ($(m-2-2.center)!.3!(m-2-1.west)$)+(0,.2) coordinate (A);
\draw ($(m-2-2.north)!.3!(m-1-2.south)$)+(-.2,0) -| (A);
\end{tikzpicture}\end{center}
If $P$ injects into $Q_1$, then $Q_2$ injects into $Q_1\oplus_P Q_2$.
\end{lem}
\begin{proof}
Since $Q_2$ is integral, we only have to show that $Q_2^\gp$ injects into $(Q_1\oplus_P Q_2)^\gp$.
We have the following amalgamated sum in the category $\catAb$ of abelian groups (see \ref{prelim:gpPreservesColim}):
\begin{center}
\begin{tikzpicture}
    \matrix (m) [matrix of math nodes, row sep=3em,
    column sep=3em, text height=1.5ex, text depth=0.25ex, inner ysep=6pt]
    {  P^\gp& Q_1^\gp \\
Q_2^\gp & (Q_1\oplus_P Q_2)^\gp \\};
    \path[arrow]
    (m-1-1) edge node[auto] {$ $} (m-1-2)
            edge node[auto,swap] {$$}  (m-2-1)
    (m-1-2) edge node[auto] {$$}(m-2-2)
(m-2-1) edge node[auto] {$$} (m-2-2);
\draw ($(m-2-2.center)!.3!(m-2-1.west)$)+(0,.2) coordinate (A);
\draw ($(m-2-2.north)!.3!(m-1-2.south)$)+(-.2,0) -| (A);
\end{tikzpicture}\end{center}
Since $P$ and $Q_1$ are integral, $P^\gp$ injects into $Q_1^\gp$. But monomorphisms in~$\catAb$ are stable under amalgamated sums (see dual statement of \cite[VIII.4 proposition 2]{160}).
\end{proof}

\begin{num}
We introduce the following notation. For a point $t$ of a fan $F$, we set $r(t) \isDef \dim M_{F,t}$, which is the height of the point~$t$.
\end{num}

\begin{prop}\label{FibresOnUt}
Let $\mcl{X}$ be a log regular scheme over~$\logs{S}$ and $\vphi\colon F'\to F$ a subdivision of its fan with $F'$ locally fs.
Let $t\in F'_s$ and $d$ be a divisor of its root index~$m$. We set $\vphi(d) = \gcd(\vphi(m),d)$, where $\vphi(m)$ is the root index of~$\vphi(t)$.
Then for every $x\in \mcl{X}(\vphi(d))^\fs$ with $\vphi(t) = \pi(x)$ we have
\[U_d(t)\cap (\vphi^*\mcl{X}(d)^\fs)_x = \mbb{G}_{m,k(x)}^{r(\vphi(t))-r(t)},\]
where $(\vphi^*\mcl{X}(d)^\fs)_x$ denotes the fibre over $x$ of the induced morphism $\vphi^*\mcl{X}(d)^\fs\to \mcl{X}(\vphi(d))^\fs$.
\end{prop}
\begin{proof}
Let $\ovl{x}$ be the image of $x$ in $\mcl{X}$. By replacing $\mcl{X}$ by an open neihbourhood of $\ovl{x}$, we may assume that it admits an fs $\mbb{N}$-chart $P\to \mcl{M}_\mcl{X}(\mcl{X})$ such that $F\isoto \Spec^\sharp P$ (see \ref{chartForFan}).
Denote by $Q$ the fibred product of monoids
\[M_{F',t}\times_{M_{F',t}^\gp} P^\gp.\]

Then by \ref{isoMonoidAlgebras}, \ref{calculUdt} and the construction in \ref{prelim:inverseSubdiv} we have the following fibred product
\begin{center}
 \usetikzlibrary{matrix}
\begin{tikzpicture}
    \matrix (m) [matrix of math nodes, row sep=3em,
    column sep=2em, text height=1.5ex, text depth=0.25ex]
    { U_d(t) &  & U_{\vphi(d)}(\vphi(t)) \\
\mbb{Z}\bigl[Q(d)^{\sat,\times}\bigr] & \mbb{Z}\bigl[P(d)^{\sat,\times}\bigr] & \mbb{Z}\bigl[P(\vphi(d))^{\sat,\times}\bigr] \\};

\path[arrow] (m-1-1) edge node[auto] {} (m-1-3)
edge node[auto,swap] {} (m-2-1)
(m-2-1) edge node[auto,swap] {} (m-2-2)
    (m-1-3) edge node[auto] {} (m-2-3)
(m-2-2) edge node[auto,swap] {} (m-2-3);

\draw ($(m-1-1.east)!.1!(m-1-3.west)$)+(0,-.12) coordinate (A);
\draw ($(m-1-1.south)!.3!(m-2-1.north)$)+(.12,0) -| (A);
\end{tikzpicture}\end{center}
Since the root index of $\mbb{N}_{\vphi(d)}\to P(\vphi(d))^\sat$ is $\vphi(m)/\vphi(d)$ (\ref{rootPdsat}), we have by \ref{PtoQsatSharp}
\[U_{\vphi(d)}(\vphi(t))\times_{\mbb{Z}[P(\vphi(d))^{\sat,\times}]} \mbb{Z}[P(d)^{\sat,\times}] \isoto U_{\vphi(d)}(\vphi(t)),\]
and taking the fibre over $x$ we get
\begin{center}
\begin{tikzpicture}
    \matrix (m) [matrix of math nodes, row sep=3em,
    column sep=3em, text height=1.5ex, text depth=0.25ex]
    {U_d(t)\cap (\vphi^*\mcl{X}(d)^\fs)_x & k(x)\\
 k(x)[Q(d)^{\sat,\times}] & k(x)[P(d)^{\sat,\times}]\\};
    \path[arrow]
    (m-1-1) edge node[auto] {$  $} (m-1-2)
            edge node[auto,swap] {$$}  (m-2-1)
    (m-1-2) edge node[auto] {$$}(m-2-2)
(m-2-1) edge node[auto] {$ $} (m-2-2);
\draw ($(m-1-1.center)!.3!(m-1-2.west)$)+(0,-.25) coordinate (A);
\draw ($(m-1-1.south)!.3!(m-2-1.north)$)+(.12,0) -| (A);
\end{tikzpicture}\end{center}

Since $P$ injects into $Q$ and $\mbb{N}\to\mbb{N}_d$ is integral, $P(d)^\sat$ injects into $Q(d)^\sat$ by \ref{integralPushout}. Furthermore, $Q(d)^{\sat,\times}/P(d)^{\sat,\times}$ is a subgroup of the free group $P(d)^\gp/P(d)^{\sat,\times}\isoto (P(d)^{\sat,\sharp})^\gp$, so that
\[Q(d)^{\sat,\times} \isoto P(d)^{\sat,\times}\oplus Q(d)^{\sat,\times}/P(d)^{\sat,\times}.\]
Furthermore,
\begin{align*}
r(\vphi(t)) - r(t) &= \dim M_{F,\pi(x)} - \dim M_{F',t}\\
&= \rk P^\gp - \rk P^\times - \rk Q^\gp + \rk Q^\times\\
&= \rk Q(d)^{\sat,\times} - P(d)^{\sat,\times}\\
&= \rk (Q(d)^{\sat,\times}/P(d)^{\sat,\times}),
\end{align*}
where we used \ref{PtoQsatSharp}, \ref{facesQsat} and the fact that $P^\gp = Q^\gp$ (see \ref{prelim:inverseSubdiv}).
We conclude that
\[U_d(t)\cap (\vphi^*\mcl{X}(d)^\fs)_x = \mbb{G}_{m,k(x)}^{r(\vphi(t))-r(t)}.\qedhere\]
\end{proof}

\begin{cor}\label{UxBlowUp}
Let $\mcl{X}\to\logs{S}$ be a log regular scheme over~$\logs{S}$, and $F$ its fan. Let $\vphi\colon F'\to F$ be a subdivision of $F$ with $F'$ locally fs and consider the induced morphism $\vphi^*\mcl{X}\to\mcl{X}$.
Then for every $t\in F'_s$, we have
\[ (\mbb{L}-1)^{r(t)-1}\, [\wtl{U}(t)] = (\mbb{L}-1)^{r(\vphi(t))-1}\, [\wtl{U}(\vphi(t))] \in K_0(\mathrm{Var}_{\mcl{X}_k}),\]
where brackets denote classes in the Grothendieck ring of $\mcl{X}_k$-varieties $K_0(\mathrm{Var}_{\mcl{X}_k})$.
\end{cor}
\begin{proof}
By \ref{FibresOnUt} applied to $d=m(t)$ and \cite[4.2.3]{48}, we see that $\wtl{U}(t)\to \wtl{U}(\vphi(t))$ is a piecewise trivial fibration of fibre $\mbb{G}_{m,k}^{r(\vphi(t))-r(t)}$.
The equality follows from \cite[4.2.2]{48}.
\end{proof}

\begin{rmq}
Applying proposition \ref{FibresOnUt} with $d=1$ instead in the proof of corollary \ref{UxBlowUp} yields the same result with each $\wtl{U}$ replaced by~$U$.
\end{rmq}

\begin{rmq}
We can more generally define $U_d(t)$ as the inverse image of $U(t)$ in the special fibre $(\mcl{X}(d)^\fs)_k$ for any integer $d\geq 1$.
Then it follows from the proposition below and the same arguments as in \ref{calculUdt} that if $d$ is a multiple of~$m$, $U_d(t)\isoto \wtl{U}(t)$ (compare with \cite[4.4]{106}).
\end{rmq}

\begin{prop}
Let $\mbb{N}\to P$ be a morphism between fs monoids of root index $1$.
Then $P(d)$ is saturated.
In particular, the maximal ideal of $P(d)^\sat = P(d)$ is generated by $P^+$ and $(0,1)$.
\end{prop}
\begin{proof}
Following the same argument as in~\ref{rootPdsat}, we can assume that $P$ is sharp with no loss of generality.

Let $(e_i)_{0\leq i \leq n}$ be a basis of $P^\gp$ with $e_0 = e_\pi$.
Then $P^\gp\isoto \mbb{Z}e_\pi\oplus \mbb{Z}^n$ and $P(d)^\gp\isoto \mbb{Z}\frac{e_\pi}{d}\oplus \mbb{Z}^n$.
The claim follows at once from the fact that the isomorphism
\[\mbb{Z}\tfrac{e_\pi}{d}\oplus \mbb{Z}^n\to \mbb{Z}e_\pi\oplus \mbb{Z}^n, \quad \tfrac{e_\pi}{d}\mapsto e_\pi\]
restricts to an isomorphism $P(d)\to P$.
The last assertion comes from~\ref{unitsInPd}.
\end{proof}


\begin{num}
If $\mcl{X}\to\logs{S}$ is log smooth, then $\mcl{X}$ is log regular by \ref{prelim:logLisseSurLogReg} and flat by \cite[4.5]{139}, since it is integral by \ref{NmorphismIntegral} and the proposition below.
Furthermore, if $\mcl{X}$ has pure relative dimension $m$, then its sheaf of log differentials $\Omega^\logr_{\mcl{X}/\logs{S}}$ is locally free of rank $m$ by \ref{prelim:OmegaLogLocFree}.
\end{num}

\begin{prop}\label{log:intgChartIntgMorph}
Let $u\colon P\to Q$ be a chart for a morphism of log schemes $f\colon X\to Y$ with $P$ and~$Q$ integral.
If $u$ is integral, then $f$ is integral.
\end{prop}
\begin{proof}
Let $x\in X$ and consider the diagram
\begin{center}
\begin{tikzpicture}
    \matrix (m) [matrix of math nodes, row sep=3em,
    column sep=3em, text height=1.5ex, text depth=0.25ex]
    { P & \mcl{M}_Y(Y) \\
Q & \mcl{M}_X(X)\\};
    \path[arrow]
    (m-1-1) edge node[auto] {$\beta_P $} (m-1-2)
            edge node[auto,swap] {$u$}  (m-2-1)
    (m-1-2) edge node[auto] {$$}(m-2-2)
(m-2-1) edge node[auto] {$\beta_Q$} (m-2-2);

\end{tikzpicture}
\end{center}
By \cite[I.4.5.3(2)]{137.2014}, it suffices to show that $\mcl{M}^\sharp_{Y,f(x)}\to \mcl{M}^\sharp_{X,x}$ is integral.
But \eqref{eq:PaSharp} shows that this morphism is canonically isomorphic to
\[P/\beta_P^{-1}(\mcl{O}^\times_{Y,f(x)})\to Q/\beta_Q^{-1}(\mcl{O}_{X,x}^\times).\]
Since $\beta_Q^{-1}(\mcl{O}_{X,x}^\times)$ is a face of~$Q$ whose inverse image in~$P$ equals $\beta_P^{-1}(\mcl{O}^\times_{Y,f(x)})$ by localness of $\mcl{M}_{Y,f(x)}\to \mcl{M}_{X,x}$, it is integral by \cite[I.4.5.3(3)]{137.2014}.
\end{proof}

\begin{num}
Let $\mcl{X}\to \logs{S}$ be a generically smooth log smooth scheme over $\logs{S}$ of pure relative dimension~$m$ and $\omega$ a volume form, i.e.\ a nowhere vanishing differential form of degree~$m$ on the generic fibre~$\mcl{X}_K$.
In particular it can be considered as an $m$-log differential over the log scheme $\mcl{X}\times_{\logs{S}}K$.
The invertible sheaf~$\Omega^{\logr,m}_{\mcl{X}/\logs{S}}$ together with the section $\omega\in \Omega^{\logr,m}_{\mcl{X}/\logs{S}}(\mcl{X}_K)$ induces a Cartier divisor $\divr(\omega)$ on~$\mcl{X}$.
Since $\omega$ is a volume form, the support of $\divr(\omega)$ is disjoint from the trivial locus of~$\mcl{X}$ and by \cite[11.8]{143}, there is a unique divisor $I_\omega$ of $F = F(\mcl{X})$ such that $I_\omega \mcl{O}_\mcl{X} = \divr(\omega)$ (a divisor of a fan is simply a sheafified version of a divisor of a monoid as seen in section~\ref{prelim_log:sec:valEtDiv}).
\end{num}

\begin{rmq}
It follow from \ref{prelim:logSmoothSmooth} that any log smooth $\logs{S}$-scheme with trivial log stucture on the generic fibre is generically smooth.
\end{rmq}

\begin{prop}\label{ordFeta}
For every $t\in F_\eta$ of height one, $v_t(I_\omega) = 1$.
\end{prop}
\begin{proof}
Let $V = \Spec A$ be an affine étale neighbourhood of $i(t)$ in $\mcl{X}_K$ on which we have a chart $\{1\}\to P$ for $\mcl{X}_K\to K$ inducing an étale morphism $h\colon V\to K[P]$.
Since $t$ is of height $1$, $\mcl{M}_{\mcl{X},i(t)}^\sharp\isoto \mbb{N}$ so that we may assume $P = \mbb{N}\oplus \mbb{Z}^r$.
Then $K[P] = K[x_0,x_1^{\pm1},\dots, x_r^{\pm 1}]$ and $(\Omega^{\logr,m}_{\mcl{X}/\logs{S}})_{i(t)}$ (resp.\ $\Omega^m_{\mcl{X},i(t)}$) is generated by $h^*(\dlog x_0\wedge\cdots\wedge\dlog x_r)$ (resp.\ $h^*(\dif{x_0}\wedge\cdots\wedge \dif{x_r})$).
Since $\omega$ generates $\Omega^m_{\mcl{X},i(t)}$, we deduce that $\omega\,(\Omega^{\logr,m}_{\mcl{X}/\logs{S}})_{i(t)} = h^*(x_0)\,(\Omega^{\logr,m}_{\mcl{X}/\logs{S}})_{i(t)}$ up to a unit, whence the proposition, because $h^*(x_0)$ generates~$\mfk{m}_{i(t)}$.
\end{proof}

\section{A formula for the volume Poincaré series}
\label{sec:formulaForPoincaré}

In this section, we set $R = k\llb \pi\rrb$, with $\car k =0$.
We keep writing~$K$ for the fraction field of~$R$ and~$S$ for its spectrum.
All $S$-schemes will be assumed separated and of finite type.
Let $\mcl{X}$ be a generically smooth $S$-scheme and $\omega$ a volume form on the generic fibre~$\mcl{X}_K$.
We denote by $\what{\mcl{X}}$ the $\pi$-adic completion of~$\mcl{X}$ and for every $d\geq 1$ we write $\omega(d)$ for the inverse image of~$\omega$ on the generic fibre of $\what{\mcl{X}}(d) \isDef \what{\mcl{X}}\times_R R(d)$.
The \emph{volume Poincaré series} of the pair $(\mcl{X},\omega)$ is defined as
\[S(\mcl{X},\omega;T) \isDef \sum_{d\geq 1} \Bigl(\int_{\what{\mcl{X}}(d)} \norme{\omega(d)}\Bigr)\, T^d \in \mcl{M}_{\mcl{X}_k}\llb T\rrb,\]
where $\mcl{M}_{\mcl{X}_k}$ is the localization $K_0(\mathrm{Var}_{\mcl{X}_k})[\mbb{L}^{-1}]$ of the Grothendieck ring of $\mcl{X}_k$-varieties and $\mbb{L}\isDef [\mbb{A}^1_{\mcl{X}_k}]$ (see~\ref{serre:déf:motInt} for a definition of the motivic integral).
We first recall a formula for $S(\mcl{X},\omega;T)$ when $\mcl{X}$ is sncd, see proposition~\ref{formulaSncd}.
The goal of this section is to extend this formula to log smooth $\logs{S}$-schemes; this is achieved in theorem~\ref{formula} and constitutes the main result of this chapter.

\subsection*{The normal crossings case}

\begin{num}\label{log:setUpSnc}
Let $\mcl{X}$ be an sncd $S$-scheme and resume the notations of \ref{prelim:logSncd}.
Let $J\sset I$ and $m_J = \gcd(N_j)_{j \in J}$.
We define a Galois cover $\wtl{E_J^\circ}$ of $E_J^\circ$ by gluing the étale covers
\[\Spec A[T]\big/\bigl(uT^{m_J}-1,(x_j)_J\bigr)\to \Spec A/(x_j)_{j\in J}.\]
\end{num}

The order of a volume form along a component of the special fibre of an sncd scheme was defined in \cite[6.8]{106}.
We recall the construction.
\begin{num}\label{ordEi}
Assume that $\mcl{X}$ has pure relative dimension~$m$ and consider a volume form~$\omega$ on~$\mcl{X}_K$.
Let $E_i$ be a component of the special fibre, $N_i$ its multiplicity and $\xi_i$ its generic point.
Then $\Omega_i = \Omega^m_{\mcl{X},\xi_i}/\torsion$ is a free rank 1 module over the discrete valuation ring $\mcl{O}_{\mcl{X},\xi_i}$.
Choose $a\geq 1$ such that $\pi^a\omega\in \Omega_\mcl{X}^m(\mcl{X})$. Then the \emph{order} of $\omega$ along  is defined to be
\[\ord_{E_i} \omega\isDef \ord_{E_i}(\pi^a\omega)- aN_i,\]
where $\ord_{E_i}(\pi^a\omega)$ is the length of the $\mcl{O}_{\mcl{X},\xi_i}$-module $\Omega_i/(\pi^a\omega\mcl{O}_{\mcl{X},\xi_i})$.
\end{num}

From \cite[7.6]{106} we have
\begin{prop}\label{formulaSncd}
Let $\mcl{X}$ be an sncd $S$-scheme and $\omega$ a volume form on~$\mcl{X}_K$.
Then
\begin{equation}\label{SXoT}
S(\mcl{X},\omega;T) = \mbb{L}^{-m} \sum_{\emptyset\neq J\sset I} (\mbb{L}-1)^{\norme{J}-1} [\wtl{E_J^\circ}] 
\sum_{k_j\geq 1, j\in J} \mbb{L}^{-\sum_J k_j\mu_j}\, T^{\sum_J k_j N_j}
\end{equation}
in $\mcl{M}_{\mcl{X}_k}\llb T\rrb$, where $\mu_j = \ord_{E_j} \omega$.
\end{prop}


The strategy of the proof of~\ref{formula} is to first interpret all the elements appearing in \eqref{SXoT} in terms of the log structure of~$\mcl{X}$, and then to show that it is stable under subdivisions of fans.
The result will subsequently follow from \ref{prelim:desingFan} and the proposition below.

\begin{prop}\label{log:equivReg}
Let $\mcl{X}$ be a log regular scheme over $\logs{S}$ and $F$ its fan. The following are equivalent:
\begin{enumerate}
\item For every $t\in F$, $M_{F,t}\isoto \mbb{N}^{r(t)}$,
\item the underlying scheme $\mcl{X}$ is regular,
\end{enumerate}
In this case, $\mcl{X}_k$ is an snc divisor and
\[\mcl{X}_k = \sum_{\substack{t\in F_s \\  r(t) = 1 }} v_t(e_\pi) \, E(t).\]
\end{prop}
\begin{proof}
The equivalence between 1. and 2. comes from \cite[7.5.35]{148}.
Assume now that they both hold and let $t\in F$ and $x\in U(t)$.
We can find an affine neighbourhood $V = \Spec A$ of~$x$ on which we have an $\mbb{N}$-chart $\phi\colon \mbb{N}^{r(t)} \oplus\mbb{Z}\to A$ local at $x$.
We denote by $x_i$ the image of the $i$th basis vector of $\mbb{N}^{r(t)}$ and by $u$ the unit $\phi(0,1)$.
Let $\mfk{p}_i$ be the prime ideal generated by $e_i$ and $v_i$ its valuation.
Those are exactly the height-one prime ideals of $\mbb{N}^{r(t)}$.
Hence we clearly have $\pi = u\prod_1^{r(t)} x_i^{v_i(e_\pi)}$.
The formula for $\mcl{X}_k$ follows from the fact that $v_i(e_\pi)\geq 1$ if and only if $\mfk{p}$ is a point of~$F_s$.
Finally, since $\mcl{M}_{\mcl{X},x}^+\mcl{O}_{\mcl{X},x} = (x_1,\dots,x_n)$, we conclude by \cite[4.2.15]{LIU} and log regularity that the $x_i$ are part of a regular system of parameters, hence $\mcl{X}_k$ has normal crossings.
\end{proof}



\begin{cor}\label{étaleCover}
Let $\mcl{X}$ be a log regular scheme over~$\logs{S}$ satisfying the equivalent conditions of \ref{log:equivReg} and denote by $F_s^1$ (resp. $F_\eta^1)$ the set of height-one points of~$F_s$ (resp. $F_\eta$).
Let $t\in F$ and set $J = \Spec M_{F,t}\cap F_s^1$ and $K = \Spec M_{F,t}\cap F_\eta^1$.
Then \[\wtl{U}(t)\isoto \wtl{E_{J}^\circ}\cap E_K^\circ,\]
where $E_K^\circ = \cap_{t\in K} E(t) \prive \bigcup_{t\in F_\eta^1\prive K} E(t)$.
\end{cor}
\begin{proof}
Let $x\in U(t)$ and $\phi\colon \mbb{N}^{r(t)} \oplus\mbb{Z}\to A$ be a chart such as in the proof of \ref{log:equivReg}.
In particular, $V\cap \bigcup_{t\in F_\eta^1\prive K} E(t) = \emptyset$.

Then, shrinking $V$ if necessary (cfr.\ \ref{chartForFan}) we have by \ref{calculUdt} and \ref{saturation}, with $m$ the root index of $t$,
\begin{align*}
\wtl{U}(t)\cap V(m)^\fs &= V\times_{\Spec \mbb{Z}[\mbb{N}^{r(t)}]} \Spec \mbb{Z}\bigl[(\mbb{N}^{r(t)}\oplus \mbb{Z})(m)^{\sat,\times}\bigr]\\
&= \Spec \frac{A[T]}{(T^{m}-u,(x_i)_{i\in J\cup K})}
\end{align*}
which gives the local expression of the étale cover $\wtl{E_{J}^\circ}\cap E_K^\circ$.
\end{proof}

\begin{lem}\label{ordEgaux}
Let $\mcl{X}$ be a log smooth $\logs{S}$-scheme of pure relative dimension~$m$ satisfying the equivalent conditions of \ref{log:equivReg}.
Let $t$ be a height-one point of $F(\mcl{X})_s$ and $E = E(t)$ the corresponding irreducible component of~$\mcl{X}_k$.
Then the canonical homomorphism $\Omega^m_{\mcl{X}, i(t)}\to \Omega^{\logr,m}_{\mcl{X}/\logs{S}, i(t)}$ induces an isomorphism
\[\Omega^m_{\mcl{X}, i(t)}/\torsion\isoto \Omega^{\logr,m}_{\mcl{X}/\logs{S}, i(t)}\]
and we have $v_t(I_\omega) = \ord_{E} \omega$.
\end{lem}
\begin{proof}
Let $x$ be a generator of the maximal ideal of $\mcl{O}_{\mcl{X},i(t)}$.
We write $\pi = u x^n$ with $u$ a unit of $\mcl{O}_{\mcl{X},i(t)}$. Since $\dlog\pi = 0$, $n\dlog x+\dlog u = 0$ and hence $\dlog x = -\dif{u}/nu$.
This implies that $\Omega^m_{\mcl{X}, i(t)}/\torsion\to \Omega^{\logr,m}_{\mcl{X}/\logs{S}, i(t)}$ is a surjective morphism between rank-one modules, hence an isomorphism.

Now let $a\geq 1$ be such that $\pi^a\omega\in \Omega_\mcl{X}^m(\mcl{X})$.
Then \cite[11.8]{143} ensures that
\[v_t(e_\pi^aI_\omega) = \ord_{E} \divr(\pi^a\omega).\]
We conclude by \ref{IsoDivM} that $v_t(I_\omega) = \ord_{E} \omega$.
\end{proof}

\subsection*{Proof of the formula}

\begin{num}
Let $\mcl{X}$ be a generically smooth log smooth scheme over $\logs{S}$ of pure relative dimension~$m$ and $F$ its fan. Let $\omega$ be a volume form on~$\mcl{X}_K$.
We set
\[S'(\mcl{X},\omega;T) \isDef \sum_{t\in F_s} (\mbb{L}-1)^{r(t)-1} [\wtl{U}(t)]
\sum_{u \in M_{F,t}^{\vee,\loc} } \mbb{L}^{-u(I_\omega)} \, T^{u(e_\pi)} \in \mcl{M}_{\mcl{X}_k}\llb T\rrb.\]
\end{num}

\begin{lem}\label{StableLogBlowUp}
Let $\mcl{X}$ be a generically smooth log smooth scheme over $\logs{S}$ of pure relative dimension~$m$, $\omega$ a volume form on $\mcl{X}_K$ and $h\colon \mcl{X}'\to\mcl{X}$ a morphism induced by a subdivision of fans $\vphi\colon F'\to F$ preserving the horizontal part. Then for every $t\in F_s$
\[\sum_{u \in M_{F,t}^{\vee,\loc}} \mbb{L}^{-u(I_\omega)} \, T^{u(e_\pi)} =
\sum_{s\in \vphi^{-1}(t)} \sum_{v \in M_{F',s}^{\vee,\loc}} \mbb{L}^{-v(I_{h^*\omega})} \, T^{v(e_\pi)} \in \mcl{M}_{\mcl{X}_k}\llb T\rrb.
\]
\end{lem}

\begin{proof}
Since $F'(\mbb{N}) = F(\mbb{N})$, for every $u\in M_{F,t}^{\vee,\loc}$, there is a unique $s\in \vphi^{-1}(t)$ and a unique $v\in M_{F',s}^{\vee,\loc}$ such that $v = u\circ \vphi_s$. In particular, $v(e_\pi) = u(e_\pi)$.

Moreover, since $\vphi$ preserves the horizontal part, $h^*\omega$ is a volume form on $\mcl{X}'_K\isoto \mcl{X}_K$ and since $h$ is log étale, $h^*\Omega^{\logr,m}_{\mcl{X}/\logs{S}} \isoto \Omega^{\logr,m}_{\mcl{X}'/\logs{S}}$ so that $\divr(h^*\omega) = h^*\divr(\omega)$.
Then $\Divr(\vphi)(I_\omega) = I_{h^*\omega}$ and $v(I_{h^*\omega}) = u(I_\omega)$.
\end{proof}

\begin{prop}\label{Sinvariant}
Let $\mcl{X}$ be a generically smooth log smooth scheme over~$\logs{S}$ of pure relative dimension~$m$ and $h\colon \mcl{X}'\to\mcl{X}$ a morphism induced by a subdivision of fans $\vphi\colon F'\to F$ preserving the horizontal part. Then
\[S'(\mcl{X}',h^*\omega;T) = S'(\mcl{X},\omega;T) \in \mcl{M}_{\mcl{X}_k}\llb T\rrb.\]
\end{prop}
\begin{proof}
Using in turn \ref{UxBlowUp} and \ref{StableLogBlowUp}, we get
\begin{align*}
S'(\mcl{X}',h^*\omega;T) &= \sum_{t\in F_s}\sum_{s\in \vphi^{-1}(t)} (\mbb{L}-1)^{r(s)-1} [\wtl{U}(s)]
\sum_{v \in M_{F',s}^{\vee,\loc} } \mbb{L}^{-v(I_{h^*\omega})} \, T^{v(e_\pi)}\\
&= \sum_{t\in F_s} (\mbb{L}-1)^{r(t)-1} [\wtl{U}(t)]
\sum_{s\in \vphi^{-1}(t)}
\sum_{v \in M_{F',s}^{\vee,\loc} } \mbb{L}^{-v(I_{h^*\omega})} \, T^{v(e_\pi)}\\
&= S'(\mcl{X},\omega;T).\qedhere
\end{align*}
\end{proof}



\begin{thm}\label{formula}
Let $\mcl{X}$ be a generically smooth log smooth scheme over $\logs{S}$ of pure relative dimension~$m$. Let $\omega$ be a volume form on $\mcl{X}_K$ and denote by $F$ the fan of $\mcl{X}$.
Then
\[S(\mcl{X},\omega;T) = \mbb{L}^{-m}\sum_{t\in F_s} (\mbb{L}-1)^{r(t)-1} [\wtl{U}(t)]
\sum_{u \in M_{F,t}^{\vee,\loc} } \mbb{L}^{-u(I_\omega)} \, T^{u(e_\pi)} \in \mcl{M}_{\mcl{X}_k}\llb T\rrb.\]
\end{thm}
\begin{proof}
We have to show that $S(\mcl{X},\omega;T) = \mbb{L}^{-m}S'(\mcl{X},\omega;T)$.
Since $\mcl{X}_K$ is smooth, $M_{F,t}\isoto\mbb{N}^{r(t)}$ for every $t\in F_\eta$ (\ref{log:equivReg}) and by \ref{prelim:desingFan} we can find a subdivision $F'\to F$ preserving the horizontal part and such that $M_{F',s}\isoto \mbb{N}^{r(s)}$ for every $s\in F'$.
Hence we may assume that $F=F'$ by \ref{Sinvariant} so that $\mcl{X}$ is sncd (cfr.\ \ref{log:equivReg}).

Let $t\in F_s$ and denote by $J_s = J_s(t)$ (resp. $J_\eta = J_\eta(t)$) the set of height-one points of~$F_s$ (resp.~$F_\eta$) belonging to $\Spec M_{F,t}$.
Then $M_{F,t} = \mbb{N}^{J_s} \oplus\mbb{N}^{J_\eta}$ and by \ref{ordFeta}
\begin{align*}
(\mbb{L}-1)^{\norme{J_\eta}}\!\!\!\! \sum_{u \in M_{F,t}^{\vee,\loc} } \!\!\!\!\mbb{L}^{-u(I_\omega)} \, T^{u(e_\pi)}
 &= (\mbb{L}-1)^{\norme{J_\eta}} \!\!\! \sum_{(k_j)\in\mbb{N}_{\geq 1}^{J_\eta}}\!\! \mbb{L}^{-\sum k_j} \!\!\!\!\!\!\! \sum_{u \in (\mbb{N}^{J_s})^{\vee,\loc} } \!\!\!\!\!\! \mbb{L}^{-u(I_\omega)} \, T^{u(e_\pi)}\\
&= \sum_{u \in (\mbb{N}^{J_s})^{\vee,\loc} } \!\!\!\!\!\mbb{L}^{-u(I_\omega)} \, T^{u(e_\pi)}.
\end{align*}

Let $I$ be the set of height-one primes of $F_s$, i.e.\ all irreducible components of~$\mcl{X}_k$.
By \ref{étaleCover}, we get
\begin{align*}
S'(\mcl{X},\omega;T) &= \sum_{t\in F_s} (\mbb{L}-1)^{\norme{J_s(t)}-1} [\wtl{U}(t)]
 \sum_{u \in (\mbb{N}^{J_s(t)})^{\vee,\loc} } \mbb{L}^{-u(I_\omega)} \, T^{u(e_\pi)}\\
 &= \sum_{\emptyset\neq J\sset I} (\mbb{L}-1)^{\norme{J}-1} [\wtl{E_J^\circ}]
 \sum_{u \in (\mbb{N}^{J})^{\vee,\loc} } \mbb{L}^{-u(I_\omega)} \, T^{u(e_\pi)}\\
 &= \sum_{\emptyset\neq J\sset I} (\mbb{L}-1)^{\norme{J}-1} [\wtl{E_J^\circ}]
 \sum_{(k_j) \in \mbb{N}_{\geq 1}^{J} } \mbb{L}^{-\sum k_j\mu_j} \, T^{\sum k_j N_j},
\end{align*}
where the $N_j$ are the multiplicities of the components of $\mcl{X}_k$ and $\mu_j$ the order of~$\omega$ along them (see \ref{ordEgaux}).
We consequently recover formula \eqref{SXoT}.
\end{proof}

\begin{rmq}\label{subdiviseEnNeron}
Let $\mcl{X}$ be an sncd $S$-scheme of pure relative dimension~$m$ and $\omega$ a volume form on $\mcl{X}_K$. Let $d$ be a positive integer.
By \cite[5.17]{107} there is a series of blow-ups of strata $h\colon \mcl{X}'\to \mcl{X}$ such that $d$ is not $\mcl{X}'_k$-linear in the sense of~\cite[2.1]{107}. 
This condition of nonlinearity allows to give an explicit expression for a Néron smoothening for $\mcl{X}'(d)\to S(d)$ (see \ref{serre:NéronSmoothXd}). It follows from \cite[6.13]{106} that the coefficient of $T^d$ in $S(\mcl{X},\omega;T)$ is 
\begin{equation}\label{coeffNotLinear}
\sum_{N'_i\mid d} \bigl[\wtl{E}_i^{\prime\circ}\bigr] \, \mbb{L}^{d/N'_i \ord_{E'_i}\omega}.
\end{equation}
If we endow $\mcl{X}$ with the log structure described in \ref{prelim:logSncd}, those blow-ups of strata correspond to log blow-ups. By \cite[7.6.18(ii)]{148} such a log blow-up is induced by a blow-up of the fan~$F$ of~$\mcl{X}$, which is actually a kind of subdivision.
So the morphism $h$ is induced by a subdivision of~$F$ and we deduce that the coefficient of~$T^d$ in $S'(\mcl{X},\omega;T)$ equals \eqref{coeffNotLinear}.
Hence we see that our proof actually doesn't rely on \cite[7.6]{106} but can be directly carried out from \cite[4.5]{106}.
\end{rmq}

\begin{cor}\label{log:poles}
Let $\mcl{X}$ be a generically smooth log smooth scheme over $\logs{S}$ of pure relative dimension~$m$, $\omega$ a volume form and $F$ its fan. Then every pole of the function $S(\mcl{X},\omega;\mbb{L}^{-s})$ is of the form $s = -v_t(I_\omega)/v_t(e_\pi)$, for some point $t\in F_s$ of height one.
\end{cor}

\begin{proof}
%
Let $t$ be a point in $F_s$ (not necessarily of height one) and $(v_i)_I$ be the valuations of $M_{F,t}$ so that the cone $\sigma(M_{F,t}^{\vee,\loc})$ is strictly positively spanned by the $v_i$ (see \ref{prelim:valuations}).
By  \cite[2.8]{161} we can partition $\sigma(M_{F,t}^{\vee,\loc})$ into cones $(\Delta_n)_{n\in N}$ such that each $\Delta_n$ is strictly positively spanned by some $(v_j)_{j\in J_n}$ linearly independent over~$\mbb{R}$, where $J_n$ is some subset of~$I$.
If we set
\[D_n^\loc \isDef \bigl\{{\textstyle \sum}_{J_n}\lambda_j v_{j}\bigm| 0\leq\lambda_j\leq1\bigr\} \cap \Delta_n \cap (M_{F,t}^\vee)^\gp,\]
then every element $u$ of $\Delta_n \cap (M_{F,t}^\vee)^\gp$ can be uniquely written as $u = u_0 + \sum_{J_n} k_j v_j$, with $k_j\in\mbb{N}$ and $u_0\in D_n^\loc$.
We compute
\begin{align*}
\sum_{u \in M_{F,t}^{\vee,\loc} } \!\!\!\!\mbb{L}^{-u(I_\omega)}\, T^{u(e_\pi)}
&= \sum_{n\in N}\sum_{u\in D_n^\loc} \!\!\mbb{L}^{-u(I_\omega)} T^{u(e_\pi)} \!\!\!\!\!\!\!\sum_{k_j\geq 0, j\in J_n} \!\!\!\!\!\mbb{L}^{-\!\sum_j \! k_j v_j(I_\omega)}\, T^{\sum_j\! k_j v_j(e_\pi)}\\
&= \sum_{n\in N}\sum_{u\in D_n^\loc} \mbb{L}^{-u(I_\omega)} T^{u(e_\pi)} \prod_{J_n} \sum_{k_j\geq 0} \Bigl(\mbb{L}^{-v_j(I_\omega)}\, T^{v_j(e_\pi)}\Bigr)^{k_j}\\
&= \sum_{n\in N} \frac{\sum_{u\in D_n^\loc}\mbb{L}^{-u(I_\omega)} T^{u(e_\pi)}}{\prod_{J_n}1-\mbb{L}^{-v_j(I_\omega)}T^{v_j(e_\pi)}}.
\end{align*}
Hence if $S(\mcl{X},\omega;\mbb{L}^{-s})$ has a pole in $s$, then $-v_j(I_\omega)-sv_j(e_\pi) = 0$ for some $j\in I$, i.e.\ $s = - v_j(I_\omega)/v_j(e_\pi)$.
Note that $v_j(e_\pi) = 0$ if and only if $v_j$ is the valuation associated to a point of $F_\eta$.
In this case $v_j(I_\omega) = 1$ by \ref{ordFeta} and the factors $(\mbb{L}-1)^{-1}$ appearing in the expression (which \emph{sensu stricto} don't belong to $\mcl{M}_{\mcl{X}_k}$) are cancelled by the coefficient $(\mbb{L}-1)^{r(t)-1}$, see also \cite[\S 10.2]{187}.
\end{proof}

\begin{num}
Let $X$ be a Calabi-Yau variety of dimension $m$ over $K = \Frac R$.
Let $\omega\in \Omega_X^m(X)$ be a volume form on~$X$. Then the \emph{motivic zeta function} of~$(X,\omega)$ is defined as
\[Z_{X,\omega}(T) \isDef \mbb{L}^m \,\sum_{d\geq 1} \Bigl(\int_{X(d)} \norme{\omega(d)}\Bigr)\, T^d \in \mcl{M}_k\llb T\rrb\]
(see \cite[6.4]{146} when $\omega$ is distinguished).
If $\mcl{X}$ is a proper $S$-model of $X$, then by \cite[7.2]{106} 
\[Z_{X,\omega}(T) = \mbb{L}^m\, S(\mcl{X},\omega;T)  \in \mcl{M}_k\llb T\rrb.\]
Hence corollary \ref{log:poles} gives us a set of candidate poles for the zeta function $Z_{X,\omega}(T)$ when $\mcl{X}$ is log smooth, and this set is much smaller than the one we would get from a desingularization of $\mcl{X}$.
Log smooth models of Calabi-Yau varieties over~$K$ appear naturally in the Gross-Siebert programme on mirror symmetry (see for example \cite{170}).
\end{num}

\begin{num}\label{computeMotivic}
Let $X$ be a smooth irreducible scheme of finite type over $k$ of dimension~$n$ and $f\colon X\to \mbb{A}^1_k = \Spec k[\pi]$ a dominant morphism.
Denote by $X_s$ the zero locus of~$f$ and set $X^* = X\prive X_s$.
Shrinking $X$ around $X_s$, we can assume that $f$ is smooth on $X^*$.
For any gauge form $\phi$ on~$X$, we can find a unique form $\alpha\in \Omega^{n-1}_{X^*/\mbb{A}^1_k}(X^*)$ such that $\alpha\wedge\dif{f}$ is the restriction of $\phi$ to $X^*$.
The induced volume form is called  the \emph{Gelfand-Leray form} and is denoted by $\tfrac{\phi}{\dif{f}}$.
By \cite[9.10]{106} the motivic zeta function of $f$ (as defined in \cite[3.2.1]{72}) can then be computed as 
\[Z_f(T) = \mbb{L}^{n-1}\, S\bigl(\mcl{X},\tfrac{\phi}{\dif{f}};\mbb{L}^{-1}T\bigr) \in \mcl{M}_{\mcl{X}_k}\llb T\rrb,\]
where $\mcl{X} = X\times_{\mbb{A}^1_k} S$.
Hence if $\mcl{Y}\to  \logs{S}$ is a log smooth $ \logs{S}$-scheme dominating~$\mcl{X}$ and such that $\what{\mcl{Y}}_K\isoto \what{\mcl{X}}_K$, \ref{formula} gives a formula for $Z_f(T)$ in terms of the model~$\mcl{Y}$.
We give an example of this in the next section.
\end{num}

\section{The motivic zeta function of a nondegenerate polynomial}
\label{log:sectionNewton}

Let $k$ be a field of characteristic $0$ and set $R = k\llb \pi\rrb$ and $K = k\llpar \pi \rrpar$.
For a tuple $\omega = (\omega_1,\dots,\omega_n)\in \mbb{N}^n$,  we denote by $x^\omega$ the product $\prod_{1\leq i\leq n} x_i^{\omega_i}$.

For the rest of this section, we fix a polynomial $f(x) = \sum_{\omega\in\Omega} a_\omega x^\omega$ in $k[x_1,\dots,x_n]$ with $f(0) = 0$, where $\Omega$ is a finite subset of $\mbb{N}^n$ and $a_\omega \neq 0$ for every $\omega\in\Omega$.
The polynomial $f$ determines a dominant morphism $\mbb{A}^n_k\to \mbb{A}^1_k$ and we set $\mcl{X} = \mbb{A}^n_k\times_{\mbb{A}^1_k} R$.
Note that the generic fibre $\mcl{X}_K$ is smooth since $\car k = 0$.

We will show in this section how to recover from~\ref{formula} a combinatorial formula for $Z_f(T)$ first proved by Guibert in \cite[2.1.3]{165}, when $f$ is a nondegenerate polynomial with respect to its Newton polyhedron (see definitions below).

\subsection*{Newton polyhedron and associated fan}

We refer to \cite[\S 2]{161} for details on Newton polyhedra.

\begin{déf}
The \emph{Newton polyhedron} $\Gamma_f$ of the polynomial~$f$ is the convex hull of $\bigcup_\Omega (\omega + \mbb{R}^n_{\geq0})$ in $\mbb{R}^n$.
\end{déf}

\begin{num}
We will assume throughout that $f$ is \emph{nondegenerate with respect to its Newton polyhedron}, which means that the varieties
\[\Spec k[x_1^{\pm 1},\dots,x_n^{\pm 1}]/(f_\tau)\] are smooth for every face $\tau$ of $\Gamma_f$, where $f_\tau\isDef \sum_{\omega\in\tau\cap\Omega} a_\omega x^\omega$.
\end{num}

\begin{num}
For every $u\in\mbb{R}^{n,\vee}_{\geq 0} = \Hom_\catMnd(\mbb{R}^n_{\geq0},\mbb{R})$ we set $m(u) \isDef \inf_{\alpha\in\Gamma_f} u(\alpha)$ and denote by $H_u$ the hyperplane of $\mbb{R}^n$ defined by $u(x) - m(u) = 0$.
Then
\[F(u) \isDef \{\alpha\in \Gamma_f\mid u(\alpha) = m(u)\} = H_u\cap\Gamma_f,\]
is a face of $\Gamma_f$ (compare with~\ref{prelim:déf:faceCone}). 
For a face $\tau$ of $\Gamma_f$ we set
\[\Delta_\tau \isDef \{u\in\mbb{R}^{n,\vee}_{\geq 0}\mid F(u) = \tau\}.\]
Its topological closure $\ovl{\Delta}_\tau$ is a strictly convex rational polyhedral cone of $\mbb{R}^{n,\vee}$ and we have (see \cite[2.6]{161})
\[\ovl{\Delta}_\tau = \{u\mid F(u)\sups \tau\}.\]
The collection $\Delta = \{\ovl{\Delta}_\tau\}_{\tau}$ forms a fan in the classical sense that is a subdivision of the fan generated by $\sigma(\mbb{N}^n) = \mbb{R}_{\geq 0}^n$ (see \cite[I.2, definition~3]{184}). \label{fanDelta}

We set
\begin{align*}
M_\tau &\isDef \ovl{\Delta}_\tau^\vee\cap \mbb{N}^n\\
&= \bigl\{\alpha\in \mbb{N}^n\bigm| v(\alpha)\geq 0\text{ for all } v\in \ovl{\Delta}_\tau\bigr\}\\
&= \bigl\{\alpha\in \mbb{N}^n\bigm| v(\alpha)\geq 0\text{ for all } v\in \Delta_\tau\bigr\}.
\end{align*}
The \emph{normal} of a facet (i.e.\ a face of codimension one) $\tau$ of $\Gamma_f$ is the unique vector $v$ of $\mbb{N}^{n,\vee}$ with coprime coordinates such that $\tau = F(v)$.
If $(\tau_i)_{1\leq i\leq r}$ are the facets containing the face $\tau$ and $(v_i)$ their normals, then
\[M_\tau = \bigl\{\alpha\in\mbb{N}^n\bigm| v_i(\alpha)\geq 0 \text{ for all } 1\leq i \leq r\bigr\}\]
(see also \ref{prelim:valuations}).
\end{num}

\begin{prop}\label{facetteMv0}
Let $\tau$ be a facet of $\Gamma_f$ and $v$ its normal. If $m(v) = 0$, then $v$ is an element of the canonical basis of~$\mbb{N}^{n,\vee}$.
\end{prop}
\begin{proof}
Without loss of generality, we may write $v = (a_1,\dots,a_r,0,\dots,0)$ for some $r\leq n$, with $a_i\geq 1$ for every $1\leq i \leq r$.
Let $\alpha = (\alpha_i)_i\in\mbb{N}^n$. We have $v(\alpha) = 0$ if and only if $\alpha_i = 0$ for every $1\leq i \leq r$ from which we get
\[\tau = \bigl\{\alpha\in \Gamma_f\bigm| v(\alpha) = m(v) = 0\bigr\} = \bigl\{\alpha\in \Gamma_f\bigm| \alpha_i = 0 \text{ for every } 1\leq i \leq r\bigr\}\]
which has dimension less or equal to $n-r$. Hence $r = 1$ because $\tau$ is a facet; this concludes the proof.
\end{proof}

\begin{prop}\label{dimMtau}
Let $\tau$ be a face of~$\Gamma_f$.
We have $\dim M_\tau = n - \dim \tau$ and $\rank (M_\tau)^\times = \dim \tau$.
\end{prop}
\begin{proof}
Since $\sigma(M_\tau) = \ovl{\Delta}_\tau^\vee$, we have $\dim M_\tau = \dim \ovl{\Delta}_\tau^\vee$ by \ref{prelim:mndCones}. The latter cone has dimension $n-\dim \tau$ by \cite[2.6]{161}.
The other equality comes from
\[\dim M_\tau = \rank (M_\tau)^\gp - \rank (M_\tau)^\times = n - \rank (M_\tau)^\times.\qedhere\]
\end{proof}

\begin{prop}\label{omega-omega0}
Let $\tau$ be a face of~$\Gamma_f$, $\omega_0\in\tau\cap\Omega$ and $\omega\in\Omega$. Then $\omega-\omega_0\in M_{\tau}$ and $\omega-\omega_0\in M_{\tau}^\times$ if and only if $\omega\in\tau$. 
\end{prop}
\begin{proof}
Let $(\tau_i)_{1\leq i\leq r}$ be the facets containing $\tau$ and $(v_i)_i$ the corresponding normals.
Since each $v_i$ reaches its minimum over $\Gamma_f$ on $\tau_i\sups\tau$, we have $v_i(\omega-\omega_0)\geq 0$ and $v_i(\omega-\omega_0) = 0$ if and only if $\omega\in\tau_i$.
In particular $\omega-\omega_0\in (M_\tau)^\times$ if and only if $\omega\in\tau_i$ for every $1\leq i \leq r$, i.e.\ $\omega\in\tau$ because $\tau = \bigcap_i \tau_i$.
\end{proof}

\subsection*{A log smooth model}

\begin{num}
The fan $\Delta$ obtained in \ref{fanDelta} induces a proper morphism $\phi\colon X_\Delta\to \mbb{A}^n_k$ as follows: for every face $\tau$ of $\Gamma_f$ we set
\[X_\tau \isDef \Spec k[M_\tau].\]
Since $\ovl{\Delta}_\tau\sset \mbb{R}^{n,\vee}_{\geq 0}$, the monoid~$M_\tau$ contains $\mbb{N}^n$, which yields a morphism $k[x_1,\dots,x_n]\to k[M_\tau]$.
Two varieties $X_\tau$ and $X_{\tau'}$ can be glued along $X_\sigma$, where $\sigma$ denotes the smallest face of~$\Gamma_f$ containing both $\tau$ and $\tau'$.
This construction yields a morphism $\phi\colon X_\Delta\to \mbb{A}^n_k$, which is proper by~\cite[1.15]{161}.
We also set $\mcl{X}_\Delta \isDef X_\Delta\times_{\mbb{A}^1_k} R$ where $X_\Delta\to \mbb{A}^1_k$ is the composition of $\phi$ with the morphism $\mbb{A}^n_k\to \mbb{A}^1_k$ determined by~$f$.
\end{num}

\begin{num}\label{logStructureXdelta}
We endow $X_\Delta$ with the log structure associated to the pre--log structure locally given by
\[M_\tau\oplus \mbb{N}\to k[M_\tau],\, (\alpha,n)\mapsto x^{\alpha}\,(f/\omega_0)^n,\]
where we choose $\omega_0\in\tau\cap\Omega$ and $f/\omega_0$ is defined as the element $\sum_\Omega a_\omega x^{\omega-\omega_0}$ of~$k[M_\tau]$.
By \ref{omega-omega0} we see that the associated log structure doesn't depend on the choice of $\omega_0$.
If we endow $\mbb{A}^1_k = \Spec k[\mbb{N}]$ with its natural log structure defined by the monoid $\mbb{N}$, then the morphism $\phi\circ f$ can be considered as a morphism of log schemes.
In the remainder, we will always implicitly use these log structures when referring to the \emph{log scheme}~$X_\Delta$ and to the morphism of \emph{log schemes} $\phi\circ f\colon X_\Delta\to\mbb{A}^1_k$.
It follows immediately from the definitions that $\mbb{N}\to M_\tau\oplus \mbb{N},\, 1\mapsto (\omega_0,1)$ is a chart for $X_\Delta\to\mbb{A}^1_k$.
\end{num}



\begin{prop}\label{XDeltaDagger}
The log scheme $X_\Delta$ is log smooth over $\mbb{A}^1_k$.
\end{prop}
\begin{proof}
Let $x\in X_\Delta$ and $\tau$ be the largest face of $\Gamma_f$ such that $x\in X_\tau$.
In particular $(M_\tau\oplus\mbb{N})^+\mcl{O}_{X_\Delta,x}\sset \mfk{m}_x\mcl{O}_{X_\Delta,x}$.
Let $r = \rank (M_\tau\oplus\mbb{N})^\times = \dim \tau$ (by~\ref{dimMtau}) and let us first assume that $r=0$.
In this case $\tau = \{\omega_0\}$ for some vertex $\omega_0\in\Omega$ and $\omega-\omega_0\in \mfk{m}_x \mcl{O}_{X_\Delta,x}$ for every $\omega\neq \omega_0$ by~\ref{omega-omega0}.
In particular,
\[u \isDef f/\omega_0 = a_{\omega_0} + \sum_{\omega\neq \omega_0} a_\omega x^{\omega-\omega_0}\]
 is a unit at $x$.
Hence $M_\tau\oplus\mbb{Z}$ is still a chart on the open $\Spec A$, with $A \isDef k[M_\tau][u^{-1}]$.
We get by (the proof of) \ref{SharpifiedChart} a chart $\mbb{N}\to M_\tau$ on the étale cover $\Spec A[u^{1/m}]$,
where $m$ denotes the root index of $\mbb{N}\to M_\tau\oplus\mbb{Z}$.
This chart induces a morphism
\[k[\mbb{N}]\otimes_{\mbb{Z}[\mbb{N}]}\mbb{Z}[M_{\tau}] = k[M_{\tau}] \to A[u^{1/m}]\]
which is étale, and we are done.

Assume now that $r>0$.
Since the root index of $\mbb{N}\to M_\tau\oplus \mbb{N}$ is one, we can find again by \ref{SharpifiedChart} a chart $\mbb{N} \to M_\tau^\sharp\oplus\mbb{N}$ for~$X_\tau\to\mbb{A}^1_k$.
The induced morphism is
\[k[M_\tau^\sharp\oplus\mbb{N}]\to \frac{k[M_\tau^\sharp\oplus\mbb{N}][M_\tau^\times]}{(f/\omega_0 - x^{(0,1)})} = k[M_\tau].\]
It is smooth if and only if
\[k[M_\tau^\sharp\oplus\mbb{N}]\to \frac{k[M_\tau^\sharp\oplus\mbb{N}][M_\tau^\times]}{(f_\tau/\omega_0)}\]
is smooth, because $\dif{(f/\omega_0-x^{(0,1)})} = \dif{(f_\tau/\omega_0)}$ over $k[M_\tau^\sharp\oplus\mbb{N}]$.
We conclude by base change to $k[M_\tau^\sharp\oplus\mbb{N}]$ since $k[\mbb{Z}^n]/(f_\tau/\omega_0)$ is smooth over~$k$ by nondegenerescence of $f$, hence also $k[M_\tau^\times]/(f_\tau/\omega_0)$.
\end{proof}

\begin{lem} \label{XDeltaIso}
The morphism $\phi\colon X_\Delta\to\mbb{A}^n_k$ induces an isomorphism
\[\phi^{-1}\bigl(\mbb{A}^n_k\prive V(f)\bigr)\toisoto \mbb{A}^n_k\prive V(f).\] 
\end{lem}
\begin{proof}
It follows from the general theory of toric varieties that $\phi$ is an isomorphism over all (varieties associated to) faces of $\mbb{N}^n$ that are not subdivided by $\Delta$. 
We will show that $f$ vanishes over all subdivided faces.

Let $\sigma$ be a face of $\mbb{N}^n$ which is subdivided by $\Delta$ and minimal with respect to this property. We may assume without loss of generality that $\sigma = \mbb{N}^r\oplus\{0\}^{n-r}$.
By assumption there is a facet $\tau$ of $\Gamma_f$ of normal $v$ belonging to the interior of $\sigma$.
Put differently, $v = (a_1,\dots,a_r,0,\dots,0)$ with $a_i\geq 1$ for every $1\leq i \leq r$.
Since faces of dimension one cannot be subdivided, $r \geq 2$ and $m(v)\neq 0$ by \ref{facetteMv0}.
This implies that every $\omega\in\Gamma_f$ has a nonzero $i$th component for some $1\leq i \leq r$.

Let $y\in \Spec k[\sigma^\vee\cap\mbb{N}^n]$. By the minimality assumption, we may assume that each coordinate function~$x_i$ belongs to~$\mfk{m}_y$ for $1\leq i \leq r$.
But then $x^\omega\in\mfk{m}_y$ for every $\omega\in \Gamma_f$ by the reasoning above, and $f\in \mfk{m}_y$, i.e.\ $f$ vanishes at~$y$.
\end{proof}

\begin{cor}\label{génériqueOK}
The $\logs{S}$-log scheme $\mcl{X}_\Delta$ endowed with the inverse log structure of $X_\Delta$ is a log smooth $\logs{S}$-model of~$\mcl{X}_K$ and $(\what{\mcl{X}}_\Delta)_K\isoto (\what{\mcl{X}})_K$.
\end{cor}
\begin{proof}
This is a direct consequence of \ref{XDeltaDagger}, \ref{XDeltaIso} and the following diagram of fibred products.
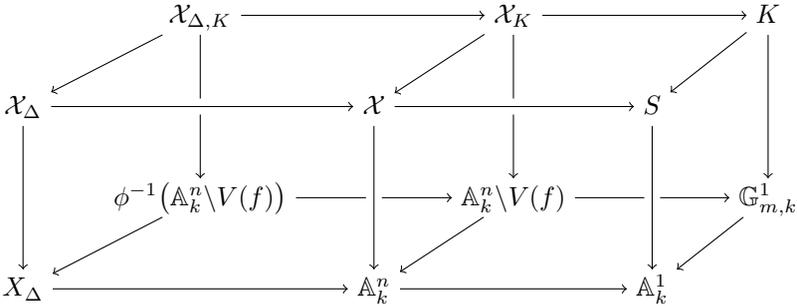
\begin{center}\begin{tikzpicture}[
cross line/.style={preaction={draw=white, -,
line width=6pt}}]
\matrix (m) [matrix of math nodes,
row sep=2em, column sep=2em,
text height=1.5ex,
text depth=0.25ex]{
& \mcl{X}_{\Delta,K} & & \mcl{X}_K & & K \\
\mcl{X}_\Delta & & \mcl{X} & & S\\
& \phi^{-1}\bigl( \mbb{A}_k^n\prive V(f)\bigr)& & \mbb{A}_k^n\prive V(f) & & \mbb{G}_{m,k}^1 \\
X_\Delta& & \mbb{A}_k^n & & \mbb{A}_k^1 \\};

\path[arrow]
(m-1-2) edge (m-1-4)
edge (m-2-1)
edge  (m-3-2)
(m-1-4) edge (m-3-4)
edge (m-2-3)
edge (m-1-6)
(m-2-1) edge [cross line] (m-2-3)
edge (m-4-1)
(m-3-2) edge (m-3-4)
edge (m-4-1)
(m-4-1) edge (m-4-3)
(m-3-4) edge (m-4-3)
edge (m-3-6)
(m-2-3) edge [cross line] (m-4-3)
edge [cross line] (m-2-5)
(m-1-6) edge (m-2-5)
edge (m-3-6)
(m-2-5) edge [cross line] (m-4-5)
(m-3-6) edge (m-4-5)
(m-4-3) edge (m-4-5);
\end{tikzpicture}\end{center}
The last part comes from the fact that the morphism $\phi\colon X_\Delta\to\mbb{A}^n_k$ is proper.
\end{proof}

This corollary allows us to compute $Z_f(T)$ from the log smooth $\logs{S}$-scheme~$\mcl{X}_\Delta$, following the construction sketched in~\ref{computeMotivic}.

\subsection*{Computing $Z_f(T)$}

\begin{num}
We will denote by $F_\Delta$ the fan associated to $\mcl{X}_\Delta$.
For a face $\tau$ of $\Gamma_f$ we also set
\[T(\tau) \isDef \bigl\{t\in F_\Delta\bigm| M_{F,t} = (M_\tau\oplus\mbb{N})^\sharp\bigr\},\]
which corresponds to $T(M_\tau^+\oplus\mbb{N}_{\geq 1})$ with the notation of~\ref{prelim:calculEt}.
Similarly
\[T(\tld{\tau}) \isDef \bigl\{t\in F_{\Delta,s}\bigm| M_{F,t} = (M_\tau\oplus\mbb{Z})^\sharp\bigr\} =T(M_{\tau}^+).\]
\end{num}

\begin{prop}
Let $\omega$ be a volume form on $\mbb{A}^n_k$ and denote by $I_\omega$ the divisor of $F_\Delta$ determined by the Gelfand-Leray form $\tfrac{\omega}{\dif{f}}$.
Let $\tau$ be a face of $\Gamma_f$ and $\omega_0\in \tau\cap \Omega$.
Then 
\[ I_\omega = \bigl(\frac{x_1\cdots x_n}{\omega_0}\bigr) \,M_{F,t}\] 
for every $t\in T(\tau)\cup T(\tld{\tau})$.
\end{prop}
\begin{proof}
We may assume that $\omega = \dif{x_1}\wedge\cdots\wedge\dif{x_n}$ and that $t$ is a point of height~one.
By \ref{prelim:OmegaLogCartesian} we may base our computations on $\Omega^\logr_{X_\Delta/\mbb{A}^1_k}$.
Assume first that $t$ corresponds to the height-one prime ideal generated by~$f/\omega_0$.
Then $M_\tau\sset \mcl{O}_{X_\Delta,t}^\times$ and in particular all $x_i$ are invertible at~$t$ so that we may assume by nondegenerescence of~$f$ that $\dif{f} = \sum_{i=1}^n a_i\dif{x_i}$ with $a_n\in\mcl{O}_{X_\Delta,t}^\times$.
In particular $(\dif{x_1},\dots,\dif{x_{n-1}},\dlog{f})$ is a basis for $\Omega_{X_\Delta/k,t}^\logr$ (cfr.\ \ref{prelim:computeOmegaFromChart}) and it follows from the locally split exact sequence
\begin{equation}\label{exctSeqOmega}
0\to f^*\Omega_{\mbb{A}^1_k/k,t}^\logr\to \Omega_{X_\Delta/k,t}^\logr\to \Omega_{X_\Delta/\mbb{A}^1_k,t}^\logr\to 0
\end{equation}
 of~\ref{prelim:OmegaLogExact} that $(\dif{x_1},\dots,\dif{x_{n-1}})$ is a basis for $\Omega_{X_\Delta/\mbb{A}^1,t}^\logr$.
Furthermore, we compute in~$\Omega^{\logr,n}_{X_\Delta/k,t}$:
\[ \dif{x_1}\wedge\cdots\wedge\dif{x_n} = a_n^{-1}\, \bigl(\dif{x_1}\wedge\cdots\wedge\dif{x_{n-1}}\bigr)\wedge \dif{f},\]
thus $a_n^{-1}\, \bigl(\dif{x_1}\wedge\cdots\wedge\dif{x_{n-1}}\bigr)$ is the Gelfand-Leray form on~$X_\Delta^*$ (notation of~\ref{computeMotivic}) and $\divr(\tfrac{\omega}{\dif{f}}) = \divr(1)$ because $\dif{x_1}\wedge\cdots\wedge\dif{x_{n-1}}$ is a basis for $\Omega^{\logr,n-1}_{X_\Delta/\mbb{A}^1,t}$.
This agrees with the claim since and $\frac{x_1\cdots x_n}{\omega_0}$ is a unit at~$t$.

If $t$ corresponds to a height-one prime ideal of $M_\tau$, then $f/\omega_0$ is invertible at~$t$, so that $(\dlog x_i)_{1\leq i \leq n}$ is a basis for $\Omega_{X_\Delta/\mbb{A}^1,t}^\logr$ because $M_\tau^\gp = \mbb{Z}^n$.
In particular $\dlog f = \sum_{i=1}^n a_i \dlog x_i$ for some $a_i \in \mcl{O}_{X_\Delta,t}$ and we may assume that $a_n$ is invertible since \eqref{exctSeqOmega} ensures that $\dlog f$ doesn't vanish at~$t$.
We compute in~$\Omega^{\logr,n}_{X_\Delta/k,t}\otimes \Frac(\mcl{O}_{X_\Delta,t})$:
\begin{align*}
\dif{x_1}\wedge\cdots\wedge\dif{x_n} &= x_1\cdots x_n\, (\dlog x_1\wedge \cdots \wedge \dlog x_n)\\
&= a_n^{-1}\, x_1\cdots x_n\, (\dlog x_1\wedge \cdots \wedge \dlog x_{n-1}\wedge \dlog f)\\
&= a_n^{-1}\, \bigl(\frac{x_1\cdots x_n}{f}\bigr) (\dlog x_1\wedge \cdots \wedge \dlog x_{n-1})\wedge \dif{f},
\end{align*}
from which we get $\divr(\tfrac{\omega}{\dif{f}}) = (\frac{x_1\cdots x_n}{f})$.
This yields the claimed result since $f/\omega_0$ is a unit at~$t$.
\end{proof}

\begin{prop}\label{Xtau0}
We have
\[\bigcup_{T(\tau)}\wtl{U}(t)\times_k \mbb{G}_m^{r(t)-1} \isoto \Spec \frac{k[\mbb{Z}^n]}{(f_\tau)}\defIs \mcl{X}_\tau(0).\]
\end{prop}
\begin{proof}
The root index of $\mbb{N}\to M_\tau\oplus\mbb{N}$ is clearly 1.
So by \ref{prelim:calculEt}
\[\bigcup_{T(\tau)}\wtl{U}(t) = \Spec \frac{k[M_\tau^\times]}{(f_\tau/\omega_0)}.\]
By \ref{dimMtau} we have $r(\tau) = \dim M_\tau\oplus\mbb{N} = n+1 - \dim \tau$ and $\rank M_\tau^\times = \dim \tau$.
Thus
\[\bigcup_{T(\tau)}\wtl{U}(t)\times_k \mbb{G}_m^{r(t)-1} = \Spec \frac{k[\mbb{Z}^{n}]}{(f_\tau/\omega_0)},\]
which completes the proof since $\omega_0$ is invertible in $\mbb{Z}^{n}$.
\end{proof}

\begin{prop}\label{Xtau1}
We have
\[\bigcup_{T(\tld{\tau})}\wtl{U}(t)\times_k \mbb{G}_m^{r(t)-1} \isoto \Spec \frac{k[\mbb{Z}^n]}{(f_\tau-1)} \defIs \mcl{X}_\tau(1).\]
\end{prop}
\begin{proof}
Let $m$ be the root index of $\mbb{N}\to M_\tau\oplus\mbb{Z}$ and $\alpha_0\in (M_\tau\oplus\mbb{Z})^\sharp = M_\tau^\sharp$ such that $\alpha_0^m = [e_\pi] = [\omega_0]$.
Choose now elements $(\alpha_1,\dots,\alpha_{r(t)-1})\in (M_\tau^\sharp)^\gp$ such that $(\alpha_0, (\alpha_i)_i)$ is a basis of $(M_\tau^\sharp)^\gp$.
This defines a section $s\colon (M_\tau^\sharp)^\gp\to M_\tau^\gp$ and yields an isomorphism
\[M_\tau\toisoto M_\tau^\sharp\oplus M_\tau^\times,\, \omega\mapsto \bigl([\omega], \omega- s([\omega])\bigr).\]
We have by \ref{calculUdt} and \ref{saturation}
\[\bigcup_{T(\tld{\tau})}\wtl{U}(t) = \Spec \frac{R[M_\tau^\sharp\oplus M_{\tau}^\times][T^{\pm 1}]}{(f - \pi,f/\omega_0 - T^m,(x^\alpha)_{\alpha\in M_\tau^\sharp})}.\]
Since $e_\pi\in M_\tau^+$ and $\pi = f = f/\omega_0\cdot x^{\omega_0}$, we necessarily have $\omega_0\in M_\tau^+$.
Thus 
\[ \bigcup_{T(\tld{\tau})}\wtl{U}(t) = \Spec \frac{k[M_{\tau}^\times][T^{\pm 1}]}{(f_\tau/\omega_0 - T^m)}.\]
We now look at the isomorphism
\[(M_\tau^\times\oplus \mbb{Z}) \oplus \mbb{Z}^{r(t)-1} \toisoto M_\tau^\times \oplus (M_\tau^\sharp)^\gp\]
that is the identity on $M_\tau^\times$, sends the $i$th basis vector of~$\mbb{Z}^{r(t)-1}$ onto $\alpha_i$ and sends the generator of~$\mbb{Z}$ (corresponding to the element $T$) onto $-\alpha_0$.
Then $T^m$ is sent onto $x^{-[\omega_0]}$ under the induced morphism of $k$-algebras, from which we deduce
\[\wtl{U}(t) \times_k \mbb{G}_m^{r(t)-1} \isoto \Spec \frac{k[M_\tau^\times \oplus (M_\tau^\sharp)^\gp]}{(f_\tau - 1)} \isoto \Spec \frac{k[\mbb{Z}^n]}{(f_\tau - 1)}.\qedhere\]
\end{proof}

\begin{prop}
We have 
\[Z_f(T) = \sum_\tau \Bigl([\mcl{X}_\tau(0)] \frac{\mbb{L}^{-1}T}{1-\mbb{L}^{-1}T} + [\mcl{X}_\tau(1)]\Bigr) \!\!\sum_{u\in\Delta_\tau\cap\mbb{N}^n} \mbb{L}^{-\sigma(u)} \, T^{m(u)},\]
where the sum is taken over all faces $\tau$ of the Newton polyhedron $\Gamma_f$ of~$f$.
\end{prop}
\begin{proof}
Observe that $(P^\sharp)^{\vee,\loc} = P^{\vee,\loc}$ for every monoid $P$. Hence
\[((M_\tau\oplus\mbb{Z})^\sharp)^{\vee,\loc} = \Delta_\tau\cap\mbb{N}^n\]
and similarly $((M_\tau\oplus\mbb{N})^\sharp)^{\vee,\loc} = (\Delta_\tau\cap\mbb{N}^n)\oplus \mbb{N}_{\geq1}$, as we see from \ref{prelim:valuations}.
For $u\in\Delta_\tau\cap\mbb{N}^n$ and $k\geq 1$ we have
\[(u,k)(e_\pi) = (u,k)(\omega_0,1) = m(u)+k\]
and
\begin{align*}
(u,k)(I_\omega) &= (u,k)(1,\dots,1,0) - (u,k)(\omega_0,0)\\
&= u(1,\dots,1) - m(u).
\end{align*}
We put $\sigma(u) \isDef u(1,\dots,1)$ and denote by $S_\tau(T)$ the sum $\sum_{u\in\Delta_\tau\cap\mbb{N}^n} \mbb{L}^{-\sigma(u)} \, T^{m(u)}$.
From \ref{Xtau0} and \ref{Xtau1} we have
\begin{align*}
\sum_{T(\tau)}\!(\mbb{L}-1)^{r(t)-1} [\wtl{U}(t)] \!\!\!\!\!\!\!\sum_{u \in M_{F,t}^{\vee,\loc} } \!\!\!\!\!\!\mbb{L}^{-u(I_\omega)} T^{u(e_\pi)}\!
&= [\mcl{X}_\tau(0)] \!\!\!\!\!\!\!\sum_{u \in \Delta_\tau\cap\mbb{N}^n}\!\sum_{k=1}^\infty \mbb{L}^{\!\!-\sigma(u) + m(u)} T^{m(u) + k}\\
&= [\mcl{X}_\tau(0)]\, S_\tau(\mbb{L}T) \frac{T}{1-T}
\end{align*}
and
\begin{align*}
\sum_{T(\tld{\tau})}(\mbb{L}-1)^{r(t)-1} [\wtl{U}(t)] \!\!\!\!\sum_{u \in M_{F,t}^{\vee,\loc} } \!\!\!\mbb{L}^{-u(I_\omega)} \, T^{u(e_\pi)}
&= [\mcl{X}_\tau(1)] \!\!\!\sum_{u \in \Delta_\tau\cap\mbb{N}^n} \!\!\!\!\mbb{L}^{-\sigma(u) + m(u)} \, T^{m(u)}\\
&= [\mcl{X}_\tau(1)]\, S_\tau(\mbb{L}T).
\end{align*}
From which we get by \ref{computeMotivic}
\begin{align*}
Z_f(T) &= \mbb{L}^{n-1}\, S\bigl(\mcl{X},\frac{\omega}{df};\mbb{L}^{-1}T\bigr)\\
&= \sum_\tau \Bigl([\mcl{X}_\tau(0)] \frac{\mbb{L}^{-1}T}{1-\mbb{L}^{-1}T} + [\mcl{X}_\tau(1)]\Bigr)\, S_\tau(T),
\end{align*}
since $\mcl{X}_K$ is of dimension $n-1$.
\end{proof}

\begin{cor}
The local motivic zeta function of~$f$ around $0$ is given by
\[Z_{f,0}(T) = \!\!\sum_{\tau \text{ compact}} \Bigl([\mcl{X}_\tau(0)] \frac{\mbb{L}^{-1}T}{1-\mbb{L}^{-1}T} + [\mcl{X}_\tau(1)]\Bigr) \!\!\sum_{u\in\Delta_\tau\cap\mbb{N}^n} \!\!\mbb{L}^{-\sigma(u)} \, T^{m(u)} \in \mcl{M}_{0}\llb T\rrb.\]
\end{cor}
\begin{proof}
The function $Z_{f,0}(T)$ is the image of $Z_f(T)$ under the base change morphism $i^*\colon \mcl{M}_{\mcl{X}_k}\to \mcl{M}_{0}$. 
We have by \ref{calculUdt}
\[i^*[\wtl{U}(t)] = \bigl[\wtl{U}(t)\times_{k[\mbb{N}^n]} k\langle\mbb{N}^n/\mbb{N}^{n,+}\rangle\bigr]
= \begin{cases}[\wtl{U}(t)] & \text{ if } \mbb{N}^n\to M_\tau \text{ is local}\\ 0 &\text{ otherwise}.\end{cases}\]

By \cite[2.3]{161} a face $\tau$ of $\Gamma_f$ is compact if and only if there is some $u\in\mbb{R}_{>0}^{n,\vee}$ such that $\tau = F(u)$, i.e.\ if there is some element~$u$ of~$\Delta_\tau$ with positive coordinates.
Since $\Delta_\tau\cap \mbb{N}^n = M_\tau^{\vee,\loc}$, we see that $\tau$ is compact if and only if there is some $u\in M_\tau^{\vee,\loc}$ such that $u(x) >0$ for every $x\in\mbb{N}^{n,+}$, which is equivalent to the morphism $\mbb{N}^n\to M_\tau$ being local.
This finishes the proof.
\end{proof}

\begin{rmq}
In \cite{165} Guibert forgot the condition $\ac(f)(\vphi) = 1$ while computing the measure of the spaces 
\[\bigl\{\vphi\in\mcl{L}(\mbb{A}^n)\bigm| \ord_t x(\vphi) = a, \ord_t f(\vphi) = m_{\Gamma(a)} + k\bigr\}\]
in lemma 2.1.1.
This explains the appearance of an extra torus in front of $[\mcl{X}_\tau(0)]$ in his proposition 2.1.3.
Note also that Guibert uses a normalization factor $\mbb{L}^{-n}$ in the definition of $Z_f(T)$.
Finally, recall that we should be more careful with the way to write down this formula as hinted at the end of the proof of~\ref{log:poles}. We refer to~\cite[10.5]{187} for more details.
\end{rmq}

\section{Polynomials in two variables}
\label{log:twoVariables}

Let $f$ be a nonconstant polynomial in two variables over a field~$k$ of characteristic zero and $Y$ be the minimal embedded resolution of $(\mbb{A}^2_k,f^{-1}(0))$.
We denote by~$Y_k$ the total transform of~$f^{-1}(0)$.
In \cite{179} Strauss showed in a special case that the components of $Y_k$ inducing poles of $Z_f(T)$ are precisely those of the strict transform and those of the exceptional divisor meeting at least three other components of~$Y_k$.
Building up on this result, Loeser proved the monodromy conjecture for the $p$-adic zeta function of~$f$ in~\cite{175}.
His proof generalizes to the motivic framework.

In \cite{178} Veys noticed that the noncontributing components of $Y_k$ could be contracted to obtain what is called the \emph{log canonical model} of the pair $(\mbb{A}^2_k,f^{-1}(0))$.
He then proceeded to show a formula for the topological zeta function computed on this model, yielding a conceptual explanation for Strauss's result.

We will show in this section that the log canonical model induces a log smooth $\logs{S}$-scheme dominating $\mbb{A}^2_k\times_{\mbb{A}^1_k} S$, so that our formula yields exactly the right poles for~$Z_f(T)$ (up to nodes, see case~3.\ below)

Let $h\colon Y\to Z$ be the contraction morphism to the log canonical model~$Z$ of~$(\mbb{A}^2_k,f^{-1}(0))$.
We denote by~$\mcl{Y}$ (resp.\ $\mcl{Z}$) the $S$-scheme $Y\times_{\mbb{A}^1_k} S$ (resp.\ $Z\times_{\mbb{A}^1_k} S$) and keep writing $h$ for the induced morphism $\mcl{Y}\to \mcl{Z}$.
Let $E_1,\dots, E_r$ be a maximal chain of exceptional curves contracted by~$h$. 
By \cite[2.5]{178} we have three possibilities:
\begin{enumerate}
\item $E_1$ and $E_r$ each meet exactly one other component of the special fibre~$\mcl{Y}_k$, denoted respectively by $E_0$ and $E_{r+1}$;
\item only $E_r$ meets exactly one other component of~$\mcl{Y}_k$, denoted by $E_{r+1}$;
\item $r=1$ and $E_1$ meets a component $E$ of $\mcl{Y}_k$ in two distinct points.
\end{enumerate}
We will first assume that we are in the first situation.
We endow $\mcl{Z}$ with the log structure induced by $\mcl{Z}_k = h(\mcl{Y}_k)$.
We only have to show that $\mcl{Z}$ is log smooth at $x = h(\bigcup_{i=1}^r E_i)$.
Following \cite[\S 3]{177}, we see that $\mcl{Z}$ is log regular at~$x$.
By \ref{SharpifiedChart} we can find an $\mbb{N}$-chart $P\to\mcl{O}(U)$ in an étale neighbourhood $U$ of~$x$ with $P$ sharp and fs.
In order to apply Kato's criterion, it suffices to show that the induced morphism
$R[P]/(f-\pi)\to \mcl{O}_{\mcl{Z},x}$ is smooth.
Consider the morphism $\phi\colon k[P]\to R[P]/(f-\pi)\to \mcl{O}_{\mcl{Z},x}$ and denote by $\mfk{m}$ (resp.\ $\mfk{n}$) the inverse image of $\mfk{m}_x$ in $k[P]$ (resp.\ $R[P]/(f-\pi)$).
We easily check that the $\mfk{m}$-adic completion of $k[P]$ is isomorphic to the $\mfk{n}$-adic completion of~$R[P]/(f-\pi)$.
Hence, we only have to show that $\phi$ is smooth. 
It is flat by \ref{prelim:logRegFlat} and since $\mfk{m} = k[P^+]$, the fibre of $\phi$ is $k\to \mcl{O}_{\mcl{Z},x}/(P^+\mcl{O}_{\mcl{Z},x})$ which is regular, hence smooth over~$k$, which finishes the proof.

To handle case~2., we choose a horizontal divisor $E_0$ meeting $E_1$ transversally at a point $x\in E_1\prive E_2$ (this is possible by \cite[21.9.11/12]{EGA-IV4}).
Then the log structure induced by $\mcl{Y}_k\cup E_0$ is still log smooth and we can proceed exactly as before.

Case~3.\ is a bit of a shame: it corresponds to the resolution of a node, and $\mcl{Z}$ is still log smooth in this case, but only if considered with the étale topology.
Unfortunately, fan theory is specific to log schemes on Zariski sites, so that we cannot apply~\ref{formula}.
Such an exceptional component actually induces the same pole as the component of the strict transform it meets.

\cleardoublepage


%






\backmatter

\includebibliography



\makebackcoverXII

\end{document}
